\documentclass[12pt]{article}

\catcode`\ =\active \def {\`e} \catcode`\ˆ =\active \def ˆ{\`a} \catcode`\ =\active \def {\`u}
\catcode`\Ž =\active \def Ž{\'e} \catcode`\ =\active \def {\c c} \def \parn{\par\noindent} 
\catcode`\ =\active \def {\^e} \catcode`\™ =\active \def ™{\^o} 
\catcode`\" =\active \def "{\^\i } \catcode`\ž =\active \def ž{\^u} 

\catcode`\: =\active \def :{\ifmmode\string:\else\unskip\kern 2pt\string: \ignorespaces\fi } 
\catcode`\; =\active \def ;{\ifmmode\string;\else\unskip\kern 2pt\string; \ignorespaces\fi } 

\font \db = msbm10 at 12 pt

\def \a{\alpha}  \def \b{\beta}   \def \g{\gamma}  \def \d{\delta}   \def \D{\Delta}  \def \G{\Gamma}
\def \la{\lambda} \def \L{\Lambda} \def \o{\omega} \def \s{\sigma}   \def \t{\theta}  
\def \e{\varepsilon}  \def \f{\varphi}    \def \rr{\varrho}  

\def \R{\mbox{\db R}}  \def \H{\mbox{\db H}} \def \E{\mbox{\db E}} \def \P{\mbox{\db P}} 
\def \N{\mbox{\db N}}

\def \S{\mbox{\db S}}      \def \Z{\mbox{\db Z}}

\def \AA{{\cal A}}    \def \CC{{\cal C}}       
\def \FF{{\cal F}}  \def \GG{{\cal G}}         
\def \LL{{\cal L}}  \def \M{{\cal M}}       \def \O{{\cal O}}

               \def \bv{\vec{b}}

\def \p{\partial}        \def \ii{\infty}  \def \sm{\setminus}  
\def \ra{\rightarrow}  \def \lra{\longrightarrow}        \def \Ra{\Rightarrow}  
   \def \LRa{\Leftrightarrow}    \def\LRA{\Longleftrightarrow}
              \def \nea{\nearrow}        \def \sea{\searrow}

\def \ub{\underbar}      \def \ol{\overline} 

   \def \ds{\displaystyle}   \def \ts{\textstyle}       
\def \rt1{\sqrt{-1}\,\,}  \def \1{^{-1}}            \def \2{^{-2}}             \def \5{{\ts {1\over 2}}}
\def \parn{\par\noindent}     \def\indf{\leavevmode\indent }

\def\cotg{{\rm cotg}\,}      \def\coth{{\rm coth}\,}        \def\th{{\rm th}\,}
\def\ch{{\rm ch}\,}       \def\sh{{\rm sh}\,}

\begin{document}

\newtheorem{defi}{Definition}
\newtheorem{theo}{Theorem}
\newtheorem{prop}{Proposition}  \newtheorem{propr}{Property}
\newtheorem{cor}{Corollary}
\newtheorem{lem}{Lemma}
\newtheorem{rem}{Remark}

\title{\bf Relativistic Diffusions \\  and Schwarzschild Geometry}

\author{Jacques FRANCHI \quad and \quad Yves LE JAN}
\date { March 2005 }
\maketitle
\vspace{-4mm} 
\begin{abstract}
The purpose of this article is to introduce and study a relativistic motion whose acceleration, in proper time, is given by a white noise. We deal with general relativity, and consider more closely the problem of the asymptotic behaviour of paths in the Schwarzschild geometry example. \par  
\end{abstract}
\bigskip 
\centerline{\small CONTENTS} \par 

   1) \ Introduction \par  
   \ref{sec.res}) \ Statement of the results \par 
   \quad  \ref{sec.relr}) \ A relativistic diffusion in Minkowski space \par 
   \quad \ref{sec.relg}) \ Extension to General Relativity \par 
   \quad \ref{sec.S}) \ The restricted Schwarzschild space ${\cal S}_0$ \par 
   \qquad  \ref{sec.syst}) \ The stochastic differential system in spherical coordinates \par 
   \qquad  \ref{sec.reduc}) \  Energy and angular momentum \par 
   \qquad  \ref{sec.comp.r}) \ Asymptotic behaviour of the relativistic diffusion \par 
   \quad \ref{sec.bh}) \ The full Schwarzschild space $\cal{S}$ \par  
   \qquad  \ref{sec.bhh}) \ Diffusion in the full space $\cal{S}$ : hitting the singularity \par 
   \qquad  \ref{sec.reg}) \ Regeneration through the singularity : entrance law after $D'$ \par 
   \qquad  \ref{sec.SDalls}) \ The relativistic diffusion for all positive proper times  \par 
   \qquad  \ref{sec.NumbC}) \ Capture of the diffusion by a neighbourhood of the hole \par 
   \ref{sec.P}) \ Proofs \par 
   \ref{sec.geod}) \ Appendix : \ Study of timelike and null geodesics \par 
   \ref{Ref}) \ References. \par\medskip 
   
\eject 

{\it \small  Corpora cum deorsum rectum per inane feruntur, ponderibus propriis incerto tempore ferme incertisque locis spatio depellere paulum, tantum quod momen mutatum dicere possis. Quod nisi declinare solerent, omnia deorsum, imbris uti guttae, caderent per inane profundum, nec foret offensus natus, nec plaga creata principiis : ita nil umquam natura creasset. \par 
   \hfill    Lucrecius   }

\section{Introduction} \indf 
   The classical theory of Brownian motion is not compatible with relativity, as it appears clearly from the fact that the heat flow propagates instantaneously to infinity. 
   A Lorentz invariant generalized Laplacian was defined by Dudley  (cf [D1]) on the tangent bundle 
of the Minkowski space, and it was shown that there is no other adequate definition than this one \footnote{Note however that some physical models of diffusion in a relativistic fluid are not Lorentz invariant since the frame of ``the fluid at rest'' plays a specific role: cf [D] and its references.}, as long as Lorentz invariance is assumed. An intuitive description of the associated diffusion (i.e. continuous Markov process) is that boosts are continuously applied in random directions of space. We show that this process is induced by a left invariant Brownian motion on the PoincarŽ group. The asymptotic behaviour of the paths of this process was studied (cf [D3]).   

   Considering the importance of heat kernels in Riemannian geometry and the extensive use that is made of their probabilistic repesentation via sample paths, it is somewhat surprising that Dudley's
first studies were not pursued and extended to the general context, namely to Lorentz manifolds.
It is indeed easy to check that the ``relativistic diffusion'' can be defined on any Lorentz manifold
 using a development, as done below. The infinitesimal generator is the generator of the geodesic flow perturbed by the vertical Laplacian. But such an extension would have little appeal, if some natural questions such as the asymptotic behaviour and the nature of harmonic functions could not be solved in some exemples of interest. \par 
   Here we provide a rather complete study of this question in the case of Schwarzschild and Kruskal-Szekeres manifolds, which are used in physics to represent ``black holes''. 
The specific interest of these manifolds comes from the vanishing of Ricci curvature, their symmetry,  
and the integrability of the geodesic flow.

   The picture that comes out in the Kruskal-Szekeres case appears quite remarkable, with paths confined in a neighborhood of the singularity, while their velocity increases, and an infinity of $SO_3\,$-invariant harmonic functions. \par 
   One difficulty of the study (and it might explain why Dudley had few followers) is that no explicit solution was found. The reason is that, even after reduction using the symmetries, the operator cannot involve less than three coordinates (even in Minkowski space), instead of one for the Laplacian on Riemann spaces of constant curvature. Estimations and comparison techniques of stochastic analysis are the main tools we use to prove our results. They do not include  yet a full determination of the Poisson boundary, but they suggest that for general Lorentz manifolds bounded harmonic functions could be characterised by classes of light rays, i.e. null geodesics. \par 

   Let us now explain more precisely the content of this article. \par  
   We consider diffusions, namely continuous strong Markov processes.  We start, in Section 
\ref{sec.relr}, with the flat case of Minkowski space $\R^{1,d}$, and therefore with the
Brownian motion of its unit pseudo-sphere, integrated then to yield the only true relativistic diffusion,
according to [D1]. We get then its asymptotic behaviour, somewhat simplifying the point of view of [D3]. \par 
   
   In Section \ref{sec.relg} below, we present an extension of the preceding construction to the framework of general relativity, that is to say of a generic Lorentz manifold. The process is first defined at the level of pseudo-orthonormal frames, with Brownian noise only in the vertical directions, and projects into a diffusion on the pseudo-unit tangent bundle. The infinitesimal generator we get in Theorem \ref{the.gen} decomposes into the sum of the vertical Laplacian and of the horizontal vector field generating the geodesic flow.  
   
   In Sections \ref{sec.S} and \ref{sec.bh}, we deal in detail with the Schwarzschild space, which is
the most classical example of curved Lorentz manifold, used in physics to model the space outside a black hole or a spherical body. \par 

   Using the symmetry, and introducing the energy and the angular momentum, which are constants of the geodesic motion, we reduce the problem to the study of a degenerate three-dimensional diffusion. We then establish in Theorem \ref{the.asymp} that almost surely the diffusion either hits the hole or wanders out to infinity, both events occurring with positive probability.
We prove also in Theorem \ref{the.asymp} that almost surely, conditionally on the non-hitting
of the hole, the relativistic diffusion goes away to infinity in some random asymptotic
direction, asymptotically with the velocity of light. Then we prove in Theorem \ref{the.bh} that almost
surely, conditionally on the hitting of the hole, the relativistic diffusion reaches the essential singularity at the center of the hole within a finite proper time, and we describe the limit.  

   We show then in Theorem \ref{the.alls} that the story can be continued further : namely, the 
Schwarzschild relativistic diffusion, a priori defined till its hitting of the center of the hole (the
essential singularity of the so-called Kruskal-Szekeres space, we also call full Schwarzschild space), can be extended to a diffusion which crosses this singularity. This extended diffusion reaches soon the restricted Schwarzschild space again, where it evolves as before, but maybe running backward in time. Thus such hole crossing can happen then again and again, but without accumulation, according to Theorem \ref{the.alls} below, so that the extended Schwarzschild relativistic diffusion is well defined for all positive proper times. \par 

   We finally study the asymptotic behaviour of this extended Schwarzschild relativistic diffusion, and show mainly in Theorem \ref{the.piege} below that there is a unique alternative possibility to the escape to infinity : there  is indeed a positive probability that the relativistic diffusion becomes endlessly confined in a spherical neighborhood of the hole, with an increasing velocity, and a trajectory becoming asymptotically planar, with an asymptotic shape. This implies the existence of an infinity of $SO_3\,$-invariant harmonic functions. 
      
\section{Statement of the results}    \label{sec.res} 

\subsection{A relativistic diffusion in Minkowski space} \label{sec.relr} \indf 
   Let us consider an integer $\,d\ge 2\,$ and the Minkowski space 
$\,\R^{1,d} := \{ \xi=(\xi^o,\vec{\xi}\,)\,\in\,\R\times\R^d\}$, endowed with the Minkowski pseudo-metric $\;\langle \xi,\xi\rangle := |\xi^o|^2-\|\vec{\xi}\|^2\,$. \par 
   Let $\,G\,$ denote the connected component of the identity in $\,O(1,d)\,$, and denote by 
$\;\H^d:=\{  \xi\in \R^{1,d}\,|\, \xi^o>0 \;\hbox{ and }\; \langle \xi,\xi\rangle =1\}$ the positive half of the unit pseudo-sphere. \par 
   The opposite of the Minkowski pseudo-metric induces a Riemannian metric on $\,\H^d$, namely the hyperbolic one, so that $\,\H^d\,$ is a model for the $d$-dimensional hyperbolic space. A convenient parametrization of $\,\H^d\,$ is $\;(\rr,\t)\in\R_+\times\S^{d-1}\,$, given by 
$\; \rr :={\rm argch} ( \xi^o)\,$ and $\;\t :=\vec{\xi}\,\Big/\sqrt{|\xi^o|^2-1}\,$. \ In these coordinates the 
hyperbolic metric writes $\; d\rr^2+\sh^2\rr\,|d\t|^2\,$, and the hyperbolic Laplacian is \quad 
${\ds \D^{\H} :={\p^2\over\p\rr^2}+(d-1)\coth\rr\,{\p\over\p\rr} +\sh^{-2}\rr\times\D_{\t}}\;$, \quad
$\D_{\t}\;$ denoting the Laplacian  of $\,\S^{d-1}\,$. \ The associated volume measure is 
$\;|\sh\rr |^{d-1} d\rr\, d\t\,$. \par 
   Note that $\,G\,$ acts isometrically on $\,\R^{1,d}\,$ and on $\,\H^d$, and that the Casimir operator on $\,G\,$ induces on $\,\H^d$ the hyperbolic Laplacian.  \par \medskip 

   Fix $\,\s >0\,$, and denote by $\,\LL^\s\,$ the $\,\s$-relativistic Laplacian, defined on  $\,\R^{1,d}\times\H^d\,$ by 
$$ \LL^\s f(\xi,p) \,:=\, p^o\,{\p f\over\p\xi^o}(\xi,p) + \sum_{j=1}^d\, p^j \,{\p f\over\p\xi^j}(\xi,p) 
+  {\ts{\s^2\over 2}}\,\D^{\H}_{(p)}f\,(\xi,p) \, ,  $$ 
that is to say \qquad ${\ds  \;\LL^\s f \,:=\, \langle p, {\rm grad}_{(\xi)}f\rangle + {\ts{\s^2\over 2}} \,\D^{\H}_{(p)}f }\,$.  \quad This is a hypoelliptic operator.  \par \medskip 

   Given any $\,(\xi_0,p_0)\in \R^{1,d}\times\H^d\,$, there exists a unique (in law) diffusion process 
\parn 
$(\xi_s,p_s)\,,s\in\R_+$, solving the $\,\LL^\s$-martingale problem, that is to say such that for any compactly supported $\,f\in C^2(\R^{1,d}\times\H^d)$, 
$\;{\ds f(\xi_s,p_s)- (\xi_0,p_0)-\int_0^s \LL^\s f(\xi_t,p_t)\, dt }\;$ is a martingale. \par 
   Note that $\,p_s\,$ is a hyperbolic Brownian motion, and that $\;{\ds \xi_s=\xi_0+\int_0^s p_t\,dt}\,$. \par 
   Note also that $\,\xi_s\,$ is parametrized by its arc length. 
Mechanically, $\,\xi_s\,$ describes the trajectory of a relativistic particle of small mass indexed by 
its proper time, submitted to a white noise acceleration (in proper time).  
Its law is invariant under any Lorentz transformation. \par 
    Note that if we denote by $\;(e_0,e_1,..,e_d)\;$ the canonical base of $\,\R^{1,d}$, and by $\,(e_j^*)\,$ the dual base (with respect to $\,\langle\;,\;\rangle$), the matrices $\, E_j:= e_0\otimes e_j^*+e_j\otimes e_0^*\,$ belong to the Lie algebra of $\,G\,$, and  generate the boost transformations.  Given $\,d\,$ independent real Wiener processes $\,w^j_s\,$, $\,p_s= (p^o_s,\vec{p}_s)\,$ can be defined by $\,p_s:= \L_se_0\,$, where the matrix $\,\L_s\in G\,$ is defined by the following stochastic differential equation : \par 
\centerline{${\ds \L_s = \L_0 + \s\sum_{j=1}^d \int_0^s \L_t\,E_j\,\circ dw^j_t\, .}$} \par \medskip 

    This means in fact that the relativistic diffusion process $(\xi_s,p_s)$ is the projection of some 
diffusion process having independent increments, namely a Brownian motion with drift, living on the PoincarŽ group. This group is the analogue in the present Lorentz-Minkowski setup of the classical group of rigid motions, and can be seen as the group of $(d+2,d+2)$ real matrices having the form $\, \pmatrix{ \L & \xi\cr 0 & 1\cr}$, with $\,\L\in G\,$, $\,\xi\in\R^{1,d}$ (written as a column), and $\,0\in\R^{1+d}$ (written as a row). Its Lie algebra is the set of matrices $\,\pmatrix{ \b & x\cr 0 & 0\cr}$, with $\,\b\in so(1,d)\,$ and $\,x\in\R^{1,d}$. The Brownian motion with drift we consider on the PoincarŽ group solves the stochastic differential equation $\; d\pmatrix{ \L_s & \xi_s\cr 0 & 1\cr} = \pmatrix{ \L_s & \xi_s\cr 0 & 1\cr}\circ d\pmatrix{ \b_s & e_0 s\cr 0 & 0\cr}$, where $\,(\b_s= \s\sum_{j=1}^d E_j\,w^j_s)\,$ is a Brownian motion on $\,so(1,d)$. This equation is equivalent to  $\; d\L_s=\L_s\circ d\b_s\,$ and $\,d\xi_s=\L_se_0\,ds\,$, so that $\,(\L_s)\,$ is a Brownian motion on $\,G\,$. On functions of $\,p=\L e_0\,$, its infinitesimal generator  $\,\sum_{j=1}^d(\LL_{E_j})^2\,$ coincides with a Casimir operator, and induces the hyperbolic Laplacian, so that $\,(p_s=\L_se_0)\,$ is a Brownian motion on $\,\H^d$, as required. 
\par\medskip 

   Then it is well known that $\;\t_s := \vec{p}_s\Big/\sqrt{|p^o_s|^2-1}\;$ converges almost surely in $\,\S^{d-1}\,$ to some random limit $\,\t_\ii\,$, and that $\,p^o_s\,$ increases to infinity. \ Set also $\; \rr_s :={\rm argch} ( p^o_s)\,$. \par \smallskip 

    The Euclidian trajectory $\,Z(t)\,$ is defined by $\,\vec{\xi}_{s(t)}\,$, where $\,s(t)\,$ is determined by $\,\xi^o_{s(t)} =t\,$. \parn 
Let us note that the Euclidian velocity ${\ds \,{d Z(t)\over dt} \,=\t_{s(t)}\,\th\rr_{s(t)}\,}$ has norm
$<1\,$ (1 is here the velocity of light). Moreover we have the following. 
\begin{rem}\label{pro.relr} \  The mean Euclidian velocity $\,{\ds{Z(t) \over t}}\,$ converges almost surely to $\,\t_\ii\in \S^{d-1}$. 
\end{rem} 
\ub{Proof}\quad  We have $\,\lim_{t\nea\ii}\limits\, s(t) = +\ii\,$, so  that $\,{\,\th\rr_{s(t)} = \sqrt{1-(p^o_{s(t)})\2}\,\,}$ goes to 1. Thus we get almost surely $\,{\ds \lim_{t\to\ii}\,{d Z(t)\over dt} \, = \t_\ii}\,$, and the result follows at once. $\;\diamond $ 

\begin{rem}\label{rem.relres}\quad {\rm The scattering amplitude, id est the law of $\,\t_\ii\,$ given $\,p_0\,$, is given by the hyperbolic harmonic measure in the unit ball of $\R^d$ (taken as model for $\H^d$), which has density proportional to $\,P(p_0,\cdot)^{d-1}$ with respect to the uniform measure of $\,\S^{d-1}$, $P$ denoting the classical Poisson kernel of the unit ball of $\R^d$. See for example ([E-F-LJ], case $\d=0$).  }
\end{rem}

\subsection{Extension to general Lorentz manifolds} \label{sec.relg} \indf 
    Let us now see how the preceding construction can be naturally extended to the framework of manifolds. \par 
    Let $\,\M\,$ be a $\,(d+1)$-dimensional manifold, equipped with a pseudo-Riemannian metric of signature $\,(+,-,..,-)\,$, together with an orientation and a time direction, and its Levi-Civita connection.  \par 
    For notational convenience, $\,T^1\M\,$ will always denote the positively oriented half of the unit tangent bundle of $\,\M\,$. As in the construction of Brownian motion on Riemannian manifolds, we have to use the frame bundle (see [M]). \par
    So let $\,G(\M)\,$ be the bundle of direct pseudo-orthonormal frames, with first element in the positive half of the unit pseudo sphere (in the tangent space), which has its fibers modelled on the special Lorentz group $\,G\,$. \ Let $\,V_j\,$ be the canonical vertical vector field associated with the preceding matrix $\,E_j\,$, and $\,H_0\,$ be the first canonical horizontal vector field. \par  
    Set \ $\;{\ds  \LL := H_0+{\ts{\s^2\over 2}}\sum_{j=1}^d V_j^2\,}$. \par 
    Let $\,\pi_1\,$ denote the canonical projection from $\,G(\M)\,$ onto the tangent bundle $\,T^1\M\,$, which to each frame associates its first element. The canonical vertical vector fields $\, V_{kl}\,$ associated with the matrices $\, E_{kl}:= e_k\otimes e_l^*-e_l\otimes e_k^*\,\in so(1,d)\,$, for $\,1\le k<l\le d\,$, generate an action of $\,SO_d\,$ on $\,G(\M)$, which leaves $\,T^1\M\,$ invariant and then allows the identification $\,T^1\M \equiv G(\M)/ SO_d\,$. The Casimir operator is $\;\CC = \sum_{j=1}^d\limits V_{j}^2 -\sum_{1\le k<l\le d}\limits V_{kl}^2\;$. \par  
    Note that the matrices $\,\{  E_{j}\,,\,E_{kl}\, ;\; 1\le j\le d\,,\, 1\le k<l\le d\}\,$ constitute a pseudo-orthonormal base of $\,so(1,d)\,$ (endowed with its Killing form). \par

\begin{lem}\label{lem.vert}\quad  The operators $\, H_0\;,\;\sum_{j=1}^d\limits V_{j}^2\;,\; \CC\;,\; \LL\;$ do act on $C^2$ functions on the pseudo-unit tangent bundle $\, T^1\M\,$, inducing respectively : the vector field $\,\LL_0\,$  generating the  geodesic flow on $\,T^1\M\,$, the so-called  
vertical Laplacian  $\,\D_v\,$,  $\,\D_v\,$ again, and the generator $\;\GG:= \LL_0 + {\ts{\s^2\over 2}}\,\D_v\,$. \par 
   More precisely, for any  test-function  $\,F\,$ on $\,T^{1}\M$, we have on $\,G(\M)$ :  
$$ (\LL_0F)\circ\pi_1\, =\, H_0(F\circ\pi_1)\; , \quad  (\D_vF)\circ\pi_1\, =\, \CC (F\circ\pi_1)\; . $$ 
Besides, in local coordinates $\;(x^i,e_j^k)\,$, with $\, e_j=e_j^k\,{\p\over\p x^k}\,$ : 
$\;V_j= e^k_j\,{\p\over\p e^k_0}+e^k_0\,{\p\over\p e^k_j}\;$, \ and 
($(g^{kl})\,$ denoting in these coordinates the inverse matrix of the  pseudo-Riemannian metric of $\,\M$) : 
$$ (\D_vF)\circ\pi_1 = \sum_{j=1}^dV_{j}^2 (F\circ\pi_1) = \Big( (e^k_0e^l_0-g^{kl})\,{\p^2\over\p e^k_0\p e^l_0} + d\,e^k_0\,{\p\over\p e^k_0}\Big) F \circ\pi_1\; . $$  
\end{lem}
\ub{Proof}\quad Let us observe that for any $\,u=(x,e_0,..,e_d)\in SO_{1,d}\M\,$, if $\,(u_s)\,$
denotes the horizontal curve such that $\,u_0=u\,$ and $\,\pi^*\dot{u}_0=\dot{x}_0=e_0\,$, then 
$\,(\pi_1 (u_s))\,$ is the geodesic generated by $\,\pi_1(u)=(x,e_0)\,$. Hence for any differentiable
function $\,F\,$ on $\,T^{1}\M\,$, we have 
$$ H_0(F\circ\pi_1)(u) = {d_o\over ds}F\circ\pi_1(u_s) ={d_o\over ds}F(\pi_1(u_s)) = \LL_0F(\pi_1
(u))\, . $$ 
Another way of expressing this is to recall that $\,H_0\,$ commutes with the rotation vertical vectors $\,V_{kl}\,$. It is also classical that the Casimir operator $\,\CC\,$ commutes with all vertical vectors $\,V_{kl}\,,V_j\,$.
\if{ whilst $\,[H_0,V_{j}]\not =0\,$ and $\,[H_0,\CC ]\not =0\,$. This explains why $\, H_0\,,\, \CC\,,\, \LL\,$ project on $\,T^1\M$, and $\, H_0\,,\, \LL\,$ do'nt project on $\,\M$. }\fi 
Moreover, since the rotation vectors $\,V_{kl}\,$ act trivially on $\,T^1\M$, the operators $\,\CC\,$ and $\,\sum_{j=1}^d\limits V_{j}^2\,$ induce the same operator $\,\D_v\,$ on $\,T^1\M$.  \par 
   Then as $\;e^k_0[e^{tV_j}(x,e)] = e^k_0(x,e)\,\ch t + e^k_j(x,e)\,\sh t\,$ and 
$\,e^k_j[e^{tV_j}(x,e)] = e^k_j(x,e)\,\ch t + e^k_0(x,e)\,\sh t\,$, we have indeed 
$\; V_jf(x,e) = {d_o\over dt}f[e^{tV_j}(x,e)] = e^k_j{\p\over\p e^k_0}f(x,e)+e^k_0{\p\over\p
e^k_j}f(x,e)\,$, id est $\; V_j= e^k_j{\p\over\p e^k_0}+e^k_0{\p\over\p e^k_j}\,$. \quad 
Using that $\; e^k_0 e^l_0 -\sum_{j=1}^d\limits e^k_je^l_j = g^{kl}\;$, we deduce immediately : 
$$ \sum_{j=1}^dV_{j}^2 = (e^k_0e^l_0-g^{kl})\,{\p^2\over\p e^k_0\p e^l_0}+\sum_{j=1}^d 
e^k_j\,{\p\over\p e^k_j}+ d\,e^k_0\,{\p\over\p e^k_0}+ 2\sum_{j=1}^d e^k_0e^l_j\,{\p^2\over\p e^l_0\p e^k_j}+\sum_{j=1}^d e^k_0e^l_0\,{\p^2\over\p e^l_j\p e^k_j}\; , $$ 
which reduces to the formula of the statement in the particular case of a function depending only on 
$\,(x,e_0)$. In accordance with the commutation relations arguments above, $\,\sum_{j=1}^dV_{j}^2 (F\circ\pi_1)\,$ is a function depending only on $\,(x,e_0)$, id est a function on $T^1\M$.  $\;\diamond$  
   
\medskip 

   Now, according to Section \ref{sec.relr}, the relativistic motion we will consider lives on $\,T^1\M\,$ and admits as infinitesimal generator the operator $\;\GG = \LL_0 + {\ts{\s^2\over 2}}\,\D_v\,$ of Lemma \ref{lem.vert} above. \par 
    If $\,\M\,$ is the Minkowski flat space of special relativity, it coincides with the diffusion defined in Section \ref{sec.relr} above.  \par 
    To construct this general relativistic diffusion, we use a kind of stochastic development to produce a stochastic flow on the bundle $\,G(\M)$, as is classically done to construct the Brownian motion on a Riemannian manifold. But we have now to project on $\,T^1\M\,$ and no longer on the base manifold $\,\M$, and to put the white noises on the acceleration, id est on the vertical vectors, and no longer on the velocity, id est on the horizontal vectors. \par 
       To proceed, let us simply fix $\,\Psi_0\in G(\M)\,$ and a $\,\R^d\,$-valued Brownian motion $\,w=(w^{j}_s)\,$, and let us consider the $\,G(\M)$-valued diffusion $\;\Psi=(\Psi_s)\in G(\M)\,$ solving the following Stratonovitch stochastic differential equation :   
$$ (*)\qquad \Psi_s = \Psi_0 + \int_0^s H_0(\Psi_t)\, dt + \s \int_0^s\, \sum_{j=1}^d\, V_{j} (\Psi_t)\circ dw^{j}_t \; .   $$ 

    By Lemma \ref{lem.vert}, the stochastic flow defined by $\,(*)\,$ commutes with the action of $\,SO_d\,$ on $\,G(\M)$, and therefore the projection $\;(\xi_s,\dot\xi_s) := (\xi_s,e_0(s)) = \pi_1(\Psi_s)\;$ defines a diffusion on $\,T^1\M$ ; namely this is the relativistic diffusion we intended to define and construct. \par  
   The following theorem defines the relativistic diffusion $(\xi_s,\dot\xi_s)$, possibly till some explosion time.  The vector field $\,\LL_0\,$ denotes the generator of the geodesic flow, which operates on the position $\,\xi$-component, and $\,\D_v\,$ denotes the vertical Laplacian (restriction to $\,T^1\M\,$ of  the Casimir operator on $\,G(\M)$), which operates on the velocity $\,\dot\xi$-component. 
\begin{theo}\label{the.gen}\quad  1) \ The $\,G(\M)$-valued Stratonovitch stochastic differential
equation 
$$ (*)\qquad d\Psi_s = H_0(\Psi_s)\, ds + \s\sum_{j=1}^dV_{j}(\Psi_s)\circ dw^{j}_s  $$ 
defines a diffusion $\,(\xi_s,\dot\xi_s) := \pi_1(\Psi_s)$ on $T^1\M$, whose infinitesimal generator
is $\,{\ds\LL_0+{\ts{\s^2\over 2}}\,\D_v}\,$.\parn  
2) \   If $\;\overleftarrow{\xi}(s) : T_{\xi_s}\M\rightarrow T_{\xi_0}\M\;$ denotes the inverse parallel transport along the $C^1$ curve $\,(\xi_{s'}\,|\,0\le s'\le s)$, then $\;\zeta_s := \overleftarrow{\xi}(s)\,\dot\xi_s\,$ is an hyperbolic Brownian motion on $\,T_{\xi_0}\M\,$. \parn
Therefore the path $(\xi_s)$ is the development of a relativistic diffusion path in the Minkowski space $\,T_{\xi_0}\M\,$. 
\end{theo}

\begin{rem}\label{rem.loc}\quad {\rm 
   In local coordinates $\;(x^i,e_j^k)\,$, with $\, e_j=e_j^k{\p\over\p x^k}\,$, $\;\Psi_s=\Big(\xi_s\,;\, e_0(s),..,e_d(s)\Big)\,$, $\,\G_{il}^k\,$ denoting as usual the Christoffel coefficients, the equation $\,(*)\,$ writes : 
\vspace{-2mm}
$$ dx^i_s = e^i_0(s) ds \; ; \quad de^k_j(s)= -\G_{il}^k(\xi_s) e^l_j(s) dx^i_s + 1_{\{ j\not= 0\}}\s 
e^k_0(s)\circ dw^{j}_s + 1_{\{ j=0\}}\s\sum_{i=1}^d e^k_i(s)\circ dw^{i}_s \; , $$ 
or equivalently in the It\^o form :
\vspace{-2mm}
$$ dx^k_s = e^k_0(s)\, ds \; ; \quad de^k_0(s)= -\G_{il}^k(\xi_s)\, e^l_0(s)\, dx^i_s + \s\sum_{i=1}^d e^k_i(s)\, dw^{i}_s + {\ts{d\,\s^2\over 2}}\, e^k_0(s)\, ds \; , \quad \mbox{ and} $$ 
\vspace{-2mm}
$$ de^k_j(s)= -\G_{il}^k(\xi_s)\, e^l_j(s)\, dx^i_s + \s\, e^k_0(s)\, dw^{j}_s + {\ts{\s^2\over 2}}\, e^k_j(s)\,ds
\quad \mbox{ for }\, j\ge 1\,,\; 0\le k\le d\; . $$ 

     Note that the martingales in the above equations for $\,e_0(s)$, that is to say the
differentials $\,{\ds dM_s^k := \s\sum_{j=1}^d e^k_j(s) dw^{j}_s}\,$, $\, 0\le k\le d\,$, have the
following quadratic covariation matrix : 
$${\ds K_s^{kl}:={\langle dM_s^k,dM_s^l\rangle\over ds}\, =\, \s^2\sum_{j=1}^d e^k_j(s)e^l_j(s)\, =
\,\s^2\, (e^k_0(s)e^l_0(s)-g^{kl}(\xi_s))} \, , $$ 
id est \quad ${\ds  K_s = \s^2 (e_0(s)\,^t\! e_0(s) - g\1(\xi_s)) \,}$, in accordance with Lemma \ref{lem.vert}. (Here $\,^t\! e_0\,$ denotes the transpose of the column-vector $\,e_0\,$, and $\,g\1\,$ denotes the inverse matrix of the pseudo-metric.)   \  $\,K_s$ has rank $d$ : we have indeed $\;K g\, e_0= (e_0\,^t\! e_0-g\1 )\, g\, e_0 = 0\,$, since $\,^t\! e_0\,g\,e_0 =1$. \par 
   Note that $\,K_s\,$ does not depend on the other frame vectors $\,e_j\,(j\ge 1)$, proving again 
that the projection $(\xi_s,\dot\xi_s)$ is a diffusion on the tangent bundle, as Theorem (\ref{the.gen},1) asserts. 
}\end{rem}

\noindent     
\ub{Proof of Theorem \ref{the.gen}} \quad 1) does not need any further proof. \par 
2) \  The process $\,(\zeta_s = \overleftarrow{\xi}(s)\,\dot\xi_s)_{s\ge 0}\,$ is continuous and lives on the fixed unit pseudo-sphere $\,T_{\xi_0}^1\M\,$. \  Let 
$\;\overrightarrow{\xi}(s) : T_{\xi_0}\M\rightarrow T_{\xi_s}\M\;$ denote the parallel transport along the $C^1$ curve $\,(\xi_{s'}\,|\,0\le s'\le s)$, and recall that 
$$ {\ts {d\over ds}}\,\overrightarrow{\xi}(s)^\ell_j\, =\, -\, \overrightarrow{\xi}(s)_j^k\,\G_{k m}^{\ell}(\xi_s) \,\dot\xi_s^m\, ,  \quad \mbox{which implies that} \quad  {\ts {d\over ds}}\, \overleftarrow{\xi}(s)^i_\ell =  \overleftarrow{\xi} (s)^i_q\,\G_{\ell m}^{q}(\xi_s)\,\dot\xi_s^m \,  ; $$ 
indeed 
$$ ({\ts {d\over ds}}\,\overleftarrow{\xi}(s)^i_\ell)\,\overrightarrow{\xi}(s)^\ell_j =  \overleftarrow{\xi} (s)^i_\ell\,\overrightarrow{\xi}(s)_j^k\,\G_{k m}^{\ell}(\xi_s)\,\dot\xi_s^m 
= \overrightarrow{\xi}(s)_j^k\,\overleftarrow{\xi} (s)^i_q\,\G_{k m}^{q}(\xi_s)\,\dot\xi_s^m \, , $$ 
whence 
$$ {\ts {d\over ds}}\,\overleftarrow{\xi}(s)^i_\ell =  \overleftarrow{\xi}(s)_\ell^j\, \overrightarrow{\xi}(s)_j^k\,\overleftarrow{\xi} (s)^i_q\,\G_{k m}^{q}(\xi_s) \,\dot\xi_s^m 
=  \overleftarrow{\xi} (s)^i_q\,\G_{\ell m}^{q}(\xi_s)\,\dot\xi_s^m \, .  $$ 
\par \smallskip \noindent 
Recall then from Remark \ref{rem.loc} that (for $\,0\le \ell\le d\,$) 
$$ d\dot\xi_s^\ell = \s\sum_{k=1}^d e^\ell_k(s)\circ dw^{k}_s -\G_{jk}^\ell (\xi_s)\,\dot\xi_s^j\,\dot\xi_s^k \, ds \, . $$ 
Therefore we get 
$$ d\zeta_s^i\; =\; (\circ\, d\,\overleftarrow{\xi}(s)^i_\ell)\,\dot\xi_s^\ell + 
\overleftarrow{\xi}(s)^i_\ell\circ d\dot\xi^\ell_s $$
$$ =\, \overleftarrow{\xi} (s)^i_q\,\G_{\ell m}^{q}(\xi_s)\,\dot\xi_s^m\,\dot\xi_s^\ell \, ds + \s\, 
\overleftarrow{\xi}(s)^i_\ell \sum_{k=1}^d e^\ell_k(s)\circ dw^{k}_s - \overleftarrow{\xi}(s)^i_\ell\,\G_{jk}^\ell (\xi_s)\,\dot\xi_s^j\,\dot\xi_s^k \, ds  $$ 
$$ =\; \s\sum_{k=1}^d \overleftarrow{\xi}(s)^i_\ell\,e^\ell_k(s)\circ dw^{k}_s 
=\; \s\sum_{k=1}^d \tilde e^i_k(s)\circ dw^{k}_s \; , $$ 
where $\;\tilde e_k(s) := \overleftarrow{\xi}(s)\,e_k(s)\,$, for $\,1\le k\le d\,$ and $\,s\ge 0\,$. \par 
   Similarly, we have \quad  
${\ds d \tilde e_k^\ell(s) = \s\, \overleftarrow{\xi}(s)^\ell_j\, \dot\xi_s^j\circ dw^k_s\, } $, \  id est   \quad  $ d \tilde e_k(s) = \s\,\zeta_s\circ dw_s^k\,$. \par 
   Observe that, for any $\,s\ge 0\,$, $\Big(\zeta_s , \tilde e_1(s) , .. , \tilde e_d(s)\Big)$ constitutes a pseudo-orthonormal basis of the fixed tangent space $\,T_{\xi_0}\M\,$. \ Hence, owing to Section \ref{sec.relr} and Remark \ref{rem.loc}, we find that the velocity process $\,\zeta_\cdot\,$ defines a hyperbolic Brownian motion on the hyperbolic space $\,T_{\xi_0}^1\M$, isometric to $\,\H^d$. \par 
   In the reverse direction, we have of course $\;\dot\xi_s = \overrightarrow{\xi}(s)\,\zeta_s\,$, meaning indeed that we recover the $\,C^1$ curve $\,\xi_\cdot\,$ as the deterministic development of the flat relativistic diffusion $\,\int_0^\cdot \zeta_s ds\,$. 
$\;\diamond$  
    
\begin{rem}\label{rem.intr}\quad {\rm The equation $\,(*)\,$ can be expressed intrinsically, in 
Stratonovitch or in It\^o form, by using the covariant differential $\,D\,$, which is defined in local
coordinates $\,(x^i,e_j^k)\,$ by 
$\;(D e_j)^k := de_j^k+\G_{\ell i}^k\,e_j^\ell\, dx^i\;$, for $\,0\le j,k\le d\,$. The equation $\,(*)\,$ is indeed equivalent to : 
\vspace{-2mm}
$$ \dot{\xi}_s = e_0(s) \quad ;\quad D e_0(s) = \s\sum_{j=1}^de_j(s)\circ dw^{j}_s \quad ;\quad 
D e_j(s) = \s\, e_0(s)\circ dw^{j}_s \quad \hbox{for } 1\le j\le d\; . $$ 
See for example ([B], page 30) or ([Em], page 427). 
}\end{rem}

\subsection{The restricted Schwarzschild space $\,{\cal S}_0\,$} \label{sec.S}\indf 
   This space is commonly used in physics to model the complement of a spherical body, star or black hole ;  see for example [DF-C], [F-N], [L-L], [M-T-W], [S]. \par    
   We take $\;\M = {\cal S}_0 := \Big\{ \xi = (t,r,\t)\in\R\times [R,+\ii [\times\S^2\Big\}\,$, where $\,R\in\R_+\,$ is a parameter of the central body, endowed with the radial pseudo-metric : 
$$ (1-{\ts{R\over r}})\, dt^2 - (1-{\ts{R\over r}})\1 dr^2 - r^2 |d\t|^2\, . $$ 
The coordinate $\,t\,$ represents the absolute time, and $\,r\,$ the distance from the origin. 

   In spherical coordinates $\,\t =(\f,\psi )\in [0,\pi ]\times (\R/2\pi\Z )$, we have 
\mbox{$\,|d\t|^2\! =d\f^2+\sin^2\!\f d\psi^2$.} \ The geodesics are associated with the Lagrangian 
$\;L(\dot\xi , \xi)\;$, where  
$$ 2\,L(\dot\xi,\xi) = (1-{\ts{R\over r}})\,\dot t^2 - (1-{\ts{R\over r}})\1 \,\dot r^2 - r^2 \dot\f^2- 
r^2\,\sin^2\f\;\dot\psi^2\; , $$   
and the non-vanishing Christoffel symbols are :
$$ \G_{rt}^t = -\G_{rr}^r = {R\over 2r(r-R)}\; ;\; \G_{tt}^r = {R(r-R)\over 2r^3}\; ;\; 
\G_{\f\f}^r = R-r\; ;\; \G_{\psi\psi}^r = (R-r)\sin^2\f\; ;\; $$ 
$$ \G_{r\f}^\f = \G_{r\psi}^\psi = r\1\; ;\; \G_{\psi\psi}^\f = -\sin\f\cos\f\; ;\;
\G_{\f\psi}^\psi = \cotg\f\; . $$ 
The Ricci tensor vanishes, the space $\,{\cal S}_0\,$ being empty.   A theorem of Birkhoff (see [M-T-W]) asserts that there is no other radial pseudo-metric  in $\,{\cal S}_0\,$ which satisfies this constraint.  \ The limiting case $\,R=0\,$ is the flat case of special relativity, considered in section \ref{sec.relr}.

\subsubsection{The stochastic differential system in spherical coordinates} \label{sec.syst}\indf 
   Let us take as local coordinates the global spherical coordinates : 
$$ \xi\equiv (\xi^0,\xi^1,\xi^2,\xi^3) := (t,r,\f,\psi )\; . $$  

   According to Remark \ref{rem.loc}, the system of It\^o stochastic differential equations governing  the relativistic diffusion $(\xi_s,\dot\xi_s)$ writes here as follows  : 
$$ dt_s = e^0_0(s)\,ds\; ,\; dr_s = e^1_0(s)\,ds\; ,\; d\f_s = e^2_0(s)\,ds\; ,\; d\psi_s =
e^3_0(s)\,ds\; , $$ 
$$ de^0_0(s) = {\ts{3\s^2\over 2}}\,e^0_0(s)\,ds -{\ts{R\over r_s(r_s-R)}}\,e^0_0(s)\,e^1_0(s)\,ds
+ dM_s^0\, , $$
\par\vspace{-2mm} 
\vbox{
$$ de^1_0(s) = {\ts{3\s^2\over 2}}\,e^1_0(s)\,ds +{\ts{R\over 2r_s(r_s-R)}}\, e^1_0(s)^2\, ds -
{\ts{R(r_s-R)\over 2r_s^3}}\,e^0_0(s)^2\, ds + (r_s-R)\,e^2_0(s)^2\, ds $$   \vspace{-6mm}
$$ \qquad +\, (r_s-R)\sin^2\f_s\,e^3_0(s)^2\, ds + dM_s^1\; , $$ }
\par\vspace{-6mm} 
$$ de^2_0(s) = {\ts{3\s^2\over 2}}\,e^2_0(s)\,ds -{\ts{2\over r_s}}\, e^1_0(s)\,e^2_0(s)\, ds
+\sin\f_s\cos\f_s\,e^3_0(s)^2\, ds + dM_s^2\; , $$
$$ de^3_0(s) = {\ts{3\s^2\over 2}}\,e^3_0(s)\,ds -{\ts{2\over r_s}}\, e^1_0(s)\, e^3_0(s)\, ds -
2\,\cotg\f_s\, e^2_0(s)\,e^3_0(s)\, ds + dM_s^3\; , $$ 
where the martingale $\, M_s := (M_s^0,M_s^1,M_s^2,M_s^3)\,$ has the following rank 3 quadratic covariation matrix : \quad ${\ds K_s = \s^2\, (e_0(s)\,^t\! e_0(s) -g\1(\xi_s))\, }.$

\subsubsection{Energy and angular momentum} \label{sec.reduc} \indf 
   We shall use widely \ the angular momentum $\quad\bv := r^2\,\t\wedge\dot\t \,$, \parn 
the energy $\quad a:= (1-{\ts{R\over r}})\, \dot t = {\p L\over \p\dot t}\,$, \ and the norm of $\,\bv$ :
$\;b:= |\bv| = r^2\,U\,$, \ with $\,U := |\dot\t|\,$. \par 

   Set also $\quad T:= \dot r\,$, \ and accordingly 
$$ T_s:=\dot{r}_s=e^1_0(s)\; ,\;\; U_s:=|\dot\t_s| = \sqrt{e^2_0(s)^2 + \sin^2\f_s\,
e^3_0(s)^2}\; ,\;\hbox{ and } \;\;  D:=\min\{ s> 0\,|\, r_s = R\}\,. $$ 

   Standard computations yield the following :  
\begin{prop}\label{pr.diff1}\quad 1) \  The unit pseudo-norm relation (which expresses that the parameter $\,s\,$ is precisely the arc length, id est the so-called proper time) writes 
$$\;T^2_s=a^2_s-(1-R/r_s)(1+b^2_s/r^2_s) \, . $$ 

   2) The process $\,(r_s,a_s,b_s,T_s)\;$ is a degenerate diffusion, with lifetime $\,D$, which solves the following system of stochastic differential equations : 
$$ dr_s = T_s\, ds\; ,  \quad dT_s =  dM^T_s + {\ts{3\,\s^2\over 2}}\, T_s\, ds + 
(r_s- {\ts{3\over 2}}R)\,{b_s^2\over r_s^4}\, ds - {R\over 2r_s^2}\, ds \; , $$ 
$$ d a_s = dM^a_s + {\ts{3\,\s^2\over 2}}\, a_s\, ds  \; , \quad  
d b_s = dM_s^b + {\ts{3\,\s^2\over 2}}\, b_s\, ds + {\s^2\, r_s^2\over 2\, b_s}\, ds \; , $$ 
with quadratic covariation matrix of the local martingale $\, (M^a, M^b, M^T)\,$ given by 
$$ K'_s := \s^2\pmatrix { a_s^2-1+{R\over r_s} & a_s\,b_s & a_s\,T_s \cr 
a_s\,b_s & b_s^2+r_s^2 & b_s\,T_s \cr a_s\,T_s & b_s\,T_s & 
T^2_s+1-{\ts{R\over r_s}}\cr } . $$ 
\end{prop}

   We get in particular the following statement, in which the dimension is reduced. 
\begin{cor}\label{cor.diff1}\quad The process $\,(r_s,b_s,T_s)\;$ is a diffusion,
with lifetime $\,D\,$ and infinitesimal generator \par \vspace{-1mm} \vbox{
$$ \GG' := T{\p\over\p r} + {\s^2\over 2}\, (b^2 + r^2) \,{\p^2\over\p b^2} + {\s^2\over 2b}\,(3b^2 + r^2)
\,{\p\over\p b} + \s^2bT\,{\p^2\over\p b\p T}  \hskip 20mm $$
\vspace{-4mm} 
$$ \hskip 30mm  + \,{\s^2\over 2}\,\Big( T^2+1-{\ts{R\over r}}\Big)\,{\p^2\over\p T^2} +
\Big({{3\,\s^2\over 2}}\, T + (r-{\ts{3\over 2}}R)\,{b^2\over r^4} - {R\over 2r^2}\,\Big){\p\over\p T} 
\; . $$}
\end{cor}

   We have the following result on the behaviour of coordinate $\,a_s\,$. 
\begin{lem}\label{lem.cv}\quad There exist a standard real Brownian motion $\,w_s\,$, and a real 
process $\,\eta_s\,$, almost surely converging in $\R$ as $s\nea D\,$, such that 
$\; a_s = \exp (\s^2\, s+\s\, w_s+\eta_s)\;$ for all $\, s\in[0,D[\,$. \ In particular $a_s$ almost
surely cannot vanish, which means that time $t_s$ is always strictly increasing. 
\end{lem} 
\ub{Proof}\quad Proposition \ref{pr.diff1} above shows that \ 
$(a_s^2-1)\,\s^2\,ds\le \langle dM^a_s\rangle\le a_s^2\,\s^2\,ds\,$, for $\,0\le s <D\,$. \parn  
So that we have almost surely (as $s\to\ii\,$, when $D=\ii$) : 
$$ \log a_s - \log a_0 \,=\; 3\s^2 s/2-\5\int_0^sa_t\2\langle dM_t^a\rangle + \int_0^sa_t\1\, dM_t^a
\;\ge\; \s^2 s + \int_0^sa_t\1\, dM_t^a $$ 
$$ = \s^2 s + o\Big( \int_0^sa_t\2\langle dM_t^a\rangle\Big) = \s^2 s + o(s) \; . $$ 
Since $\,(1-{\ts{R\over r_s}})\le a_s^2\,$, this implies 
$\;{\ds\int_0^D (1-{\ts{R\over r_s}})\, a_s\2\, ds\; <\ii}\;$ almost surely. \par 
   Consider then the standard real Brownian motion $\,w\,$ defined by  \parn  
$ dM_s^a = \s\, \sqrt{a_s^2-(1-R/r_s)}\,dw_s\,$, and the process $\,\eta\,$ defined by the formula in the statement. \parn 
We have : 
$$ d\eta_s = d(\log a_s) -\s^2\,ds-\s\,dw_s = \5\,(1-R/r_s)\, a_s\2\,\s^2\, ds +
\Big(\sqrt{1-(1-R/r_s)a_s\2 }-1\Big)\,\s\, dw_s \, , $$ 
and then for any $\,s<D\,$ : 
$$ \eta_s = \eta_0 + {\ts{\s^2\over 2}}\int_0^s (1-R/r_t)\, a_t\2\, dt - \s\int_0^s {(1-R/r_t)\,
a_t\2\over 1+ \sqrt{1-(1-R/r_t)\, a_t\2 }}\, dw_t \; , $$ 
which converges almost surely to a finite limit as $s\nea D\,$, since almost surely for all $\,s\in ]0,D[ \,$ : 
$$ \Big\langle\int_0^s {(1-R/r_t)\, a_t\2\over 1+ \sqrt{1-(1-R/r_t)\, a_t\2 }}\, dw_t\Big\rangle \le 
\int_0^s \Big( (1-{\ts{R\over r_t}})\, a_t\2\Big)^2\, dt < \int_0^D (1-{\ts{R\over r_t}})\, a_t\2\,
ds < \ii\; . \;\;\diamond $$ 

\subsubsection{Asymptotic behaviour of the relativistic diffusion $\,(\xi_s,\dot\xi_s)\,$}
\label{sec.comp.r} \indf 
   We see in the appendix (Section \ref{Tlgeod}) that in the geodesic case $\,\s=0\,$ five types of behaviour can occur,  owing to the trajectory of $\,(r_s)\,$ ; it can be : 
\parn 
- running from $R$ to $+\ii$, or in the opposite direction ; \parn   
- running from $R$ to $R$ in finite proper time ; \qquad  - running from $+\ii$ to $+\ii$ ; \parn 
- running from $R$ to some $R_1$ or from $R_1$ to $+\ii$, or idem in the opposite direction ; \parn  
- running for endlessly in a bounded region away from $R$. \par  

   More detailed results can be found in [L-L] and mainly [M-T-W]. A full treatment is given, for future reference, in the appendix below (Section \ref{sec.geod}). \par 

   The stochastic case $\,\s\not= 0\,$ can be seen as a perturbation of the geodesic case $\,\s=0\,$ ;
however the asymptotic behaviour classification is quite different.  

\begin{theo}\label{the.asymp}\quad 1) \ For any initial condition, the radial process $\,(r_s)\,$ almost
surely reaches $R$ within a finite proper time $D$ or goes to $+\ii$ as $s\to +\ii$ (equivalently : as 
$\,t(s)\to +\ii$ if $\,a_0>0$, and \ as $\,t(s)\to -\ii$ if $\,a_0<0$).
\parn 
2) \ Both events in 1) above occur with positive probability, from any initial condition. \parn 
3) \ Conditionally on the event $\,\{ D=\ii\}\,$ of non-reaching the
central body, the Schwarzschild relativistic diffusion $\,(\xi_s,\dot\xi_s)\,$ goes almost surely to
infinity in some  random asymptotic direction of $\,\R^3$, asymptotically with the velocity of light. 
\end{theo} 

   Note in particular that the relativistic diffusion almost surely cannot explode before the finite proper time $\,D\,$. \par 
   
   The proof we give for this theorem is rather long. It is postponed till Section \ref{sec.P} below. 

\subsection{The full Schwarzschild space $\cal{S}$}\label{sec.bh} \indf 
   The full Schwarzschild space $\,\cal{S}$, also known as the Kruskal-Szekeres space (see
[DF-C], [F-N], [S], and especially [M-T-W] and its historical account page 822),  can be defined by 
extending the  previous restricted Schwarzschild space ${\cal S}_0$ as follows. \quad 
On $\,{\cal S}_0\,$, \ set 
$$ u :=\sqrt{{\ts{r\over R}}-1}\times e^{r/(2R)}\times \ch({\ts{t\over 2 R}}) \quad \hbox{ and } \quad 
v :=\sqrt{{\ts{r\over R}}-1}\times e^{r/(2R)}\times \sh({\ts{t\over 2 R}}) \; . $$ 
Note that $\quad {\ds \Big({{r\over R}}-1\Big)\times e^{r/R} = u^2-v^2\,}$, and that 
the Schwarzschild  pseudo-metric expresses in the Kruskal-Szekeres coordinates $\,(u,v,\t)\,$ as : 
$$ (**) \qquad {{4R^3\over r}}\,\,e^{r/R} \,( dv^2 - du^2) - r^2 |d\t|^2\, , $$ 
where $\,r = {\bf r}(u^2-v^2)\,$,  $\, {\bf r}\,$ denoting the inverse function of $\,[r\mapsto ({{r\over R}}-1)\, e^{r/R}]\,$ (which is an increasing diffeomorphism from $\R_+$ onto $[-1,+\ii[\,$). \par 
   In those Kruskal-Szekeres coordinates, we have \ 
$\;{\cal S}_0 = \Big\{ (u,v,\t)\in \R^2\times\S^2 \,\Big|\, u > |v|\Big\}$. \par    
   The full Schwarzschild space $\,\cal{S}\,$ is now defined as 
$$ {\cal S} := \Big\{ (u,v,\t)\in \R^2\times\S^2 \,\Big|\, u^2-v^2 > -1\Big\}\,, $$ 
and is equipped with the pseudo-metric defined by $\,(**)\,$ above and by $\,r = {\bf r}(u^2-v^2)\,$. 
\par 
  $\,\cal{S}\,$ contains $\;{\cal S}_0\,$, $\,- {\cal S}_0\,$  (isometric to ${\cal S}_0$), two isometric copies of the hole : \parn
  \centerline{$\,{\cal H} := \Big\{ (u,v,\t)\in {\cal S}\,\Big|\, v> |u|\Big\}\;$ and $\;- {\cal H}\,$,} \parn 
and the boundary between $\,\pm {\cal S}_0\,$ and $\,\pm{\cal H}\,$, which is $\,\{ r=R\} = \{ u=\pm v\}$. \par 
  It is a Lorentz manifold, to which our general construction of section \ref{sec.relg} applies.  \par 

   The energy and angular momentum are extended to $\,T^1\cal{S}\,$ by setting 
$$ a := {2 R^2\over r}\,\,e^{-r/R} \, (u \dot v - v\dot u )\quad  
\hbox{ and }\quad \bv := r^2\,\t\wedge\dot\t \; . $$ 

   As before we set \quad $ T:= \dot r = \,{2 R^2\over r}\, e^{-r/R}\,(u\dot u-v\dot v)$ \quad and \quad $ b:=|\bv| = r^2\,U\,$. \ Recall that $\;a\,$ and $\,\bv\;$ are constant along geodesics. The unit pseudo-norm relation writes as before : 
$$ a^2 - T^2 + ({\ts{R\over r}}-1)\,\Big(\,{b^2\over r^2}+1\Big) = 0\, , $$  
or equivalently \quad ${\ds \dot u^2 - \dot  v^2 + {r\,e^{-r/R}\over 4 R^3}\,\Big(\,{b^2\over r^2}+1\Big) = 0\, ,}$ \  which implies $\,|\dot  v|>|\dot u|\,$, whence $\,\dot  v>|\dot u|\,$ along a timelike path. 
 The following correspondences between a line-element $(\xi,\dot\xi)\in T^1{\cal S}\,$ and its projection $\,\xi\in{\cal S}\,$ are easily deduced : \ $(\xi,\dot\xi)\in \{ r>R\,;\,a>0\}\LRa \xi\in {\cal S}_0 \,$, \par  
$(\xi,\dot\xi)\in \{ r>R\,;\,a<0\}\LRa \xi\in -{\cal S}_0\,$, \ $(\xi,\dot\xi)\in \{ r<R\,;\,T<0\}\LRa \xi\in -{\cal H}\,$, \par  
$(\xi,\dot\xi)\in \{ r<R\,;\,T>0\}\LRa \xi\in {\cal H}\,$, \ $(\xi,\dot\xi)\in \{ r=R\,;\,a=T=0\}\LRa \xi\in \{u=v=0\}$. 
\par \medskip 

    Now two other coordinates, namely the so-called inward and outward Eddington-Finkelstein coordinates $\,u^-\,$ and $\,u^+$,  prove to be very convenient for performing calculations on $\,\cal{S}\,$. They are defined, not on the whole $\,\cal{S}\,$, but  $\,u^-\,$ on $\;\if{\overline}\fi {\cal S}\cap\{ u+v\not= 0\}\,$  and $\,u^+\,$ on $\;\if{\overline}\fi {\cal S}\cap\{ u-v\not= 0\}$, by : 
 
$$ u^- := 2R\,\log |u+v|\quad \hbox{ and } \quad u^+ := -\, 2R\,\log |u-v| \; .  $$ 

   In those Eddington-Finkelstein  coordinates, the metric expresses as :  
$$ (1-{\ts{R\over r}})\, (du^{\pm})^2 \pm 2\,du^{\pm}\,dr\, - r^2 |d\t|^2\; ,  $$ 
and the energy expresses as
$$ a = (1-{\ts{R\over r}})\,\dot u^{\pm} \pm \dot r = {\p L\over\p\dot u^\pm}\; . $$ 

   We shall need (from Section \ref{sec.reg}) to complete the space $\,{\cal S}$, in $\,\overline{\cal S} := {\cal S}\sqcup \p{\cal S}$, where 
   $$ \p{\cal S} := \Big\{ \pm (u,v,\t)\,\Big|\, (u,v,\t)\in \R^2\times\S^2 \;;\; u^2-v^2 = -1\Big\}\, , $$ meaning that we identify the opposite points above the singularity $\{ r=0\}$. Note that $\,r, u^+,u^-\,$ are naturally continued to $\,\overline{\cal S}$, that \ $\; r=0\LRA v^2=u^2+1\LRA u^+=u^-\;$ in $\,\overline{\cal S}\,$, \quad and that  
$$  u^+ + r + R\,\log |{\ts{r\over R}}-1|\; =\; u^- - r - R\,\log |{\ts{r\over R}}-1| \quad \hbox{ on } 
\,\;\overline{\cal S}\cap\{ |u|\not= |v|\}\, , $$  
this quantity being equal to $\,t\,$ on $\,\pm\,{\cal S}_0\,$. \ The $\overline{\R}$-valued absolute time $\,t\,$ is continuously extended to $\,\if{\overline}\fi {\cal S}\cap\{ u=v=0\}^{\bf c}$, by setting 
$\,t=+\ii\,$ on $\,\{ u=v\not= 0\}$,  $\,t=-\ii\,$ on $\,\{ u=-v\not= 0\}$, \ 
$\,t= 2R\, {\rm argth} (u/v)\,$ on $\,\pm{\cal H}$, \ and $\,t=  2R\, {\rm argth} (v/u)\,$ on $\,-{\cal S}_0\,$ (as on $\,{\cal S}_0$). \par   
   Note that the region $\,\AA := \{ u=v=0\}$ appears right away as exceptional, as the only part of 
$\,\overline{\cal S}\,$ where $\,t\,$ cannot be continued, and the only part of $\,\overline{\cal S}\,$ where both $\,u^{\pm}\,$ explode. \par \smallskip 
\if{   Note (from the expressions of $\,a,T\,$ above) that the exceptional region $\{ u=v=0\}$ of $\,{\cal S}\,$ identifies with the part  $\,\AA :=\{ r-R =\ T =a=0 \} \,$ of the tangent space $\,T^1{\cal S}$.  \par \smallskip 
}\fi
   In the Kruskal-Szekeres coordinates, a path $(u_s, v_s, \t_s)$ is timelike if and only if  
$\;\dot v_s > |\dot u_s|\,$. 
This implies that any timelike path started in $\,{\cal H}\,$ has to hit  ${\ds \{ v = \sqrt{u^2+1}\,\} \subset \{ r=0\}}$, and that any timelike path started in $\,-{\cal H}\,$ has to hit $\{ r=R\}$. In particular, a timelike path started in $\,-{\cal H}\,$ and avoiding $\,\AA\,$ has to enter the region $\{ r>R\}$ through $\,\{ t=-\ii\}$, appearing either as a particle born and then evolving in $\,{\cal S}_0\,$ (for ever, or entering then $\,{\cal H}\,$ through $\,\{ t=+\ii\}$), or the analogue through $\,-{\cal S}_0\,$, which could be viewed as the case of an antiparticle through $\,{\cal S}_0\,$.  
\par 
      These dynamics of the full Schwarzschild space $\,\cal{S}\,$ show that the inward Eddington-Finkelstein  coordinate $\,u^-\,$ is appropriate to the study of timelike paths started in $\,{\cal S}_0\cup{\cal H}$, till they hit $\{ r=0\}$ (and even to extend them further, see section \ref{sec.reg} below), and that the outward Eddington-Finkelstein  coordinate $\,u^+\,$ is appropriate to the study of timelike paths started in $\,-\overline{\cal H}\,$ and entering $\,\pm{\cal S}_0\,$, till they hit $\,{\cal H}$. 

\subsubsection{Diffusion in the full space $\cal{S}$ : hitting the singularity}\label{sec.bhh} \indf 
   We follow the same route as in the restricted Schwarzschild space, to express the relativistic
diffusion on $\,T^1\cal{S}$. Let us proceed, using the Eddington-Finkelstein coordinates $\,u^\pm$. \ The Lagrangian $\;L(\dot\xi,\xi)\;$ writes 
$$ 2\,L(\dot\xi,\xi) = (1-{\ts{R\over r}})\,(\dot u^{\pm})^2\pm 2\,\dot u^{\pm}\,\dot r - r^2 \dot\f^2- 
r^2\,\sin^2\f\;\dot\psi^2  \, . $$   

   We apply then Remark \ref{rem.loc}, to get the It™ stochastic differential equations of the
relativistic diffusion in the full Schwarzschild space : 
$$ d\dot u^{\pm}_s = \s\,dM^{\pm}_s+{\ts{3\s^2\over 2}}\,\dot u^{\pm}_s\,ds \pm{\ts{R\over 2r_s^2}}\,(\dot
u^{\pm}_s)^2\,ds \mp r_s\,(\dot\f^2_s+\sin^2\f_s\;\dot\psi^2_s)\,ds \; ; $$ 
$$ d\dot r_s = \s\,dM^r_s+{\ts{3\s^2\over 2}}\,\dot r_s\,ds + ({\ts{R\over r_s}}-1){\ts{R\over
2r_s^2}}\,(\dot u^{\pm}_s)^2\,ds \mp {\ts{R\over r_s^2}}\,\dot u^{\pm}_s\,\dot r_s\,ds 
+ (r_s-R)\,(\dot\f^2_s+\sin^2\f_s\;\dot\psi^2_s)\,ds \; ; $$ 
$$ d\dot \f_s = \s\,dM^\f_s+{\ts{3\s^2\over 2}}\,\dot \f_s\,ds - {\ts{2\over r_s}}\,\dot r_s\,\dot
\f_s\,ds  + \sin\f_s\,\cos\f_s\;\dot\psi^2_s\,ds \; ; $$ 
$$ d\dot \psi_s = \s\,dM^\psi_s+{\ts{3\s^2\over 2}}\,\dot \psi_s\,ds - {\ts{2\over r_s}}\,\dot r_s\,\dot
\psi_s\,ds  - 2\, \cotg\f_s\;\dot\f_s\dot\psi_s\,ds \, , $$ 
for some continuous local martingale $\,(M^{\pm}_\cdot,M^r_\cdot,M^\f_\cdot,M^\psi_\cdot)\,$, having quadratic covariation  matrix (according to Lemma \ref{lem.vert} and Section \ref{sec.syst}) : 
$$ \pmatrix { (\dot u^{\pm}_s)^2 & \dot u^{\pm}_s\,\dot r_s\mp 1 & \dot u^{\pm}_s\,\dot \f_s &
\dot u^{\pm}_s\,\dot\psi_s \cr \dot u^{\pm}_s\,\dot r_s\mp 1 & \dot r_s^2+1-{R\over r_s} & \dot
r_s\,\dot \f_s & \dot r_s\,\dot \psi_s \cr \dot u^{\pm}_s\,\dot \f_s & \dot r_s\,\dot \f_s & \dot
\f_s^2+r_s\2 & \dot \f_s\,\dot\psi_s \cr \dot u^{\pm}_s\,\dot \psi_s & \dot r_s\,\dot \psi_s & \dot
\f_s\,\dot \psi_s & \dot\psi^2_s+(r_s\sin\f_s)\2\cr } . $$  

   We have again a reduced diffusion $\,(r_s, b_s, T_s)\,$ (with minimal dimension), solving the same system of It™ stochastic differential equations as before. \par 

   This system of stochastic differential equations has been derived using Eddington-Finkelstein coordinates, so that it is valid a priori outside $\{ u=v=0\}$. But the smooth functions $\,(r, a,b,T)\,$ of the relativistic diffusion have an It\^o decomposition with continuous coefficients, so that the formulas involving them hold without restriction. \par 
\if{
$$ dr_s=T_s\, ds\; ;\quad da_s= \s\,dM^a_s + {\ts{3\s^2\over 2}}\,a_s\,ds \; ; \quad 
db_s = \s\,dM^b_s + {\ts{3\s^2\over 2}}\,b_s\,ds +{\ts{\s^2\,r_s^2\over 2\,b_s}}\,ds\; ; $$ 
$$ dT_s = \s\,dM^T_s + {\ts{3\s^2\over 2}}\,T_s\,ds + (r_s-{\ts{3R\over 2}})\,r_s^{-4}\, b_s^2\,ds - {\ts{R\over 2\,r_s^2}}\,ds \, , $$ 
for some continuous local martingale $\,(M^a_\cdot,M^b_\cdot,M^T_\cdot)\,$ having quadratic covariation  matrix 
$$ K''_s := \pmatrix { a_s^2+{R\over r_s}-1 & a_s\,b_s & a_s\,T_s \cr 
a_s\,b_s & b_s^2+r_s^2 & b_s\,T_s \cr a_s\,T_s & b_s\,T_s & T^2_s+1-{\ts{R\over r_s}}\cr } . $$ 
We have then the same generator as before for the degenerate diffusion $\,(r_s,a_s,b_s,T_s)$, and a fortiori for $\,(r_s,a_s,T_s)$, $(r_s, b_s,T_s)$, which are again autonomous degenerate diffusions
(with minimal dimension). \par \smallskip  
}\fi
   From the pseudo-norm relation, we see that $\,T_s\,$ cannot vanish in the region $\,\{ r<R\}$. \ As
$\,r_s\,$ enters this region necessarily almost surely with derivative $\,T_D<0\,$ (indeed
$|T_D|=|a_D|>0$ by Lemma \ref{lem.cv}),  $\,r_s\,$ is then necessarily strictly decreasing. Precisely, we have the following. 

\begin{theo}\label{the.bh}\quad  The relativistic diffusion in $\,T^1{\cal S}\,$ either escapes to infinity, or enters above  $\,{\cal H}\,$ at time $\,D\,$ and converges to the singularity within some finite proper time $\,D'$. \ Moreover, in the second case we have almost surely :   \par 
$1)$ \ for $\,s\ge D$, $\,r_s\,$ decreases and hits 0 at proper time $D'$, with $\,D<D'\le D+{\pi\over
2}\,R\,$ ; moreover $\;\lim_{s\nea D'}\,T_s = -\ii\,$ ; \par 
$2)$ \ $\,\t_s\,,\, u^\pm_s\,,\, \bv_s\,$ and $\,a_s\,$ converge to finite limits as $\,s\nea D'$, and
$\;b_{D'}\,$ cannot vanish ; \par 
$3)$ \ as $\,s\nea D'$, we have the following equivalents : 
$$ r_s \sim [{\ts{5\over 2}}\,b_{D'}\sqrt{R}\,(D'-s)]^{2/5}\quad \hbox{ and }
\quad T_s \sim -b_{D'}\sqrt{R}\times [{\ts{5\over 2}}\,b_{D'}\sqrt{R}\,(D'-s)]^{-3/5}\, . $$ 
\end{theo} 
\smallskip 

\begin{rem}\label{rem.loglog}\quad {\rm  The equivalents in Theorem (\ref{the.bh},3) above 
can be specified further. Indeed, we have almost surely, as $\,s\sea 0\,$ : 
$$ r_{D'-s} = [{\ts{5\over 2}}\,b_{D'}\sqrt{R}\,s]^{2/5}\times \Big(
 1 + ({\ts{5b_{D'}\over 2R^2}}\,s)^{2/5}+ \O({\ts\sqrt{s\log |\log s|}}\,)\Big) \, , $$
and 
$$ T_{D'-s} = -b_{D'}\sqrt{R}\times [{\ts{5\over 2}}\,b_{D'}\sqrt{R}\,s]^{-3/5}\times \Big(
 1 - 2 ({\ts{5b_{D'}\over 2R^2}}\,s)^{2/5}+ \O({\ts\sqrt{s\log |\log s|}}\,)\Big)\, . $$ 
Indeed, using the stochastic differential equation of $\,T_s\,$ and the iterated logarithm law, together with the equivalents in Theorem (\ref{the.bh},3), we deduce easily these more precise asymptotic expansions near $D'$. 
} \end{rem} 
 
\begin{rem}\label{rem.-S}\quad {\rm    We know from Theorem \ref{the.gen} that the relativistic diffusion can start from any initial condition in the full space $\,T^1S\,$. When it starts above  $\,-{\cal H}\,$, the pseudo-norm relation forbids any vanishing of $\,T_s\,$ (which has then to remain $\,>0\,$), till the level $\{r=R\}$ is hit, which takes a proper time less than $\,\pi R/2\,$, for the very same reason as in Theorem (\ref{the.bh},1). \ When the diffusion starts above $\{ r=R\}$, it enters $\{ r\not= R\}$ at once, as any timelike path.     Note that above $\,\AA\equiv \{u=v=0\}\subset \{r=R\}$, we have necessarily $\,T=a=0\,$. Moreover it can be proved that $\,T^1\AA\,$ is polar for the relativistic diffusion. \ So, when starting above $\,-{\cal H}$, the relativistic diffusion enters then above $\,\pm{\cal S}_0\,$, before possibly entering 
later above $\,{\cal H}$. 
}\end{rem} 
 
    We postpone the proof of Theorem \ref{the.bh} till Section \ref{sec.P} below. \par \smallskip 
    A part of this proof is based on the following proposition, which allows to recover the whole relativistic diffusion from the reduced relativistic diffusion $\,(r_s, b_s, T_s)$. 
\begin{prop}\label{pro.systeqst}\quad The spherical coordinate $\,\t_s\,$ satisfies the following
stochastic differential  equation (conditionally on the reduced relativistic diffusion $\,(r_s, b_s,
T_s)$) :  
$$ d\Big({\dot\t_s\over U_s}\Big) = \Big({r_s\over b_s}\,\s\,d\b_s\Big)\, \t_s\wedge{\dot\t_s\over U_s}\, -\,\Big({b_s\over r^2_s}\,ds\Big)\,\t_s - \Big({\s^2\,r_s^2\over 2\,b_s^2}\,ds\Big)\, {\dot\t_s\over
U_s}\; , $$  
for some standard real Brownian motion $\,\b_\cdot\,$, which is independent of 
$\,(r_\cdot, b_\cdot,T_\cdot)\,$. \par 
   Moreover, $\;\dot\t_s/U_s\;$ converges in $\,\S^2$ as $\,s\nea D'$, almost surely. 
\end{prop}

 We postpone the proof of this proposition till Section \ref{sec.P} below. \par 

\begin{cor}\label{cor.bh}\quad The curve in the full space $\cal{S}$ defined by the image of the
trajectory \parn $\{(r_s,u_s^-,\t_s)\,|\,s\le D'\}\,$ admits almost surely a semi-tangent at the center
of the hole. 
\end{cor}
\ub{Proof}\quad  Using the strict monotonicity of $\,r_s\,$ near the singularity, we see that it is
sufficient to verify the left-differentiability of the curve $\,(r_s\mapsto (u^-_s,r_s\t_s))\,$ at
$s=D'$. \ Now using Theorem (\ref{cor.bh},3), as $\,s\nea D'\,$, on one hand we have \par 
$\p u^-_s/\p r_s = \dot u^-_s/T_s = (a_s/T_s+1)\,r_s/(r_s-R)\sim -r_s/R \ra 0\,$, \ and on the other hand \par 
$\p (r_s\t_s)/\p r_s = \t_s+(\dot\t_sr_s/T_s) \ra \t_{D'}\in\S^2\,$,  \ since \par 
$|\dot\t_sr_s/T_s| = b_s/(r_sT_s) \sim -R^{-1/2}\,[{\ts{5\over 2}}\,b_{D'}\sqrt{R}\,(D'-s)]^{1/5}
\ra 0\,$. $\;\diamond $ 

\begin{rem}\label{rem.expl}\quad {\rm Since we have $\, |\dot\t_s|=U_s=b_s/r_s^2\,$, we observe the explosion of the spherical speed $\dot\t_s$, as well as of the radial speed $T_s$, and as of the speed $\,\dot u^-_s$, at proper time $D'$, id est at the hitting of the singularity $\{ r=0\}$. 
\ However we also just saw that $\,(r_s,\bv_s,a_s,\t_s,u^-_s)\,$ is left-continuous at $D'$, and that
moreover the curve in the space $\cal{S}$ defined by the image of the trajectory
$\,\{(r_s,u_s^-,\t_s)\,|\,s\le D'\}\,$ admits almost surely a tangent at any point. Indeed it happens
that the explosion of the derivatives does not forbid to define a continuation of the relativistic
diffusion after the finite hitting proper time $D'$. This will be indeed the purpose of Section
\ref{sec.reg} below.    }\end{rem}

\subsubsection{Regeneration through the singularity : entrance law after $D'$}\label{sec.reg} \indf
   We see from Theorem \ref{the.bh} and Remark \ref{rem.expl} that the set of endpoints of the relativistic paths (at proper hitting time $\,D'$ of $\{ r=0\}$) identifies with the boundary $\,\p {\cal S}$, which we defined in Section \ref{sec.bh}, identifying pointwise the outward $\p{\cal H}\,$ and inward $\p (-{\cal H})$ boundaries. Thus we have indeed  $\,\p {\cal S}\equiv \p\, T^1{\cal S}$.  \par 
    In this identification, a differentiable inward path ending at $\{ r=0\}$ can be continued by a differentiable outward path, so that the $\R^3$-valued curve $\,r\t\,$ is differentiable at any point. \ 
   In particular, geodesics are thus well defined for any proper time, and there are geodesics which cross endlessly the singularity $\{ r=0\}$, namely those which are described in case 1.2 (and are met also in cases 2.2.1, 2.2.2, 2.5.2, and 2.6), completed by Remark \ref{rem.fin}, of Section \ref{sec.geod}. For generic values of parameters $(a,\bv,k)$, such geodesics are dense in some disk of $\R^3$ (centred at $0$). \par 
      Since the diffusion can hit the singularity $\{ r=0\}$ in finite time, it is natural to look for an entrance law, allowing to continue it after time $\,D'\,$.  Clearly it as to enter above $\,-{\cal H}$. Thus we have to define for the diffusion on $\,T^1S\,$ a family of entrance laws above the singularity $\{ r=0\}$, and more precisely on the boundary $\,\p {\cal S}$.      

   Let $\,\GG\,$ denote the generator of the relativistic diffusion, acting on $\,C^2(T^1{\cal S})\,$.
Theorem \ref{the.bh} allows us to extend the relativistic diffusion to a continuous strong Markov process on $\,T^1{\cal S}\cup\p\, T^1{\cal S}$, provided we establish the following proposition.
\begin{prop}\label{the.martp}\quad The martingale problem associated with $\,\GG$ has a unique continuous solution, starting from any point of $\,\p {\cal S}$.  
\end{prop} 

   We postpone the proof of this proposition till Section \ref{sec.P} below.  

\subsubsection{The relativistic diffusion for all positive proper times}\label{sec.SDalls} \indf 
   The existence and uniqueness of the entrance law (in Proposition \ref{the.martp}) allows to prove the first assertion of the following. \par 
   \vbox{
\begin{theo}\label{the.alls}\quad  There exists a unique continuous strong Markov process on $\,\ol{\cal S}\,$ with positive lifetime inducing the relativistic diffusion on $\,{\cal S}$. \par 
   The lifetime of this extended relativistic diffusion is almost surely infinite.  
\end{theo} } \par 

   For such a process, we define an increasing sequence of hitting proper times $\,D_k\,$ as follows. \    Let $\,D_0:=D\in [0,\ii]\,$ denote the hitting time of $\{r\le R\}$, $\,D_1:=D'\,$ denote the hitting time of $\{r= 0\}$, and set by induction, for any $\,n\in\N^*$ : \par 
$\qquad D_{3n}:= \inf\{ s>D_{3n-1}\,|\,r_s=R\}\,,\quad D_{3n+1}:= \inf\{ s>D_{3n}\,|\,r_s=0\}\,$, \quad
and \par 
$\qquad D_{3n+2}:= \inf\{ s>D_{3n+1}\,|\,r_s=R\}$.  \quad Finally consider $\;D_\ii:=\sup_n\,D_n\,$. 
\par \smallskip 
   This is obviously an increasing sequence of stopping times, strictly increasing as long as it is finite. The preceding section \ref{sec.reg} extends in a unique way the law of the relativistic diffusion to the proper time interval $\,[0,D_\ii[\,$.  It is clear that the process cannot be extended continuously beyond $\,D_\ii\,$. \par 
   This proves the first assertion of Theorem \ref{the.alls}, within the lifetime $\,D_\ii\,$.    We postpone till Section \ref{sec.P} the proof of the second assertion of this theorem : 
$\,D_\ii\,$ is almost surely infinite. 

 \begin{rem}\label{rem.Ds}\quad {\rm  We saw in Theorem (\ref{the.bh},1) that $\,D_1\le D_0+{\pi\over 2}\,R\,$. Now exactly the same reason shows that $\; D_{3n+1}-D_{3n}\,$ and $\,D_{3n+2}-D_{3n+1}\;$ are $\,\le \,\pi\,R/2\,$, as long as these times are finite. \ The time intervals 
$[D_{3n-1},D_{3n}]$ correspond to the excursions outside the hole, and the time intervals $[D_{3n},D_{3n+2}]$ correspond to the excursions inside the hole. Moreover we see that $\,D_k\,$ becomes infinite if and only if $\,k=3n\,$ and the process escapes to infinity during its $\,n$-th excursion outside the hole. \par 
   Recall that (according to Lemma \ref{lem.cv}) during every excursion outside the hole $\{r\le R\}$,
the diffusion can have its absolute time coordinate $\,t=t_s\,$ strictly increasing from $-\ii$ to $+\ii$
or strictly decreasing from $+\ii$ to $-\ii$ (this case can be seen as the antiparticle case), depending on the sign of $\,a\,$ at the exit of the hole. \par 
   Moreover, the $\R^3$-valued curve $\,(s\mapsto r_s\,\t_s)\,$ is differentiable at any proper time $s$,  whereas the $\R$-valued curves $\,(s\mapsto u^\pm_s)\,$ present a cusp (with half a tangent : this appears in the proof of Theorem \ref{the.bh}, where we saw that $\,\dot u\-_{}\sim -rT/R \ra \pm\ii$ near $D'$) at proper times $D_{3n+1}$ (and are differentiable at any time $s\not= D_{3n+1}$). 
}\end{rem}

 \begin{rem}\label{rem.mesinv}\quad {\rm The Liouville measure is invariant for the extended relativistic diffusion on $\,\ol{\cal S}\,$. \ It induces the invariant measure $\; dr\,da\,dT\;$ of the autonomous  diffusion $\,(r_s,a_s,T_s)\,$.
 }\end{rem}

\subsubsection{Capture of the diffusion by a neighbourhood of the hole}\label{sec.NumbC}\indf 
   The preceding section leads naturally to the following question : can the extended
relativistic diffusion cross infinitely many times the hole, as some geodesics do ? \par 

   Note first that it was clear from the preceding sections \ref{sec.comp.r}, \ref{sec.reg}, and
\ref{sec.SDalls}, that there is, for any $\,n\in\N\,$ and any initial condition, a positive probability that the extended relativistic diffusion crosses exactly $\,n\,$ times the hole and thereafter goes away to infinity. So the following theorem asserts essentially that the limiting case of $\,n=\ii\,$ crossings of the hole can also happen, thereby completing the picture of all possible asymptotic behaviours of the extended relativistic diffusion. \par 
   Moreover it says that this last case corresponds to an asymptotic confinement of the relativistic diffusion in the vicinity of the hole.

\begin{theo}\label{the.piege} \  Almost surely, from any initial condition, the extended relativistic diffusion can have only two types of asymptotic behaviour, each occurring with positive probability : \\ 
Either \par 
   1) \ $\,r_\cdot\,$ and $\,|a_\cdot|\,$ go away to infinity, and $\,\t_\cdot\,$ converges ; \\ 
or \par  
   2) \ $\rr\,:=\,\limsup_{s\to\ii}\limits\,r_s\,\in [R,3R/2]\,$, $\;\liminf_{s\to\ii}\limits\,r_s\,=0\,$, $\;b_\cdot\,$ goes away to infinity, \par \quad \ 
$\vec{b}_\cdot/b_\cdot\,$ converges, and $\;a_\cdot/b_\cdot\,$ converges to $\;\ell = \pm\,{1\over\rr} \, \sqrt{1-{R\over\rr}}\,$. \par

   Moreover \par 
    $(i)$ \ For any $\,\e>0$, if $\,r_0>3R/2\,$ and if $\,T_0\,$ is large enough, then the probability
that the relativistic diffusion goes away to infinity is at least $\,1-\e$. \par 
    $(ii)$ \ For any $\,\e>0$, if $\,r_0\in ]R,3R/2[\,$, $\,T_0=0\,$, and $\,b_0\,$ is large enough, then 
\parn \centerline{$\P(|\rr-r_0|<\e)>1-\e\,$.} 
\end{theo}

   The proof of this theorem is rather delicate. It will be presented in Section \ref{sec.P} below. 

\begin{rem}\label{rem.harm}\quad   {\rm  For any open interval $]x,y[\,\subset\, ]R,3R/2[\,$, $\,\P_{\bf\cdot}(x<\rr <y)\,$ is a non-trivial $\,SO_3\,$-invariant $(\LL_0+{\ts{\s^2\over 2}}\D_v)$-harmonic function on $\,T^1\M\,$. 
}\end{rem}  

   The following result describes more precisely what happens when the diffusion is captured by a neighbourhood of the hole : while (according to Theorem \ref{the.piege}) they are asymptotically planar, they exhibit progressively another type of regularity. 
\if{    : as if       which can evoke roughly  become asymptotically regular,   
and (for $\,\rr\,$ smaller than some intermediate radius in $]R,3R/2[\,$) resemble the picture drawn by a fixed point (at distance $\,\rr/2\,$ of the center) of a disk (of radius say $\,r$) rolling (without slipping) within a circle (which is centred at 0 and has radius $\,r+\rr/2\,$, both being orthogonal to the limiting direction of the angular momentum $\,\vec{b}\,$). 
}\fi 
\begin{cor}\label{cor.piege}\quad In the second case of Theorem \ref{the.piege}, and more precisely conditionally on the event $\,\rr<3R/2$, the times $\, D_{n}'\,$ of the first maxima 
of the radius $\,r_\cdot\,$ at each excursion out of the hole are such that $\;\lim_{n\to\ii}\limits
r_{D_{n}'} =\rr\;$ and that $\;{\ds \int_{D_{3n+1}}^{D_{n}'} d\t_s }\;$ converges as $\,n\to\ii$, towards $\;{\ds \pm\, \int_0^{\rr}\! {dr\over\sqrt{[R-r+\ell^2 r^3]\, r}} }\,$ (which is well defined if and only if $\,\rr<3R/2$). 
\end{cor}

   The proof of this corollary will close Section \ref{sec.P} below. 
   
 \begin{rem} \label {rem.desc} \quad {\rm The result of Corollary \ref{cor.piege} concerns  the time intervals $\,[D_{3n+1}, D_{n}'] $, that is to say the upcrossings from the singularity to the successive tops of the limiting trajectories.  It is very likely that the same result is valid as well for the downcrossings, that is to say the time intervals $\,[D_{n-1}', D_{3n+1}] $, yielding the same angular random limit (the sign of $\,T_s\,$ compensating for the interchange of the bounds $\,D_{n-1}'\,$ and $\,D_{3n+1}\,$ in the integral). \parn 
So another statement in the spirit of Corollary \ref{cor.piege}, but which demands some more work, should be : almost surely \quad 
${\ds  \lim_{n\to\ii}\,  \int^{D_{n}'}_{D_{n-1}'} d\t_s\,  =\,  \pm\, 2\,\Psi \;}$, \ where  
$$ \Psi \, := \int_{0}^{\rr}  {dr\over \sqrt{r\,[R-r+\ell^2 \,r^3]}} =  \int_{0}^{1}  {dr\over \sqrt{r\,[(1-r^3)(R/\rr) -(1-r^2)r]}}\ $$ 
is a strictly increasing continuous function of $\,\rr/R\,$, from $[1, {\ts{3\over 2}}[\,$ onto $[{\ts{\pi\over 2}}, \ii [\,$. \par 
         It is thus likely that the shape of the excursions should approach more and more the null geodesics, id est the light rays. See the appendix, Section \ref{Ngeod}. 
\if{ The picture evoked just above Corollary \ref{cor.piege} is appropriate for $\,\psi <\pi\,$, in which case the radius $\,r\,$ of the disk is $\,r= \5\rr\psi/(\pi-\psi )$. 
}\fi     }
\end{rem} 

\section{Proofs} \label{sec.P}

\subsection{Proof of Theorem \ref{the.asymp}}\indf \label{sec.PT2}
   In the proof of this theorem, we shall use the following very simple lemma.
\begin{lem}\label{lem.mart}\quad Let $\,M_\cdot\,$ be a continuous local martingale, and $\,A_\cdot\,$ 
a process such that \parn 
$\liminf_{s\to\ii}\limits\; A_s/\langle M\rangle_s >0\,$ almost surely on $\,\{ \langle M\rangle_\ii
=\ii\}$. \ Then $\;\lim_{s\to\ii}\limits\, (M_s+A_s) = +\ii\;$ almost surely on $\,\{ \langle
M\rangle_\ii =\ii\}$. 
\end{lem}
\ub{Proof} \quad Writing $\;M_s = W(\langle M\rangle_s)\;$, for some real Brownian motion $\,W$, we find
almost surely some $\,\e >0\,$ and some $\,s_0\ge 0\,$ such that $\;A_s\ge 2\e\,\langle M\rangle_s\;$ and
$\;|M_s|\le\e\,\langle M\rangle_s\;$ for $\,s\ge s_0\,$. Whence $\;M_s+A_s \ge\e\;\langle M\rangle_s\;$
for $\,s\ge s_0\,$. $\;\diamond $ 
\par \medskip  
 
   We prove now successively the 3 assertions of Theorem \ref{the.asymp}. \par 
\smallskip 
   1) \ {\it Almost sure convergence on $\,\{ D=\ii\}$ of $\,\,r_s\,$ to $\,\ii$.}
\par\smallskip 
   This proof will be split into six parts. \par     
   Let us denote by $\,A\,$ the set of paths with infinite lifetime $\,D\,$ such that the radius $\,r_s\,$
does not go to infinity. \ We have to show that it is negligible for any initial condition 
$\,x=(r,b,T)=(r_0,b_0,T_0)\,$ belonging to the state space $\,[R,\ii[\times\R_+\times\R\,$. 
\par\smallskip 

   The cylinder $\;\{ r=3R/2\}\;$ plays a remarquable r\^ole in Schwarzschild geometry. In particular, it
contains light lines. We see in the following first part of proof that we have to deal with this
cylinder.  \par\smallskip 

   $(i)$ \ {\it $\,r_s\,$ must converge to $\,3R/2\,$, almost surely on $\,A\,$.} \par \smallskip 

   Observe from the unit pseudo-norm relation (Proposition \ref{pr.diff1}, 1) that $\,{|T_s|/a_s}\,$
is bounded by 1. \ Let us apply It\^o's formula to $\;{\ds Y_s:= (1-{\ts{3R\over 2\,r_s}})\,{T_s/
a_s}}\,$ : 
$$ Y_s = M_s + {\ts{3R\over 2}}\int_0^s {T^2_t\over a_t\,r^2_t}\, dt + \int_0^s (1-{\ts{3R\over 2r_t}})^2
{b^2_t\over a_t\,r^3_t}\, dt  - \,\s^2\int_0^s (1-{\ts{R\over r_t}})\,{Y_t\over a_t^2}\, dt - 
\int_0^s (1-{\ts{3R\over 2r_t}})\,{R\over 2a_t\,r_t^2}\, dt \; , $$  
with some local martingale $\,M\,$ having quadratic variation : 
$$ \langle dM_s\rangle\, =\, (1-{\ts{3R\over 2r_s}})^2(1-{\ts{R\over r_s}})\Big( 1-{T_s^2\over
a_s^2}\Big){\s^2\over a_s^2}\; ds \;\le\; \s^2\,a_s\2\,ds\, . $$ 

   Now $\,{|Y_s|}\,$ is also bounded by 1. Hence Lemma \ref{lem.cv} implies that the last two terms in the
expression of $\,Y_s\,$ above have almost surely finite limits  as $\,s\to\ii\,$. Idem for $\,\langle
M_s\rangle\,$, and then for $\,M_s\,$. Moreover the two remaining bounded variation terms in the
expression of $\,Y_s\,$ above increase. As a consequence, we get that $\,Y_s\,$, $\,{\ds\int_0^s
{T^2_t\over a_t\,r^2_t}\, dt}\,$, and $\,{\ds\int_0^s(1-{\ts{3R\over 2r_t}})^2\, {b^2_t\over a_t\,r^3_t}\,
dt}\,$ converge almost surely in
$\R$ as $\,s\to\ii\,$. So does also $\,{\ds\int_0^s {dt\over a_t\,r^2_t}}\,$. \par 

   Now using that $\; {\ds {a\over r^2}\le \Big( a+{R\over r}\,\Big({b^2\over a\,r^2}+{1\over
a}\Big)\Big)\, r\2 = {T^2\over a\,r^2}+{b^2\over a\,r^4}+{1\over a\,r^2} }\;$ by the unit pseudo-norm
relation, we deduce that almost surely 
$$ \int_0^\ii(1-{\ts{3R\over 2r_t}})^2\, \Big|{d\over dt}\,(1/r_t)\Big|\, dt =
\int_0^\ii(1-{\ts{3R\over 2r_t}})^2\, {|T_t|\over r^2_t}\, dt \le \int_0^\ii  (1-{\ts{3R\over 2r_t}})^2\,
{a_t\over r^2_t}\, dt <\ii\; . $$ 
This implies the almost sure convergence of $\;{ \Big( 1-{3R\over 2r_s}+{3R^2\over 4r_s^2}\Big)\Big/
r_s }\;$, and therefore of $\,(1/r_s)\,$. \par 
  Since $\,\lim_{s\to\ii}\limits(1/r_s)\,$ cannot be 0 on $A\,$, we have necessarily 
$\,\lim_{s\to\ii}\limits r_s = 3R/2\,$ almost surely on $A\,$, from the convergence of 
$\,{\ds\int_0^\ii  (1-{\ts{3R\over 2r_t}})^2\,{a_t\over r^2_t}\, dt }\;$. \par \bigskip

   $(ii)$ \ {\it $\;{b_s/a_s}\,$ converges to $\,3R\sqrt{3}/2\,$, and $\,T_s/b_s\,$ goes to 0, almost
surely on $A\,$.} \par \smallskip 
   Indeed, It™'s formula gives (for some real Brownian motion $w$) 
$$ {b^2_s\over a^2_s} = {b^2_0\over a^2_0}+ 2\s\int_0^s {b_s^2\over a_s^2}\,\sqrt{{r_s^2\over
b_s^2}-{1-{R\over r_s}\over a_s^2}}\, dw_s +2\s^2\int_0^s{r_s^2\over a_s^2}\,ds - 3\s^2\int_0^s (1-{R\over
r_s})\,{b^2_s\over a^4_s}\,ds \; . $$  
Since by the unit pseudo-norm relation we have $\;{b^2_s\over
a^2_s} < r_s^2/(1-{R\over r_s})\;$,  whence $\;{b^2_s/a^2_s}\;$ bounded on $A$, the above formula and
Lemma \ref{lem.cv} imply the almost sure convergence of $\;{b^2_s/a^2_s}\,$ on $A\,$. Indeed the bounded
variation terms converge, and as $\,{b^2_s/a^2_s}\,$ is positive, the martingale part has to converge also.
Using the unit pseudo-norm relation again, we deduce that 
$\;{\ds {T^2_s\over a^2_s} = 1 - (1-{\ts{R\over r_s}})\Big({b^2_s\over a^2_sr_s^2}+{1\over a^2_s}\Big)}\,$
has also to converge, necessarily to $0$, since otherwise we would have an infinite limit for $T_s\,$,
which is clearly impossible on $A$. The value of the limit of $\,b_s/a_s\,$ follows now directly from
this and from $(i)$.  
\par \medskip 

   $(iii)$ \ {\it We have almost surely on $A\,$ : $\;{\ds\int_0^\ii (r_t-{\ts{3R\over 2}})^2\,
b^{2}_t\,dt <\ii }\,$, and \ $\;{\ds\int_0^\ii T_t^2\,dt<\ii }\,$. 
} \par \smallskip 

   Let us write It\^o's formula for  
$\; Z_s:= (r_s-{\ts{3R\over 2}})\,T_s = \5\,{d\over d_s}\,(r_s-{\ts{3R\over 2}})^2\,$ : 
$$ Z_s = Z_0 + M_s + {\ts{3\s^2\over 4}}(r_s-r_0)(r_s+r_0-3R) +
\int_0^s T_t^2\,dt + \int_0^s (r_t-{\ts{3R\over 2}})^2 b^{2}_t\,{dt\over r^4_t} 
- {\ts{R\over 2}}\int_0^s (r_t-{\ts{3R\over 2}})\,{dt\over r^2_t} \, , $$  
where $\,M_\cdot\,$ is a local martingale having quadratic variation given by : 
$$\; \langle M\rangle_s = \s^2\int_0^s (r_t-{\ts{3R\over 2}})^2\, \Big( 1-{\ts{R\over
r_t}} + T_t^2\Big)\,dt \, .  $$ 

   Note that if $\,\langle M\rangle_\ii = \ii\,$, then by $\,(ii)\,$ above 
$\;{\ds\lim_{s\to\ii}\limits\; \int_0^s (r_t-{\ts{3R\over 2}})^2\, b^{2}_t\,{dt\over r^4_t}\Big/\langle
M\rangle_s =\ii}\,$. \parn  
Note moreover that in this case 
$\;{\ds \int_0^s |r_t-{\ts{3R\over 2}}|\,{dt\over r^2_t}\,\le  
\sqrt{\int_0^s (r_t-{\ts{3R\over 2}})^2\, b^{2}_t\,{dt\over r^4_t}}\times\sqrt{\int_0^s{dt\over b^2_t}}
}\;$ is also negligible with respect to $\;{\ds\int_0^s (r_t-{\ts{3R\over 2}})^2\,
b^{2}_t\,{dt\over r^4_t} }\;$. \par 
   On the other hand, we must have 
$\;{\ds\liminf_{s\to\ii}\limits\; |Z_s|= 0}\;$ on $A\,$. \par

   Therefore we deduce from Lemma \ref{lem.mart} that necessarily $\,\langle M\rangle_\ii < \ii\,$, and
then that $\,M_s\,$ has to converge, almost surely on $A\,$. \par 
   Using again that $\;{\ds\liminf_{s\to\ii}\limits\; |Z_s|= 0}\,$, we deduce the almost sure
boundedness and convergence on $A\,$ of $\;{\ds\int_0^\ii T_t^2\,dt}\;$ and of $\;{\ds\int_0^\ii
(r_t-{\ts{3R\over 2}})^2\, b^{2}_t\,{dt\over r^4_t} }\;$.   
\par\medskip 

   $(iv)$ \ {\it $(r_s-{\ts{3R\over 2}})^2 b_s\,$ and  $\,{T^2_s/b_s}\;$ go to 0 as $s\to\ii\,$, 
almost surely on $\,A\,$.} \par\smallskip  
   Indeed, on one hand we deduce from $(iii)$ that (for some real Brownian motion $\,W_\cdot$)
$$ (r_s-{\ts{3R\over 2}})^2 b_s= \s W\Big[\int_0^s (r_t-{\ts{3R\over 2}})^4(b_t^2+r_t^2)dt\Big] + 
2\!\int_0^s (r_t-{\ts{3R\over 2}})T_tb_tdt +{\ts{\s^2\over 2}}\!\int_0^s (r_t-{\ts{3R\over 2}})^2
(3b_t+{\ts{r_t^2\over b_t}})dt $$ 
has to converge almost surely on $\,A\,$ as $s\to\ii\,$, necessarily to 0 since it is integrable
with respect to $s\,$.  \par 
  On the other hand we have for some real Brownian motion $\,W'_\cdot$, by It™ formula : 
\parn 
\vbox{
$$ {T^2_s\over b_s} = {T^2_0\over b_0} + \s\, W'\Big[\int_0^s \Big({T^4_t\over b^2_t}+r_t^2\,{T^4_t\over
b_t^4}+4(1-{\ts{R\over r_t}}){T^2_t\over b^2_t}\Big) dt\Big] + {\ts{\s^2\over 2}}\int_0^s {T^2_t\over
b_t}\,dt + 2\int_0^s (r_t-{\ts{3R\over 2}})\,T_t\,b_t\,{dt\over r^4_t} \hskip 6mm  $$ 
$$ \hskip 50mm   - \int_0^s {R\,T_t\over r^2_t\,b_t}\,dt + \s^2\!\int_0^s(1-{\ts{R\over r_t}})\,{dt\over
b_t} +  {\ts{\s^2\over 2}}\int_0^s {r_t^2\,T_t^2\over b_t^3}\,dt \; . $$ }
Recall from $(i)$ that $\,{T_t\over b_t}\to 0\,$ and that $\,b_t\sim {3R\sqrt{3}\over 2}\,a_t\,$. 
Thus using $(iii)$ we see easily that all integrals in the above formula converge. Hence we deduce the
almost sure convergence of $\,{\ds s\mapsto {T^2_s/b_s}}\,$ on $\,A\,$, necessarily to $0\,$, since it
is integrable.  \par\medskip 

   $(v)$ \ {\it It is sufficient to show that $\;{\ds\int_0^\ii |r_t-{\ts{3R\over 2}}|\, |T_t|\,
b^{2}_t\,dt <\ii}\;$, and that $\;{\ds \int_0^\ii T_t^4\,dt <\ii}\;$, almost surely on $\,A\,$. }
\par\smallskip  
   Indeed, assuming that these 2 integrals are finite, It™'s formula shows that we have for some real
Brownian motion $\,W''_\cdot$ : \parn 
\vbox{
$$ T_s^2 = T_0^2 + 2\s\, W''\Big[\int_0^s (T_t^2+1-{\ts{R\over r_t}})\,T_t^2\,dt\Big] +
4\s^2\!\int_0^sT_t^2\,dt + 2\!\int_0^s (r_t-{\ts{3R\over 2}})\,T_t\,b^{2}_t\,{dt\over r^4_t}\hskip 10mm $$
\vspace{-3mm} 
$$ \hskip 50mm + \,\s^2\!\int_0^s(1-{\ts{R\over r_t}})\,dt - R\!\int_0^s{T_t\over r_t^2}\,dt\;  $$ }
$$ = \g_s + \s^2\int_0^s(1-{\ts{R\over r_t}})\,dt - R\!\int_0^s{T_t\over r_t^2}\,dt 
= \g_s + \int_0^s\Big[ {\ts{1\over 3}}+ {\ts{2\over 3r_t}}(r_t-{\ts{3R\over 2}})\Big] dt 
+ {R\over r_s}-{R\over r_0}\, = \g'_s + s/3\, , $$  
where $\,\g_\cdot\;,\;\g'_\cdot\,$ (since $\,|r_t-{\ts{3R\over 2}}|=o(b_t^{-1/2})=o(a_t^{-1/2})\,$ by
$(iv)$ and $(ii)$) are bounded converging processes on $\,A\,$. \ Whence 
$\;\lim_{s\to\ii}\limits T^2_s=\ii\;$ almost surely on $A\,$, which with $(iii)$ above implies 
that $\,A\,$ must be negligible.  \par\medskip  

   $(vi)$ \ {\it End of the proof of the convergence of $\,\,r_s\,$ to $\,\ii$ on $\,\{ D=\ii\}\,$. }
\par\smallskip 
   By Schwarz inequality, the first bound in $(v)$ above will follow from 
$\;{\ds\int_0^s T_t^2b_t\,dt<\ii}\;$ and from 
$\;{\ds\int_0^s (r_t-{\ts{3R\over 2}})^2 b^{3}_t\,dt <\ii}\;$. Now these two terms appear in the It™
expression for \ $\;Z_s^1:= (r_s-{\ts{3R\over 2}})\,T_s\,b_s\;$ :  
$$ Z_s^1 = Z_0^1 + M_s^1 + {\ts{\s^2\over 2}}\int_0^s \!\Big[ 8+{r_t^2\over b_t^2}\Big] Z_t^1 dt +
\int_0^s T_t^2b_t\,dt + \int_0^s (r_t-{\ts{3R\over 2}})^2 b^{3}_t\,{dt\over r^4_t} 
- {\ts{R\over 2}}\!\int_0^s (r_t-{\ts{3R\over 2}})b_t\,{dt\over r^2_t} \;, $$ 
with a local martingale $\,M^1_\cdot\,$ having quadratic variation :  
$$ \langle M^1\rangle_s = \s^2\int_0^s (r_t-{\ts{3R\over 2}})^2\,b^{2}_t\times \Big( 1-{\ts{R\over r_t}}
+[4+ r_t^2b_t\2]\,T_t^2\Big)\,dt\; .  $$ 

      Note that by Schwarz inequality, $(iii)$ above implies that $\;{\ds\int_0^\ii
|Z^1_t|\,dt <\ii }\;$, and then that $\;{\ds \int_0^s \Big[ 8+{r_t^2\over b_t^2}\Big] Z_t^1 dt}\;$ is 
bounded and converges, almost surely on $\,A\,$, as $s\to\ii\,$. \par 
   Using the first assertion of $(iv)$, observe that 
$\;{\ds\lim_{s\to\ii}\limits\; {{\ds\int_0^s} (r_t-{\ts{3R\over 2}})^2\, b^{3}_t\,{dt\over r^4_t}
+{\ds\int_0^s} T_t^2b_t\,dt\over\langle M^1\rangle_s} =\ii}\;$ if $\,\langle M^1\rangle_\ii = \ii\,$.
\quad  Note moreover that in this case 
$$ \Big|\int_0^s (r_t-{\ts{3R\over 2}})\,b_t\,{dt\over r^2_t}\,\Big| \le 
\sqrt{\int_0^s (r_t-{\ts{3R\over 2}})^2\, b^{3}_t\,{dt\over r^4_t}}\times\sqrt{\int_0^s{dt\over b_t}} $$ 
is also negligible with respect to $\;{\ds\int_0^s (r_t-{\ts{3R\over 2}})^2\,
b^{3}_t\,{dt\over r^4_t}+{\ds\int_0^s} T_t^2\,b_t\,dt }\;$. \par 
   Therefore we deduce from Lemma \ref{lem.mart} and from the integrability of $\,t\mapsto|Z^1_t|\,$, that
necessarily $\,\langle M^1\rangle_\ii < \ii\,$, and then that $\,M^1_s\,$ has to converge, almost surely
on $A\,$. \par 
   Hence $\,Z^1_\cdot\,$ must have a limit almost surely on $A\,$, which must be 0, owing to the 
integrability of $\,Z^1_\cdot\,$. This forces clearly  
$\;{\ds\int_0^\ii (r_t-{\ts{3R\over 2}})^2\, b^{3}_t\,{dt\over r^4_t}+{\ds\int_0^\ii} T_t^2\,b_t\,dt }\;$
to be finite, almost surely on $A\,$, showing the first bound in $(v)$ above. \par \smallskip 

   Finally, the integrability of $\,T_t^2\,b_t\,\,$ and the second convergence of $(iv)$ imply 
the second bound in $(v)$ above : $\;{\ds \int_0^\ii T^4_t\,dt <\ii}\;$ almost surely on $\,A\,$. \par 
   This concludes the proof of the first assertion in Theorem \ref{the.asymp}.
\par\bigskip 

   2) \ {\it $r_s\to R\,$ and $\,r_s\to\ii\,$ occur both with positive probability, from any initial
condition. }  \par\smallskip    

   Let us use the support theorem of Stroock and Varadhan (see for example ([I-W], Theorem VI.8.1)) to show
that the diffusion $\,(r_\cdot,b_\cdot,T_\cdot)\,$ of Corollary \ref{cor.diff1} is irreducible.  
Since we can decompose further the equations given in Proposition \ref{pr.diff1} for
$\,(r_\cdot,b_\cdot,T_\cdot)\,$, using a standard Brownian motion $\,(w_\cdot,\b_\cdot,\g_\cdot)\in\R^3\,$,
as follows :
$$ dr_s = T_s\, ds\;,\quad d b_s = \s\,b_s\,dw_s + \s\,r_s\,d\b_s + {\ts{3\,\s^2\over 2}}\, b_s\,ds + 
{\s^2\, r_s^2\over 2\, b_s}\,ds\; , $$ 
$$ dT_s =  \s\,T_s\,dw_s + \s\,\sqrt{1-{\ts{R\over r_s}}}\,d\g_s + {\ts{3\,\s^2\over 2}}\, T_s\, ds + 
(r_s-{\ts{3\over 2}}R)\,{b_s^2\over r_s^4}\, ds - {R\over 2r_s^2}\, ds \; , $$ 
we see that trajectories moving the coordinate $\,b_\cdot\,$ without changing the others, and trajectories
moving the coordinate $\,T_\cdot\,$ without changing the others, belong to the support of 
$\,(r_\cdot,b_\cdot,T_\cdot)\,$. Moreover we see from Section \ref{sec.geod} that there are
timelike geodesics,  and then trajectories in the support, which link $\,r\,$ to $\,r'\,$, and then
considering the velocities also, which link say  $\,(r,b'',T'')\,$ to $\,(r',b'',T''')\,$. \ So, for given
$\,(r,b,T)\,$ and $\,(r',b',T')\,$ in the state space, we can, within the support of
$\,(r_\cdot,b_\cdot,T_\cdot)\,$, move
$\,(r,b,T)\,$ to $\,(r,b'',T'')\,$, then $\,(r,b'',T'')\,$ to $\,(r',b'',T''')\,$, and finally move 
$\,(r',b'',T''')\,$ to $\,(r',b',T')\,$, thereby showing the irreducibility of
$\,(r_\cdot,b_\cdot,T_\cdot)\,$. 
\par 

   This implies that it is enough to show that for large enough $\,r_0,T_0\,$, the convergence to
$\ii$ occurs with probability $\ge 1/2\,$, and that for $\,r_0\,$ close enough from $R$ and $\,T_0\,$
negative enough, the convergence to $R$ occurs with probability $\ge 1/2\,$ as well. Now this can be done
by a classical supermartingale argument using the process $\,1/|T_s|\,$, stopped at some hitting time.
Indeed we see from Proposition \ref{pr.diff1}  that 
$$ {1\over |T_s|} + \int_0^s \Big( {\ts{\s^2\over 2}}\,T_t^2 - \s^2\,(1-{\ts{R\over r_s}}) 
- {R\,T_t\over 2\,r_t^2}\, + (2r_t-3R)\,{b_t^2\,T_t\over r_t^4}\,\Big)\, {dt\over |T_t|^3} $$  
is a local martingale. \par 
   Take first $\;r_0\ge 3R/2\,$, $\,T_0\ge 4+{4\over R\s^2}\,$, and $\;\tau :=\inf\Big\{ s>0\,\Big|\,
T_s=2+{2\over R\s^2}\,\Big\}$ : \ $r_s$ increases on $\,\{0\le s<\tau\}\,$ and then we see that 
$\,1/|T_{s\wedge\tau}|\,$ is a supermartingale, which implies that \parn 
$(2+{2\over R\s^2})\1\P (\tau <\ii )\le \liminf_{s\to\ii}\limits\, \E ({1\over |T_{s\wedge\tau}|}
1_{\{\tau <\ii\}}) \le \liminf_{s\to\ii}\limits\, 
\E ({1\over |T_{s\wedge\tau}|}) \le \E ({1\over T_{0}}) \le (4+{4\over R\s^2})\1\,$, \parn 
and then that $\;\P (\lim_{s\to\ii}\limits r_s = +\ii ) \ge \P (\tau =\ii ) \ge 1/2\,$.  \par 

   Conversely take $\;r_0\le 3R/2\,$, $\,T_0\le -{2}\,$, and $\;\tau':=\inf\Big\{ s>0\,\Big|\,
T_s=-\sqrt{2}\,\Big\}$ : \ $r_s$ decreases on $\,\{0\le s<\tau'\}\,$ and then we see that 
$\,1/|T_{s\wedge\tau'}|\,$ is a supermartingale, which implies that \parn 
\centerline{$2^{-1/2}\,\P (\tau'<\ii )\le \liminf_{s\to\ii}\limits\, \E ({1\over |T_{s\wedge\tau'}|}
\,1_{\{\tau'<\ii\}}) \le \liminf_{s\to\ii}\limits\, 
\E ({1\over |T_{s\wedge\tau'}|}) \le \E ({1\over |T_{0}|}) \le {1/2}\,$,} \parn 
and then that $\;\P (D<\ii ) \ge \P (\tau'=\ii ) \ge 1/\sqrt{2}\,$. \par 
   This concludes the proof of the second assertion in Theorem \ref{the.asymp}.\par\bigskip   

    3) {\it Existence of an asymptotic direction for the relativistic diffusion, on $\,\{ D=\ii\}$.} 
\par\smallskip  

   We want to generalize the observation made in Section \ref{sec.relr} for $R=0$, see 
\mbox{Remark \ref{pro.relr}}. Recall from Lemma \ref{lem.cv} that it does not matter for this
asymptotic behaviour whether we consider the trajectories as function of $\,s\,$ or of $t(s)\,$ (id est as
viewed from a fixed point). \par   
   We shall use Remark \ref{pro.relr} and  proceed by comparison between the flat Minkowski case $R=0$ and the Schwarzschild case $R>0$. \    Let us split this proof into four parts. \par \medskip  

 $(i)$ \ {\it  We have $\;{\ds\int_0^\ii {a_t\over r^2_t}\,dt <\ii}\;$ and $\;{\ds\int_0^\ii {U_t\over
r_t}\,dt <\ii}\;$, almost surely on $\,\{ D=\ii\}$.}\par\smallskip 

     We know from 1) above that $\,r_s\to\ii\,$ almost surely on $\,\{ D=\ii\}$. \par 

   The very beginning of this proof remains valid : Using (1,$i$) again, we have almost surely 
$\;{\ds\int_0^\ii {T^2_t\over a_t\,r^2_t}\,dt}\,$ and 
$\,{\ds\int_0^\ii (1-{\ts{3R\over 2r_t}})^2\, {b^2_t\over a_t\,r^3_t}\,dt}\;$ finite, whence 
$\;{\ds\int_0^\ii {b^2_t\over a_t\,r^3_t}\,dt}\;$ finite, and then, since 
$\; {\ds {a\over r^2}\le {T^2\over a\,r^2}+{b^2\over a\,r^4}+{1\over a\,r^2} }\;$, 
also $\;{\ds\int_0^\ii {a_t\over r^2_t}\,dt}\;$ finite, almost surely on $\,\{ D=\ii\}$. \par 
   Now by the unit pseudo-norm relation, we have 
$\;{\ds {U\over r}  = {b\over r^3} \le {a\over r^2\,\sqrt{1-{\ts{R\over r}}}}}\;$,  \ whence \parn 
$\;{\ds\int_0^\ii {U_t\over r_t}\,dt}\;$ finite, almost surely on $\,\{ D=\ii\}$.
\par\medskip 

$(ii)$ \ {\it  The perturbation of the Christoffel symbols due to $R$ is $\,\O(r\2)\,$.}
\par\smallskip 

   Recall from the beginning of Section \ref{sec.S} the values of the Christoffel symbols $\,\G^i_{jk}\,$. 
Denote by $\,\tilde\G^i_{jk}\,$ the difference between these symbols and their analogues for $R=0$, which
is a tensor, has only five non-vanishing components in spherical coordinates, and then is easily
computed in Euclidian coordinates
$\;(x_1=r\,\sin\f\,\cos\psi\,;\,x_2=r\,\sin\f\,\sin\psi\,;\,x_3=r\,\cos\f)\,$ : \ we find 
$$ \tilde\G^{x_i}_{x_j,x_k} = {\p x_i\over \p r}\times \Big( {\p r\over \p x_j}{\p r\over \p x_k}\G^r_{rr}
+ {\p \f\over \p x_j}{\p \f\over \p x_k}\G^r_{\f\f} + {\p \psi\over \p x_j}{\p \psi\over \p x_k}
\G^r_{\psi\psi}\Big) $$ 
$$ = {x_i\over r}\times \Big( {-R\over 2r(r-R)}{x_j\over r}{x_k\over r} + R{\p \f\over \p x_j} 
{\p \f\over\p x_k} + R\sin^2\f\,{\p \psi\over \p x_j}{\p \psi\over \p x_k}\Big) = \O(r\2)  $$ 
since $\;\Big|{\p\f\over\p x_j}\Big|\le 1/r\;$ and  $\;\Big|{\p\psi\over\p x_j}\Big|\le 1/(r\sin\f)\;$. 
The same is valid directly for the remaining components $\;\tilde\G^{x_i}_{t,t}\,$ and 
$\,\tilde\G^{t}_{x_j,t}$.  \par \medskip 

 $(iii)$ \ {\it The parallel transport converges, almost surely on $\,\{ D=\ii\}$.}
\par\smallskip 

   Recall from Theorem \ref{the.gen} (in Section \ref{sec.relg}) that the inverse parallel
transport $\,\overleftarrow{\xi}(s)\,$ along the $C^1$ curve $\,(\xi_{s'}\,|\,0\le s'\le s)$ satisfies 
$$ {\ts{d\over ds}}\,\overleftarrow{\xi}(s)^i_j \, =\, \overleftarrow{\xi}(s)^i_k\times\G_{j\ell}^{k}(\xi_s)
\times\dot\xi_s^\ell\; , $$ 
so that, using $(ii)$ above and $\,|T_s|\le a_s\,$ : 
$$ \overleftarrow{\xi}(s)\, =\, \int_0^s \O\Big( r_v\2\times |\dot\xi_v|\Big)\,dv = 
\int_0^s \O\Big( |\dot t_v| + |\dot r_v| + r_v\,|\dot\t_v|\Big)\,r_v\2\,dv 
= \int_0^s\O\Big( 2\,\,{a_v\over r_v^2}+{U_v\over r_v}\Big)\, dv\, . $$ 
Hence, using $(i)$ above, we conclude that the parallel transport (as its inverse $\,\overleftarrow{\xi}(s)$) admits a finite limit as $s\to\ii\,$, almost surely on $\,\{ D=\ii\}$.
\par\medskip 

 $(iv)$ \ {\it End of the proof.} \par\smallskip 

   Recall from Theorem \ref{the.gen} (in Section \ref{sec.relg}) that  the continuous process $(\,\zeta_s = \overleftarrow{\xi}(s)\,\dot\xi_s\,)$ is a hyperbolic Brownian motion on the hyperbolic space $\,T_{\xi_0}^1\M$, isometric to $\,\H^3$. \par 
   Now, according to Section \ref{sec.relr}, where merely $\,a_s=p^o_s\,$, we know that $\;\zeta_s/a_s\;$ converges almost surely as $\,s\to\ii\,$ towards $\,(1,1,\t_\ii)\,$ (in coordinates $(t,r,\t)$), for some random $\,\t_\ii\in\S^{2}$. \  Using $(iii)$ above, we deduce that $\;\dot\xi_s/a_s\;$ converges almost surely as $\,s\to\ii\,$ towards $(1,1,\hat\t_\ii)$, for some random $\,\hat\t_\ii\in\S^{2}$. \ 
This means also that the velocity $\,dZ_t/dt\,$ of the trajectory $\;Z_\cdot := (t_s\mapsto (r_s,\t_s))$ 
converges almost surely towards $\,(1,\hat\t_\ii)$, 1 being here the velocity of light. \  Merely integrating this, as in the proof of Remark \ref{pro.relr}, we get finally the  generalization of Remark \ref{pro.relr} to the relativistic diffusion. This ends the whole proof of Theorem \ref{the.asymp}. \ $\;\diamond $ 
   \par\medskip 
   
   Note moreover that the part (3,$(iv)$) of the above proof shows that, since the hyperbolic Brownian motion does not explode, there is no explosion at the level of the fibre. This is a general fact for relativistic diffusions.  Finally if the radius $\,r_s\,$ could explode within some finite proper time, $\,T_s\,$, and by the unit pseudo-norm relation the energy $\,a_s\,$, would explode as well ; but this is clearly forbidden by the simple stochastic differential equation governing $\,a_s\,$. This proves the non-explosion of the Schwarzschild relativistic diffusion, quoted directly after the statement of Theorem \ref{the.asymp}. 

\subsection{Proof of Theorem \ref{the.bh}}\indf \label{sec.PT3}
     $(i)$ \  The first sentence is clear from Theorem \ref{the.asymp} when the diffusion starts above $\,{\cal S}_0\,$, except the finiteness of $\,D'$, proved in $(ii)$ below. This is the same when the diffusion starts above $\,-{\cal S}_0\,$, Theorem \ref{the.asymp} being valid as well in this very similar case (it is sufficient to change the signs of $\,a\,$ and $\,t$).  The other cases are reviewed in Remark \ref{rem.-S}. \ So we just have to establish the assertions $(1), (2), (3)$ of the statement,
assuming $\,D\,$ finite. \par \smallskip 

    $(ii)$ \ Let us prove $(1)$ first. Since $\,T_D<0\,$ and $\,T_s\,$ cannot vanish in the region 
$\,\{ r<R\}$, $\,r_s\,$ must decrease and then converge to some $\,\rr\in[0,R[\,$, with 
$\;{\ds\limsup T_\cdot \le -\sqrt{{\ts{R\over\rr}}-1}\,<0}\;$, which in turn forces $\,\rr=0\,$, 
hit within a finite proper time $\,D'\,$, and $\,\lim_{s\nea D'}\,T_s = -\ii\,$.

   Only the upper bound for $D'$ remains to be proved. It is a consequence of the pseudo-norm
equation, which implies $\,T^2\ge {R\over r}-1\,$ and then $\;T_s\le -\sqrt{{R\over r_s}-1}\,$ on 
$\{ D\le s\le D'\}\,$. Indeed, consider 
$$ g(r):= {\ts{R\over 2}}\,{\rm Arcsin}\Big[{\ts{2\over R}}\sqrt{r(R-r)}\,\Big] - \sqrt{r(R-r)}
\quad \hbox{ if }\quad 0\le r\le {\ts{R\over 2}}\; , $$ 
$$ \hbox{and }\qquad 
g(r):= {\ts{R\over 2}}\Big(\pi -{\rm Arcsin}\Big[{\ts{2\over R}}\sqrt{r(R-r)}\,\Big]\Big) -
\sqrt{r(R-r)} \quad \hbox { if }\quad {\ts{R\over 2}}\le r\le R\; . \quad $$ 
We have \quad $ g'(r) = ({R\over r}-1)^{-1/2}\,$, $\; g''(r) = {R\over 2r^2}\times({R\over r}-1)^{-3/2}\,$,
$\; g(0)=g'(0)=0\;$, \par $\; g(R)={\ts{\pi\over 2}}\,R\;,\; g'(R)=+\ii\,$. 
\ Hence $\; g'(r_s)\times T_s \le -1\;$ implies by integration on $[D,s]\,$, for $ D\le s\le D'$ : \quad 
$ g(r_s)\le {\ts{\pi\over 2}}\,R +D-s\,$ ; taking $\,s=D'\,$, this proves $(1)$. \par \smallskip 

   $(iii)$ \ Let us then prove the non-explosion of $\,\log b_\cdot\,$ at $\,D'$. \ 
For that, let us apply the comparison theorem  (see ([I-W], VI, th 4.1)) to $\,b_{[D,D'[}\,$ : we get so a real diffusion process $\,\b_\cdot\,$ solving 
$$ \b_s = b_D + \s \int_D^s \sqrt{\b_s^2+R^2}\,dw_s +{\ts{ 3\s^2\over 2}} \int_D^s (\b_s^2+R^2)\, {ds\over \b_s} \, , $$ 
and such that almost surely\quad ${\ds\sup_{D\le s<D'} b_s\le \sup_{D\le s<D'}\b_s}\,$. \ Indeed the ratio of the drift coefficient and of the squared diffusion coefficient of $\,b_\cdot\,$ is maximal for $\,r_s=0\,$. \ Now there is a real Brownian motion $\,B_\cdot\,$ such that \ 
$$\b_{D+s'} = B^{-1/2} \Big[ \inf\Big\{ s \,\Big|\, \int_0^s {dt\over (1+R^2B_t)B_t^2} > 4 s'\Big\}\Big] , $$
so that $\,\b_\cdot\,$ cannot diverge at a finite time ; indeed it is immediately seen that  \parn  
${\ds \int_0^{s} {dt\over B_t^2} = 2\log (B_0/B_s) + 2 \int_0^s {dB_t\over B_t}}\,$ must almost surely diverge as $\,s\,$ approaches the hitting time of 0 by $\,B\,$. \  
   This proves that $\,b_\cdot\,$ almost surely cannot explode at $\,D'\,$, and thus must be continuous on $\,[0,D']\,$. \par 
    Moreover, since the differential equation governing $\,b_\cdot\,$ can be written \parn 
${\ds\log (b_s/b_0) = \s\,W\Big[\int_0^s \Big( 1+{r_t^2\over b_t^2}\Big)\,dt\,\Big] +\s^2\, s}\;$, \ 
for some standard Brownian motion $\,W_\cdot\,$, we see that $\,b_{D'}=0\,$ would imply $\;{\ds\int^{D'}_0 {r_t^2\over b_t^2}\,dt = \ii}\,$, and then $\,\limsup_{s\nea D'}\limits\, b_s = +\ii\,$, a contradiction. \par

   $(iv)$ \ Let us now prove the statement (3). For that, let us write again the unit pseudo-norm relation : near $D'$ it writes \par\smallskip   
$\, T_s\,r_s^{3/2} = - \sqrt{(R-r_s)b_s^2+Rr_s^2+(a_s^2-1)r_s^3}\;\mathop{\lra}_{s\nea D'}\limits\,
-b_{D'}\sqrt{R}\:$. \ So that for any $\,\e>0\,$  we have : \quad
$-b_{D'}\sqrt{R}\,(1+\e)\le T_s\,r_s^{3/2}\le -b_{D'}\sqrt{R}\,(1-\e)\;$, \  for $\,s\,$ sufficiently
close to $D'$. \parn 
Integrating this, we get immediately, \ for $\,s\,$ sufficiently close to $D'$ : 
$$ -b_{D'}\sqrt{R}\,(1+\e)(D'-s)\le -{\ts{2\over 5}}\,r_s^{5/2}\le -b_{D'}\sqrt{R}\,(1-\e)(D'-s)\, , $$ 
which means the first equivalent in the statement $(3)$. The second equivalent follows at once by using again the unit pseudo-norm relation. \ This proves $(3)$. \par\smallskip  
   As a consequence, we deduce at once that $\,\,|\dot\t_s|=U_s=b_s/r_s^2\,\,$ is integrable near $D'$, which implies the convergence of $\,\t_s\,$ in $\S^2$ as  $\,s\nea D'$. \par \smallskip  
   
   $(v)$ \ End of the proof of the statement (2). \par 
   To establish the non-explosion of $a_s\,$, let us rewrite its stochastic differential equation
with two independent real standard Brownian motions $\, w,\b\,$ : \par 
$\;d a_s = \s\,a_s\,dw_s + \s\,({\ts{R\over r_s}}-1)\,d\b_s + {\ts{3\,\s^2\over 2}}\,a_s\,ds\,$, \ 
and consider $\, X_s:=a_s\times e^{-\s\,w_{s-D}-\s^2(s-D)}$. \parn 
 We have $\; d X_s= e^{-\s\,w_{s-D}-\s^2(s-D)}\,\s\,({\ts{R\over r_s}}-1)\,d\b_s\,$, \ whence for $\,D\le s\le D'\,$ : \parn   
$$ a_s = e^{\s\,w_{s-D}+\s^2(s-D)}\; \Big( a_D + \s\int_D^s e^{-\s\,w_{t-D}-\s^2(t-D)}\,
({\ts{R\over r_t}}-1)\,d\b_t\Big) \qquad $$
$$ \qquad = e^{\s\,w_{s-D}+\s^2(s-D)} \Big( a_D + \s\, W\Big[\int_D^s e^{-2\s\,w_{t-D}-2\s^2(t-D)}\, 
({\ts{R\over r_t}}-1)^2\,dt\Big]\Big)\,, $$ 
for some real standard Brownian motion $\,W_\cdot$. \ The equivalent seen above for $\,r_s\,$ near $D'$ shows the almost sure convergence of the integral $\,\int_D^{D'} {dt\over r_t^2}\,$, and then 
by the above formula, of $\,a_s\,$ as $\nea D'$. \par 
   The same equivalent again shows the almost sure integrability of $\,\dot u^-_s\,$ near $D'$ : indeed 
$$ \dot u^-_s = (a_s+T_s)/(1-R/r_s) \sim [2\,b_{D'}^4/(5\,R^3\,(D'-s))]^{-1/5}\, , $$ 
which finally proves the convergence of $\,u^-_s\,$ as $\nea D'$. Likewise for $\,u^+_s\,$. \par
   Finally it remains to show the convergence of $\,\bv_s\,$, or equivalently of
$\,\bv_s/b_s = \t_s\wedge\dot\t_s/U_s\,$, since we already saw above the convergence of $\,b_s\,$ in $\R_+^*$. Since we also saw the convergence of $\,\t_s\,$ in $\S^2$, it remains to get the convergence of $\,\dot\t_s/U_s\,$. Now this is the last assertion of Proposition \ref{pro.systeqst}. $\;\diamond$ 

\subsection{Proof of Proposition \ref{pro.systeqst}}\indf \label{sec.PP3}
   We have 
$$ d\Big({\dot\t_s\over U_s}\Big) 
= 2\,r_s\,b_s\1\,T_s\,\dot\t_s + r_s^2\,b_s\1\,d\dot\t_s - (r_s^2\,b_s\2\,db_s)\,\dot\t_s + 
r_s^2\,b_s^{-3}\,\langle db_s\rangle\,\dot\t_s - r_s^2\,b_s\2\,\langle db_s,d\dot\t_s\rangle\; . $$ 
\indent   To perform the computations, let us use the basis $\,(u,v,k)\,$ of $\,\R^3\,$ defined by : 
$$ u=(\cos\psi, \sin\psi , 0)\;;\quad v=(-\sin\psi, \cos\psi , 0)\;;\quad k=(0,0,1)\; , $$ 
so that \quad $\t = u\sin\f +k\cos\f \;$, \ and  
$$  \dot\t = (u\cos\f -k\sin\f)\,\dot\f + v\,\dot\psi\sin\f \quad ;\quad  
\t\wedge\dot\t = v\,\dot\f - (u\cos\f -k\sin\f)\,\dot\psi\sin\f \; . $$ 
Then \quad 
${\ds  d\dot\t = (u\cos\f -k\sin\f)\,d\dot\f + v\,\sin\f\,d\dot\psi - \dot\f^2\t
+\,2\,v\,\dot\f\dot\psi\,\cos\f - u\,\dot\psi^2\,\sin\f\; , }$ \quad 
$$ {\rm whence }\hskip 12mm \langle d\dot\t, db\rangle = (u\cos\f -k\sin\f)\,\langle d\dot\f, db\rangle 
+ v\,\sin\f\, \langle d\dot\psi, db\rangle \hskip 40mm $$
$$ \hskip 24mm = \Big( (u\cos\f -k\sin\f)\dot\f + v\,\sin\f\,\dot\psi\Big)\,\s^2\,b\1\,(b^2+r^2) \, =\; 
 b\1\,\langle db\rangle\,\dot\t\; , $$ 
so that the last two terms in the expression of $\, d\Big({\dot\t_s\over U_s}\Big)\,$ above cancel. 
\quad Hence \par
\vbox{  
$$ d\Big({\dot\t_s\over U_s}\Big)  = \; 2\,{{r_s\over b_s}}\,T_s\,\dot\t_s -
r_s^2\,b_s\2\,\Big(\s\,dM^b_s + {\ts{3\s^2\over 2}}\,b_s\,ds +{\ts{\s^2\,r_s^2\over 2\,b_s}}\,ds\Big)
\,\dot\t_s  \hskip 44mm $$ \vspace{-3mm} 
$$  \hskip 8mm +\, U_s\1\,(u_s\cos\f_s -k\,\sin\f_s)\, 
\Big( \s\,dM^\f_s+{\ts{3\s^2\over 2}}\,\dot \f_s\,ds - {\ts{2\over r_s}}\,T_s\,\dot
\f_s\,ds  + \sin\f_s\,\cos\f_s\;\dot\psi^2_s\,ds \Big) $$ 
$$ +\, U_s\1\,v_s\,\sin\f_s\,\Big( \s\,dM^\psi_s+{\ts{3\s^2\over 2}}\,\dot \psi_s\,ds -
{\ts{2\over r_s}}\,T_s\,\dot \psi_s\,ds  - 2\, \cotg\f_s\;\dot\f_s\dot\psi_s\,ds \Big) $$ 
$$ +\, U_s\1\Big( 2\,v_s\,\dot\f_s\dot\psi_s\,\cos\f_s - \dot\f_s^2\,\t_s
- u_s\,\dot\psi^2_s\,\sin\f_s\Big) ds  $$ 
} \vspace{-7mm} 
$$ =\; {{\s\,r_s^2\over b_s^2}}\,\Big( -\dot\t_s\,dM^b_s + (u_s\cos\f_s -k\,\sin\f_s)\, b_s\,dM^\f_s +
v_s\,b_s\,\sin\f_s\,dM^\psi_s\Big) \qquad $$ 
$$ \qquad -\,{{\s^2\,r_s^4\over 2\,b_s^3}}\,\dot\t_s\,ds - U_s\1\,\dot\f_s^2\,\t_s\,ds +
U_s\1\,\sin\f_s\,\Big( (u_s\cos\f_s -k\,\sin\f_s)\,\cos\f_s - u_s\Big) \dot\psi^2_s\,ds \,\, . $$ 
Observe now that the definition of $\,b=r^2\,U\,$ implies 
$$ dM^b_s =\,r_s^2\,U_s\1\,(\dot\f_s\,dM^\f_s + \sin^2\f_s\,\dot\psi_s\,dM^\psi_s)\; . $$ 
Therefore, expressing all in the basis $\,(\t_s,\,\dot\t_s/U_s,\,\t_s\wedge\dot\t_s/U_s)\,$, we get 
$$ d\Big({\dot\t_s\over U_s}\Big) =\; \s\,U_s\2\,\sin\f_s\,\Big( \dot\f_s\,dM^\psi_s - 
\dot\psi_s\,dM^\f_s\Big)\, \t_s\wedge{\dot\t_s\over U_s} - {{\s^2\,r_s^4\over 2\,b_s^3}}\,\dot\t_s\,ds - 
U_s\,\t_s\,ds \; , $$ 
that is to say the formula of the statement, with 
$\; d\b_s:= r_s\,U_s\1\sin\f_s\,\Big( \dot\f_s\,dM^\psi_s - \dot\psi_s\,dM^\f_s\Big)$. \par 

   Now it is straightforward to verify (from the covariation matrix given about the beginning of Section
\ref{sec.bhh}, before Theorem \ref{the.bh}) that 
$$ \langle d\b_s\rangle\, = ds\; ,\quad \langle d\b_s, dM_s^b\rangle\, = 0\; ,\quad \langle d\b_s,
dM_s^T\rangle\, = 0\; , $$ 
which shows that indeed $\,\b_\cdot\,$ is a standard Brownian motion and is independent from $\,(r_\cdot,
b_\cdot,T_\cdot)$. \par 
   To establish the last assertion of the statement, we deduce from this expression for  
$\,d\Big({\dot\t_s\over U_s}\Big)\,$ that the bounded variation part of $\;{\dot\t_s/U_s}\;$ is not
larger than $\;{b_s\,r\2_s} +{\s^2\,b_s\2\,r_s^2}\,$, which (as $\,r_s\2$, recall
Theorem (\ref{cor.bh},3)) is almost surely integrable on $\,[0,D']\,$, while its martingale part has
quadratic variation not larger than ${\s^2\,b_s\2\,r_s^2}\,$, and thus is almost surely integrable on
$\,[0,D']\,$ as well by Theorem (\ref{cor.bh},2). $\;\diamond$ 

\subsection{Proof of Proposition \ref{the.martp}}\indf \label{sec.PT4}
  1) \ To establish this, we need a system of stochastic
differential equations relative to the whole relativistic diffusion $\,(\xi_s,\dot\xi_s)$, and not only
to its projection $\,(r_s, b_s, T_s)\,$ on the coordinates $\,(r, b,T)$. Recall that the whole
relativistic diffusion $\,(\xi_s,\dot\xi_s)\,$ lives in a 7-dimensional space. For our purpose the
following system of coordinates is convenient : $\,(r, \t, u^-, a, T, \dot\t)\in
\R_+\times\S^2\times\R\times\R\times\R\times T_{\t}\S^2\,$. \ Recall that \parn 
\centerline{$ b= r^2\,|\dot\t|= r^2\,U = r^2\,\sqrt{\dot\f^2+\dot\psi^2\,\sin^2\f}\,$.}\par 

   Now we have the following simplification : we can recover the coordinate $\,u^-\,$ from the other
ones by integration, since its value is fixed at the origin. Hence, to recover the whole relativistic diffusion $\,(\xi_s,\dot\xi_s)$, it is sufficient to get the reduced diffusion $\,(r_s, a_s, b_s, T_s)$, and to recover $\,(\t_s, \dot\t_s)\,$ conditionally from $\,(r_s, a_s, b_s, T_s)$.
\par\smallskip 

   2) \ Let us deal with the reduced diffusion $(r_s, a_s,
b_s)\,$ : we have to show first that for any $\,(a_0,b_0)\in\R\times\R_+^*\,$ there exists a unique
(in law) diffusion process $\,(r_s, a_s, b_s)\,\in [0,R]\times\R\times\R_+^*\,$ defined up to proper
time $\,D'':=\inf\{s\,|\,r_s=R\}\,$, starting with initial condition $\,(0,a_0,b_0)\,$, and having
infinitesimal generator (see Proposition \ref{pr.diff1}) 
$$ \GG'' := T{\p\over\p r} + {\s^2\over 2}\,\Big( 
(a^2-1+{\ts{R\over r}})\,{\p^2\over\p a^2} + (b^2 + r^2) \,{\p^2\over\p b^2} + 2a b\,{\p^2\over\p a\p
b} + 3a\,{\p\over\p a} + (3b + {r^2\over b})\,{\p\over\p b}\Big) \, , $$  
where $\quad T= T(r,a,b) := \sqrt{a^2+({\ts{R\over r}}-1)({\ts{b^2\over r^2}}+1)}\,\,$ is
chosen positive. \par 
   Clearly $\,r_s\,$ must increase strictly as $\,s\nea D''\,$, so that $\,s\,$ can be expressed as random continuous strictly increasing function of $\,r\in [0,R]\,$ : $\,s=s(r)\,$. Let us set 
$\,\tilde{a}_r:= a_{s(r)}\,$ and $\,\tilde{b}_r:= b_{s(r)}\,$. In other words, we consider here the 
radial coordinate $r$ as an alternative time coordinate. Then $\,(\tilde{a}_r,\tilde{b}_r)\,$ has
to be a time inhomogeneous diffusion on $\,\R\times\R_+^*\,$ started from $\,(a_0,b_0)\,$ and with
infinitesimal generator 
$$ \tilde{\GG} := {\s^2\over 2\,T(r,\tilde{a},\tilde{b})}\,\Big( 
(\tilde{a}^2-1+{\ts{R\over r}})\,{\p^2\over\p\tilde{a}^2} + (\tilde{b}^2 + r^2) \,{\p^2\over\p
\tilde{b}^2} + 2\tilde{a} \tilde{b}\,{\p^2\over\p\tilde{a}\p \tilde{b}} + 3\tilde{a}
\,{\p\over\p\tilde{a}} + (3\tilde{b} + {r^2\over \tilde{b}})\,{\p\over\p \tilde{b}}\Big) \, . $$  

   Conversely, given such a diffusion, the inverse time change yields a diffusion with generator
$\,\GG''\,$. Therefore it is enough to prove existence and uniqueness for the $\,\tilde{\GG}$-diffusion stopped at $\,r=R\,$. Now it is well known that a sufficient condition is locally boundedness
and continuity of the coefficients of the associated stochastic differential equation, together with a
local Lipschitz condition on these coefficients with respect to $(\tilde{a},\tilde{b})\,$. \par 
   Now the term $\,{r^2/\tilde{b}}\,$ causes no trouble since $\,\tilde{b}_0\not= 0\,$ and since the
proof of Theorem (\ref{cor.bh},2) insures that $\,\tilde{b}_r\,$ stays in $\,\R_+^*\,$, and then 
we observe that $\;(\tilde{a}^2-1+{\ts{R\over r}})\Big/ T\,$ goes to 0 as $\,r\sea 0\,$, and that 
$\,T^{-1/2}\,$, $\,T\1\,$, and their derivatives with respect to $(\tilde{a},\tilde{b})\,$, stay bounded
as well ; indeed  (for example) \quad 
$ T^{-3/2}\,|\,{\p T\over\p\tilde{b}}\,| = 2\, T^{-5/2}\,({\ts{R\over r}}-1)\,{\ts{\tilde{b}\over r^2}}  
\le 2\,T^{-1/2}/\tilde{b}\;$ is bounded. \par \smallskip 
   3) \ It remains to prove that we can recover $\,(\t_s,\dot\t_s)\,$ conditionally from $\,(r_s, b_s,
T_s)$, once $\,(\t_0,\dot\t_0)\,$ is fixed, in a unique way. \ 
   Now Proposition \ref{pro.systeqst} displays the equation we have to solve. To solve it, let us
complete the equation of Proposition \ref{pro.systeqst} into a linear system in the variables 
$\;V_s:=\Big(\t_s, (\dot\t_s/U_s), \t_s\wedge (\dot\t_s/U_s)\Big)$, by adjoining the equation
$\; d\t_s = U_s\,(\dot\t_s/U_s)\,ds\,$, \ and the following one (immediately deduced from Proposition \ref{pro.systeqst}) : 
$$ d\Big(\t_s\wedge{\dot\t_s\over U_s}\Big) = -\Big({r_s\over b_s}\,\s\,d\b_s\Big)\,
{\dot\t_s\over U_s} - \Big({\s^2\,r_s^2\over 2\,b_s^2}\,ds\Big)\, \t_s\wedge{\dot\t_s\over U_s}\; . $$  
So we get a linear differential system : $\; dV_s = V_s\,dA_s\,$, \ where the matrix-valued differential
$\,dA_s\,$ is given by : 
$$ dA_s := \pmatrix { 0 & -(b_s/r^2_s)\, ds & 0 \cr 
(b_s/r_s^2) ds & -({\s^2r_s^2/2b_s^2})\,ds & -({\s r_s/b_s})\,d\b_s \cr 
0 & ({\s r_s/b_s})\,d\b_s & -({\s^2r_s^2/2b_s^2})\,ds \cr } . $$ 
Note that this differential system can be equivalently written in the Stratonovitch form $\; dV_s =
V_s\,\circ\,d\tilde{A}_s\,$, \ where the matrix-valued differential
$\,d\tilde{A}_s\,$ is given by : 
$$ d\tilde{A}_s := \pmatrix { 0 & -(b_s/r^2_s)\, ds & 0 \cr 
(b_s/r_s^2) ds & 0 & -({\s r_s/b_s})\,d\b_s \cr  0 & ({\s r_s/b_s})\,d\b_s & 0 \cr } , $$ 
so that any solution takes its values in the rotation group. \par 

   Let us solve this linear equation by means of the following series : \parn 
\centerline{${\ds V_s=V_{D'}\,\Big(1+\sum_{k\in\N^*}J_k(s)\Big)}\,$, \ where for each $k$\quad 
${\ds J_k(s) := \int_{\{ D'<s_1<..<s_k<s\} } dA_{s_1}\times ..\times dA_{s_k}\, . } $} \par 
To justify that, let us choose on the space of $(3,3)$-matrices the Euclidian operator norm, \ and fix
$\;C=C(\o)\ge 1\,$, measurable with respect to
$\,(r_\cdot, b_\cdot, T_\cdot)$, such that for $\,D'\le s\le \min\{ D'',D'+1\}\,$ we have : 
$$ b_s/r_s^2\le C\,(s-D')^{-4/5}\; ,\; \s\,r_s/b_s\le C\,(s-D')^{2/5}\; ,\; 
\s^2r_s^2/b_s^2\le C\; ; $$ 
this is possible by Theorem \ref{cor.bh}. Let us suppose that 
$\; \E\Big[\Vert J_k(s)\Vert^2 \Big|\FF^o\Big] \le {(5C)^{2k}\over k!}\times (s-D')^{2k/5}$, \parn 
where $\,\FF^o\,$ denotes the $\s$-field generated by the reduced diffusion 
$\,(r_\cdot, b_\cdot, T_\cdot)\,$. Then we have 
$$ \E\Big[\Vert J_{k+1}(s)\Vert^2 \Big|\FF^o\Big] \; = \;\E\Big[\Big\Vert \int_{D'}^s J_k(s')\,dA_{s'}
\Big\Vert^2 \Big|\FF^o\Big] \qquad $$ 
$$ \le 2 \int_{D'}^s \E\Big[\Vert J_k(s')\Vert^2\Big|\FF^o\Big]  ({\s\,r_{s'}\over b_{s'}})^2\,ds'
+ 2 \int_{D'}^s {b_{s'}\over r_{s'}^2}\, ds'\times \int_{D'}^s \E\Big[\Vert J_k(s')\Vert^2\Big|\FF^o\Big] 
\; {b_{s'}\over r_{s'}^2}\, ds' $$ 
$$ + 2\,(s-D')\times\int_{D'}^s \E\Big[\Vert J_k(s')\Vert^2\Big|\FF^o\Big] 
\;({\s^2r_s^2\over 2b_s^2})^2\, ds' $$
$$ \le 2\,C^2\,{(5C)^{2k}\over k!}\left[ \,\int_{D'}^s (s'-D')^{(2k+4)/5}\, ds' + 5 (s-D')^{1/5} 
\int_{D'}^s (s'-D')^{(2k-4)/5}\, ds' \right. \hskip 10mm $$
$$ \hskip 70mm \left.+ \,(s-D') \int_{D'}^s (s'-D')^{2k/5}\, ds'\right ] $$ 
$$ = 2\,C^2\,{(5C)^{2k}\over k!}\left[ {(s-D')^{(2k+9)/5}\over 2k+9}  + 5\, {(s-D')^{(2k+2)/5}
\over 2k+1 } + {(s-D')^{(2k+10)/5}\over 2k+5}\right ] $$
$$ \le {(5C)^{2(k+1)}\over (k+1)!}\times (s-D')^{2(k+1)/5}\, . $$ 
Since the series $\;{\ds \sum_k {(5C)^{k}\over (k!)^{1/2}}}\,$ converges, this shows by induction that the
series $\,{\ds\sum_{k}J_k(s)}\,$ almost surely converges in $L^2\,$, conditionally with respect to the
reduced diffusion $\,(r_\cdot, b_\cdot, T_\cdot)\,$. \par 
   As to the unicity of this solution, note that any other solution $\,V'_s\,$ must satisfy for any
$\,n\in\N^*$ : 
$$ V'_s = V_{D'}\,\Big(1+\sum_{k=1}^n J_k(s)\Big) + \int_{\{ D'<s_0<s_1<..<s_n<s\} }
V'_{s_0}\times dA_{s_0}\times dA_{s_1}\times ..\times dA_{s_n}\, , $$ 
so that, using that $\,\Vert V'_{s_0}\Vert = \Vert V_{D'}\Vert\,$ is constant, we get (with the
$L^2$-norms understood conditionally as above) :  
$$ \Vert V'_{s}-V_s\Vert_2 \le \Big\Vert \int_{\{ D'<s_0<s_1<..<s_n<s\} }
V'_{s_0} dA_{s_0} .. dA_{s_n}\Big\Vert_2 + \Big\Vert\sum_{k>n} J_k(s)\Big\Vert_2
\le \Vert J_{n+1}(s)\Vert_2 + \sum_{k > n} \Vert J_k(s)\Vert_2 $$ 
$$ \le 2\sum_{k > n} {(5C)^{k}\over (k!)^{1/2}}\quad {\rm for \ any }\; n\in\N^*\,, \quad {\rm showing
\ that } \quad  V'_{s}= V_s \quad  {\rm almost \ surely. } $$ 

   This concludes the whole proof of Proposition \ref{the.martp}. $\;\diamond $ 

\subsection{Proof of the second assertion of Theorem \ref{the.alls} : $\,D_\ii=\ii\,$ a. s.}\indf \label{sec.PT5}
  For any $\,\e\in ]0,1[\,$, \ set $\;A_\e := \Big\{ \sup_{s<D_\ii}\limits b\, \le \e\1
\Big\}\,$. \quad By continuity of $\,b\,$, we have \ $\P (\{ D_\ii<\ii\}\cap A_\e^c) <\e'\,$, for any fixed $\,\e'>0\,$ and small enough fixed $\,\e\,$. \parn   For any $n\in\N$, \ set $\;\tau'_n :=\inf\{s>D_{3n+2}\,|\,T_s=0\}\,$ and $\;\tau_n:=\inf\{s>D_{3n+3}\,|\,r_s=R/2\,\}\,$. \parn 
Note that $\; D_{3n+2}<\tau'_n<D_{3n+3}<\tau_n<D_{3n+4}\;$ on $\{ D_\ii<\ii\}$,  
and $\;\sum_n\limits\,(\tau_n-\tau'_n) <D_\ii\,$. \parn 
Moreover $\;\tau_n-\tau'_n = \tau\circ\Theta^{\tau'_n}\,$, where $\,\tau\,$ denotes the hitting time of 
$\,\{ r=R/2\,\}$. \parn 
   Setting $\;A'_\e:= \{\tau< c\e\}\cap \Big\{ \sup_{s\le D'}\limits b\, \le \e\1\Big\}\,$ (for some constant $\,c=c(\s,R) >0\,$), \ we have
$$ \P ( D_\ii<\ii)-\e' < \P (\{ D_\ii<\ii\}\cap A_\e)\le \P\Big(\sum_n\,(1_{A_\e}\times\tau)\circ\Theta^{\tau'_n} <\ii\Big) $$
$$ \le \P\Big(\liminf_n\, (\Theta^{\tau'_n})\1(A'_\e)\Big) \le \sum_n\, 
\lim_{p\to\ii}\P\Big(\bigcap_{m=n}^p\limits\, (\Theta^{\tau'_m})\1(A'_\e)\Big) $$ 
$$ = \sum_n\, \lim_{p\to\ii}\E\Big(
\P_{\xi_{\tau'_p}}(A'_\e)\times\prod_{m=n}^{p-1}\limits\,(1_{A'_\e}\circ\Theta^{\tau'_m})\Big) $$ 
by the strong Markov property. Hence we see that this proof is achieved if we show that 
$\;\lim_{p\to\ii}\limits\E\Big(\P_{\xi_{\tau'_p}}(A'_\e)\times\prod_{m=n}^{p-1} \limits\,
(1_{A'_\e}\circ\Theta^{\tau'_m})\Big) = 0\,$, for any $n$. Now this follows immediately, by induction, from the following lemma. $\;\diamond $ 
\begin{lem}\label{lem.maj}\quad For small enough $\,\e>0\,$, and for any initial condition $\,r_0>R\,$
and $\,b_0>0$, provided $\,T_0=0$, $\,A'_\e\,$ being as in the proof of Theorem \ref{the.alls} above, we have $\;\P_{}(A'_\e) <1/2\,$.  
\end{lem} 
\ub{Proof}\quad  Let us write again the stochastic equation of $\,T_s\,$, under the following form,
for some real standard Brownian motion $\,W$ :
$$ T_s= \s\,W\Big[\int_0^s(T_t^2+1-{\ts{R\over r_t}})\, dt\Big] +{\ts{3\s^2\over 2}}\int_0^sT_t\,dt 
+\int_0^s(r_t-{\ts{3R\over 2}})\,{\ts{b_t^2\over r_t^4}}\,dt-\int_0^s{\ts{R\over 2r_t^2}}\,dt\; . $$ 

   Consider $\;\tau:=\inf\{ s\,|\,r_s=R/2\}\;$ and $\;\s'_\e:=\inf\{ s\,|\,|T_s|=\e\1\}\,$, and some
constant $\,q\,$. \par 
   On the event $\;A''_\e:=\Big\{ \max |W| ([0,(\e\2 +1)q]) < (2\s\e)\1\Big\}\cap \Big\{ \sup_{s\le
D'}\limits b\,\le \e\1\Big\}\,$, we have for $\,0\le s\le \min\{ s,\s'_\e,\tau\}\,$ : \quad 
${\ds r_s=r_0+\int_0^sT_t\,dt > R -\e\1 s}\;$, \  and by the equation of $T_s$ : 
$$ -{\ts{1\over 2\e}}-{\ts{3\s^2\over 2\e}}\,q - 3\,(2/R)^{3}\e\2\,q - (2/R)\,q <T_s <
{\ts{1\over 2\e}}+ {\ts{3\s^2\over 2\e}}\,q + (2/R)^{3}\e\2\,q \; , $$
whence \quad 
${\ds |T_s| < {{1\over 2\e}}+{{1\over 2\e}}\,\Big( 3\s^2 +4\,R\1\e + 48\,R^{-3}\e\1\Big)\,q \, . }$ \parn 
Hence $\;q <\min\{\s'_\e,\tau\}\,$ on $\,A''_\e\;$ if $\; q\le \Big( 3\s^2 +4\,R\1\e +
48\,R^{-3}\e\1\Big)\1\,$ and $\; q\le R\e/2\,$. \parn  
Thus there exists a constant $\, c=c(\s,R)>0\;$ such that $\;\min\{\s'_\e,\tau\}> q=c\,\e\;$ on
$\,A''_\e\,$.
\par  
   Therefore $\;\Big\{ \max |W|([0,(\e\2 +1)q]) < (2\s\e)\1\Big\}\cap A'_\e =\emptyset \;$, and then 
$$ \P(A'_\e)\le \P\Big[\max |W|([0,(\e\2+1)q]) \ge (2\s\e)\1\Big]  = \P\Big[\max
|W|([0,q])\ge (2\s\,\sqrt{\e^2 +1}\,)\1\Big]  $$ 
$$ \le 2\, \P\Big[\max W([0,q])\ge (2\s\,\sqrt{\e^2 +1}\,)\1\Big] = 
2\,\P\Big[ |W_1|\ge (2\s\,\sqrt{\e^2 +1}\,)\1/\sqrt{c\,\e}\,\Big] $$ 
$$ =\; 4\int_{1/(2\s\sqrt{c\,\e\,(1+\e^2)}\,)}^\ii e^{-x^2/2}{dx\over\sqrt{2\pi}} \;\le \;
{12\,\s\over\sqrt{2\pi}}\,\sqrt{c\,\e\,}\,e^{-1/(12\,\s^2c\,\e)}\,. \;\;\diamond $$ 

\subsection{Proof of Theorem \ref{the.piege}}\indf \label{sec.PT6}
   The proof of this theorem will be split into 8 parts and a series of lemmas. \par\smallskip 

   {\bf 1)} \ Let us begin by the dichotomy of the first assertion : there is no other possibility
than the obvious one : $\,r_\cdot\,$ goes to infinity, and the confinement exhibited here :
$\,r_\cdot\,$  remains endlessly bounded.
\begin{lem}\label{lem.dicho}\quad Almost surely, if $\;r_\cdot\,$ does not go to infinity, \ then it is
bounded. 
\end{lem}
\ub{Proof}\quad  Consider the double sequence of hitting times $\,\{\L_n, \L'_n\,|\,n\in\N\}\,$ defined by $\;\L_0=\L_0'=0\,$, $\;\L'_n:=\inf\{ s>\L_n\,|\,r_{s}<R\}\,$,  $\;{\L_n} := \inf\{
s>\L'_{n-1}\,|\,\, r_{s}>e^{n}\;\, \hbox{ and }\;\; T_s=0\}$, and set $\;E_n:=\{ \L_n<\ii\}$, for
$\,n\in\N^*$. \par 
   By Theorem (\ref{the.asymp},1) the event $\,\cap_n\,E_{n}\,$ contains all trajectories such that $\,r_s\,$ does not go to infinity, but is unbounded. Thus we want precisely to prove that 
$\;\P(\cap_n\,E_{n}) =0\,$. \par 
    Let us apply the comparison theorem (see for example ([I-W], Theorem 4.1)) : there exist comparison processes $\,X^+_s\,$ and $\,X^-_s\,$ on the same probability space, such that almost surely on $\,E_n\,$ : 
$$ \max_{\L_n\le s\le \L_n+1}\limits T_{s}\ge \max_{\L_n\le s\le \L_n+1}\limits X^+_s\,,\;  X^+_{\L_n}=0\quad\hbox{ and } \; \min_{\L_n\le s\le \L_n+1}\limits T_{s}\ge \min_{\L_n\le s\le \L_n+1}\limits X^-_s\,, \;X^-_{\L_n}=0\,, $$ 
where $\,X^+_s\,$ and $\,X^-_s\,$ are real diffusions given by their stochastic differential equations,
which are deduced from the equation governing $\,{T}_s\,$ by bringing down the ratio between the drift term and the diffusion coefficient, and by bringing the diffusion coefficient down (in the $X^+$ case) or up (in the $X^-$ case). \  Recall that 
$$ d{T_s} = \s\,\sqrt{T_s^2+1-{\ts{R\over r_s}}}\,\, dw_s + {\ts{3\s^2\over 2}}\,T_s\,ds +
 (r_s-{\ts{3R\over 2}})\,{b_s^2\over r_s^4}\,ds - {R\over 2 r_s^2}\,ds \, . $$ 
Thus, as long as $\,r_{s}\ge 3R/2\,$, we can take $\,{\ts{1\over 3}}\le T^2+{\ts{1\over 3}}\le
T^2+1-{\ts{R\over r}}\,\le  T^2+1\,\,$ and then : 
$$ dX^+_s= {\ts{\s\over \sqrt{3}}}\,dw_s+\,{\ts{3\s^2\over 2}}\,(X^+_s\wedge 0)\,ds -{\ts{2\over 9R}} \,ds\,, $$  
and 
$$ dX^-_s = \s\,\sqrt{(X^-_s)^2+1}\, dw_s + \,\Big(
{\ts{9\s^2\over 2}}\,(X^-_s\wedge 0)  - {\ts{2\over 3R}}\Big)\, ds\, .  $$ 

   Now, for $\,A:= 4+{\ts{4\over R\s^2}}\,$, \ fix $\;{\ds \e := \P\Big(\max_{\L_n\le s\le \L_n+1}\limits X^+_s>  A\Big) }\,$,  which is $\,>0\,$, \ $\,N>0\,$ such that 
$\;{\ds \P\Big(\min_{\L_n\le s\le \L_n+1}\limits X^-_s< -N\Big) <\e/2 }\,$, \ and $\,n\,$ such that 
$\,r_{\L_n}>N+3R\,$.  \ Note that the event 
$\; E'_n := \Big\{ \max_{\L_n\le s\le \L_n+1}\limits X^+_s> A\Big\}\cap 
\Big\{ \min_{\L_n\le s\le \L_n+1}\limits X^-_s\ge -N\Big\}\,$ has probability $\,>\e/2\,$. \par 

   Setting $\;\L''_n:=\inf\{ s>\L_n\,|\,r_{s}\le 3R/2\}\,$ and applying the comparison theorem, we get : 
$$ \min_{\L_n\le\, s\,\le\, (\L_n+1)\wedge\L''_n}\limits T_{s}\,\ge\, \min_{\L_n\le s\le \L_n+1}\limits X^-_s \,\ge -N \quad  \hbox{ on } E_n\cap E'_n\, , $$ 
and then 
$$ \min_{\L_n\le\, s\,\le\, (\L_n+1)\wedge\L''_n}\limits r_{s}\,\ge \, r_{\L_n} -N \ge 3R \quad 
\hbox{ on } E_n\cap E_n'\, , $$ 
showing that $\;\L''_n >1\;$ on $\,E_n\cap E_n'\,$. \ Hence, applying the comparison theorem again, we get also : 
$$ \max_{\L_n\le s\le \L_n+1}\limits T_{s}\,\ge\, \max_{\L_n\le s\le \L_n+1}\limits X^+_s \,> A
\quad  \hbox{ on } E_n\cap E'_n\, . $$ 

   Thus, using the strong Markov property, we find that for large enough $\,n\,$ :
$$ \P\Big( (\exists \; s>\L_n)\quad T_s > 4+{\ts{4\over R\s^2}} \;\;\hbox{ and }\;\; r_s>3R/2
\,\Big/\,E_n\Big) \,>\,\e/2\, . $$ 

   Let us use now the proof of Theorem (\ref{the.asymp},2), where we proved that for $\,r_0\ge 3R/2\,$ and $\, T_0\ge 4+{\ts{4\over R\s^2}}\,$, then $\;\P(\lim_{s\to\ii}\limits r_s=\ii)\ge 1/2\,$,
together with the Markov property : we obtain 
$$ \P\Big( E_{n+1}\,\Big/\,E_n\Big)\,\le\,\P\Big(\L'_n<\ii\,\Big/\,E_n\Big)\, <\, 1- \e/4\, ,\;\hbox{
for } n \hbox{ larger than some } n_0\,. $$   
Therefore we get $\;\P( E_{n_0+k+1}) < (1- \e/4)^k\,$ and then $\;\P(\cap_n\,E_{n}) =0\,$, as wanted.  $\;\diamond$ 
\par\medskip  \medskip

   {\bf 2)} \ The first case in the theorem was already handled in Lemma \ref{lem.cv} and Theorem
\ref{the.asymp}. The irreducibility of the relativistic diffusion is clear from Theorem (\ref{the.asymp},2). \par 
   Let us now focus on the second case in the theorem, supposing therefore that $\,r_\cdot\,$ is bounded. \  Set $\;\tau_M:=\inf\{ s\,|\,r_s= M\}$, for $\,M\,$ of the form $\,M= k\times 3R/2\,$, with $\,k\in\N^*$. \par 
   The proof is divided in several distinct lemmas. We shall always let the relativistic diffusion  start at proper time $\,D_{-1}=0\,$ from level $\{ r=R\}$ with $\,T_0>0\,$. We let the hitting times $\,D_j\,$ be as in Section \ref{sec.SDalls}. \ Fix also some $\,\e\,$ in $\,]0,{1\over 3}[\,$. \par\medskip  \medskip 

   {\bf 3)} \  {\it Estimates related to} $\,b_\cdot\,$ \par\medskip  

   Let us begin by proving that (when $\,r_\cdot\,$ is bounded) $\,b_\cdot\,$ has to go to infinity, 
which by means of the Markov property will allow then to consider only large enough $\,b_0\,$. 
\begin{lem}\label{lem.bdiv}\quad  Almost surely, if $\;r_s\,$ is bounded, \ then 
$\;{\ds \lim_{s\to\ii}\,b_s\,e^{-(1-\e)\s^2s} = +\ii \,}$.  
\end{lem}
\ub{Proof}\quad  Note that $\{\sup_{s\ge 0}\limits\,r_{s} <\ii\}=\bigcup_{M}\limits\,\{\tau_{M}=\ii\}$. 
\parn 
Let us recall the logarithmic form of the stochastic differential equation of $\,b\,$ :  \parn 
${\ds \log (b_{s}\,e^{-(1-\e)\s^2s}) = \log b_{0}\,+ \e\s^2\,s +\s\, W\Big[ s+\int_{0}^{s}{\ts r_t^2\over
b_t^2}\,dt\Big]}$ \quad (for some Brownian motion $\,W$). \parn
A straightforward consequence is that $\;\limsup_{s\to\ii}\limits (b_{s}\,e^{-(1-\e)\s^2s}) = +\ii\;$  
almost surely, so that  the stopping time $\,t_{n}:=\inf\{ s>0\,|\, b_{s}\,e^{-(1-\e)\s^2s}>n\}$ is finite
for any $\,n\in\N\,$. \par 
   Then let us write the equation governing $\,b_\cdot\,$ in the following linearized form 
(for some other Brownian motion $\,\tilde{W}$) : 
$$ b_{s}\,e^{-(1-\e)\s^2s}\, =\, e^{\s\,w_{s}+\e\,\s^2 s} \Big( b_{0} + \s\tilde{W}\Big[ 
\int_{0}^{s} e^{-2\s\,w_{t}-2\s^2 t}\, r_t^2\,dt\Big] + \int_{0}^{s}
e^{-\s\,w_{t}-\s^2 t}\, {\ts{\s^2\, r_t^2\over 2\, b_t}}\,dt\,\Big) . $$  
It is clear from this expression that for any fixed $\,M\,$ the probability of the asymptotic event  
$\;A_{M}:= \{\tau_{M}=\ii\}\cap \{b_{s}\,e^{-(1-\e)\s^2s}\,\hbox{ does not go to }\ii\}\,$ 
can be made arbitrarily close to 0 by choosing $\,b_{0}\,$ large enough. \ 
Using that $\,A_{M}= \Theta_{t_{n}}\1(A_{M})\,$, we apply the strong Markov property
at the sequence of stopping times $\,\{t_{n}\}\,$ to conclude that $\,\P(A_{M})=0\,$, which yields the result. 
$\;\diamond$     \par\medskip   

   We estimate then the increase of $\,b\,$, when started from some large value $\,b_0\,$. 
\begin{lem}\label{lem.estim2b}\quad  Fix $\,\e>0$, and recall that $\;\tau_M:=\inf\{ s\,|\,r_s= M\}$.
Then there exists a lower bound $\,b(M,\e)\,$ such that for $\,b_0\ge b(M,\e)\,$ we have : 
$$ \P\Big( b_0^{1-(\e/2)}\, e^{(1-\e)\s^2s}\le b_s \le b_0^{1+(\e/2)}\, e^{(1+\e)\s^2s} \quad\hbox{ for
all }\; s\le \tau_M\Big) >1- 2\, b_0^{-\e^2/2} >\, 1-\e\, . $$  
\end{lem}
\ub{Proof}\quad  For any $\,q\in]0,1[\,$, set \quad $\nu_q:=\inf\{ s\,|\,
b_s\, e^{-(1-\e)\s^2s}\le qb_0\}\in ]0,\ii]\,$. \parn
Recall that the equation governing $\,b_\cdot\,$ can be expressed in the following form : \parn 
\centerline{${\ds\log (b_s/b_0) = \s\,W\Big[\int_0^s \Big( 1+{r_t^2\over b_t^2}\Big)\,dt\,\Big] +\s^2\,
s}\;$, \  for some standard Brownian motion $\,W_\cdot\,$.} \parn 
Hence for $\,0\le s\le\nu_q\wedge\tau_M\;$ we have : 
$$ \log (b_{s}\,e^{-(1-\e)\s^2s}/b_0) \ge \e\s^2s\,+\,\s\,\min\,W[s,(1+(M/qb_0)^2)s]\, , $$ 
whence
$$\P(\nu_q <\tau_M)\,\le\,\P\Big( (\exists\,s)\;\;\e\,s+\min\,W[s,(1+(M/qb_0)^2)s]\le\log q\Big)\qquad $$ 
$$ = \,\P\Big( (\exists\,s')\;\; {\e\, s'\over {\ts 1+(M/qb_0)^2}} +
\min\,W[{\ts{s'\over 1+(M/qb_0)^2}},s']\le \log q\Big) $$ 
$$ \qquad \le \,\P\Big( (\exists\,s)\;\; {\e\, s\over 1+(M/qb_0)^2} + W_s \le \log q\Big)\, 
= \,\exp\Big( {2\e \, \log q\over 1+(M/qb_0)^2 }\Big) \, . $$ 
Taking $\;q=b_0^{-\e/2}\,$, this yields : 
$\;\P(\nu_q <\tau_M)\le b_0^{-\e^2/[1+M^2b_0^{\e-2}]} < b_0^{-\e^2/2}\,$
for $\,b_0\ge b(R,\e)\,$. \parn 
This proves the lower control of the statement. The upper control is obtained exactly in the same way.
$\;\diamond$ 
\medskip \medskip

   {\bf 4)} \ We prove now that $\;\rr\,:=\,\limsup_{s\to\ii}\limits\,r_s\;$ cannot be strictly between
$\,3R/2\,$ and $+\ii\,$. 

\begin{lem}\label{lem.limsup}\quad  Almost surely, if $\;r_s\,$ is bounded, \ then 
$\quad \rr\,:=\,\limsup_{s\to\ii}\limits\, r_s \le 3R/2\,$. 
\end{lem}
\ub{Proof}\quad  Let us proceed somewhat as for \mbox{Lemma \ref{lem.dicho}.} Fix $\,\e\in]0,1[\,$, 
$\,M> 3\e+3R\,$, and consider the double sequence of hitting times $\,\{\L_n, \L'_n\,|\,n\in\N\}\,$
defined by $\quad\L_0=\L_0'=0\,$, \parn 
\centerline{${\L_n} := \inf\{ s>\L'_{n-1}\,|\,\, r_{s}>2\e+3R/2\;\, 
\hbox{ and }\;\; T_s=0\;\, \hbox{ and }\;\;  b_s> e^{2n}\}$,}\parn 
\centerline{$\L'_n:=\inf\{ s>\L_n\,|\,\,r_s < \e+3R/2 \;\,\hbox{ or }\;\; b_s< e^n\}$.}  \parn 
Recall that $\quad\tau_M:=\inf\{ s>0\,|\,r_s>M\}\,$, and set $\quad E_n:=\{\L_n<\tau_M\}$, \ for $\,n\in\N^*$. \par 

  By Theorem (\ref{the.asymp},1) and by Lemma \ref{lem.bdiv}, the event $\,\cap_n\,E_{n}\,$ contains all trajectories such that  $\;r_\cdot\,$ is bounded by $\,M$  and $\;\limsup_{s\to\ii}\limits\, r_s>
2\e+3R/2\,$. Thus we want precisely to prove that  $\;\P(\cap_n\,E_{n}) =0\,$. \par 

    Let us apply the comparison theorem (see for example ([I-W], Theorem 4.1)) : there exist comparison processes $\,X^+_s\,$, $\,X^-_s\,$ and $\,Y^-_s\,$ on the same probability space, such that \parn 
almost surely on $\,E_n\,$ : $\quad \min_{\L_n\le\, s\,\le\, \L_n+1}\limits b_{s}\ge \min_{\L_n\le\, s\,\le\, \L_n+1}\limits Y^-_s\;, \;Y^-_{\L_n}=b_{\L_n}>e^{2n}\,,$ \quad and 
$$ \max_{\L_n\le\, s\,\le\, \L_n+1}\limits T_{s}\ge \max_{\L_n\le\, s\,\le\, \L_n+1}\limits X^+_s\;, \;X^+_{\L_n}=0\quad\hbox{ and } 
\quad \min_{\L_n\le\, s\,\le\, \L_n+1}\limits T_{s}\ge \min_{0\le s\le 1}\limits X^-_s\;, \;X^-_{\L_n}=0\,, $$ 
where $\,X^+_s\,$ and $\,X^-_s\,$ are real diffusions given by their stochastic differential equations,
which are deduced from the equation governing $\,{T}_s\,$ by bringing down the ratio between the drift
term and the diffusion coefficient, and by bringing the diffusion coefficient down (in the $X^+$ case) or
up (in the $X^-$ case), and similarly for $\,b_s,Y^-_s\,$. \  Recall that 
$$ d{b_s} = \s\,\sqrt{b_s^2+r_s^2}\,\, dw'_s + {\ts{3\s^2\over 2}}\,b_s\,ds +
 \,{\s^2\,r_s^2\over 2\,b_s}\,ds\, , \quad \hbox{ and } $$ 
$$ d{T_s} = \s\,\sqrt{T_s^2+1-{\ts{R\over r_s}}}\,\, dw_s + {\ts{3\s^2\over 2}}\,T_s\,ds +
 (r_s-{\ts{3R\over 2}})\,{b_s^2\over r_s^4}\,ds - {R\over 2 r_s^2}\,ds \, . $$ 
Thus, on $\,[\L_n,\L_n']\,$ we can use $\,b^2+r^2\,\le \,b^2+M^2\,\,$ and 
$\,T^2+{\ts{1\over 3}}\le T^2+1-{\ts{R\over r}}\,\le T^2+1\,$, and then, for $\,n\,$ such
that  
$\;\e M^{-4}\,e^{2n} - {\ts{2\over 3R}} > 3\,e^n\,$ : 
$$ dY^-_s = \s\sqrt{(Y^-_s)^2+M^2}\,\, dw'_s + \,{\ts{3\s^2\over 2}}\,Y^-_s\,ds\, , $$
$$  dX^+_s = \s\,\sqrt{(X^+_s)^2+{\ts{1\over 3}}}\,\,dw_s + \, {\ts{3\s^2\over 2}}\,(X^+_s\wedge 0)\,ds  
+ e^{n}\, ds \, , $$  
and 
$$ dX^-_s = \s\sqrt{(X^-_s)^2+1}\,\,dw_s + \,{\ts{9\s^2\over 2}}\,(X^-_s\wedge 0)\,ds + 3\, e^{n}\,
ds \, . $$ 
Now $\;{\ds \P\Big(\min_{\L_n\le s\le \L_n+1}\limits Y^-_s> e^n\Big) }\,$,  
$\;{\ds \P\Big(\min_{\L_n\le s\le \L_n+1}\limits X^-_s > -\e\Big) }\,$, and 
(for $\,A:= 4+{\ts{4\over R\s^2}}\,$) \parn 
${\ds \P\Big(\max_{\L_n\le s\le \L_n+1}\limits X^+_s > A\Big) }\,$
are (for large enough $\,n$) arbitrary near from 1, so that the event \parn 
$F_n := \Big\{ \min_{\L_n\le s\le \L_n+1}\limits Y^-_s> e^n\Big\}\cap \Big\{ \min_{\L_n\le s\le \L_n+1} \limits X^-_s > -\e\Big\} \cap \Big\{ \max_{\L_n\le s\le \L_n+1}\limits X^+_s > A\Big\}\,$ has probability arbitrary near from 1. \par 

   Hence, applying the comparison theorem, we get : 
$$ \min_{\L_n\le\, s\,\le\, (\L_n+1)\wedge\L'_n\wedge\tau_M}\limits b_{s}\,\ge\, 
\min_{\L_n\le s\le \L_n+1}\limits Y^-_s \,> e^n \quad  \hbox{ on } F_n\cap E_n\, , $$ 
$$ \min_{\L_n\le\, s\,\le\, (\L_n+1)\wedge\L'_n\wedge\tau_M}\limits T_{s}\,\ge\, 
\min_{\L_n\le s\le \L_n+1}\limits X^-_s\,> -\e \quad  \hbox{ on } F_n\cap E_n\, , $$ 
and then 
$$ \min_{\L_n\le\, s\,\le\, (\L_n+1)\wedge\L'_n\wedge\tau_M}\limits r_{\L_n+s}\,> \, r_{\L_n} -\e > 
\e+3R/2 \quad \hbox{ on } F_n\cap E_n\, , $$ 
showing that $\;\L'_n>\L_n+1\;$ on $\,F_n\cap\{ \tau_M \ge \L'_n\}\cap E_n\,$. \ Hence, applying the
comparison theorem again, we get also : 
$$ \max_{\L_n\le s\le \L_n+1}\limits T_{s}\,\ge\, \max_{\L_n\le s\le \L_n+1}\limits X^+_s \,> A
\quad  \hbox{ on } F_n\cap\{ \tau_M \ge\L'_n\}\cap E_n\, . $$ 

   Thus, using the strong Markov property, we find that for large enough $\,n\,$ :
$$ \P\Big( \tau_M <\L'_n \;\;\hbox{ or }\;\;\Big[ (\exists \; s\in [\L_n,\L'_n[\,)\quad T_s >
4+{\ts{4\over R\s^2}} \;\;\hbox{ and }\;\; r_s>3R/2\Big] \,\Big/\,E_n\Big) \,>\,1-\e\, . $$ 

   Let us use now the proof of Theorem (\ref{the.asymp},2), where we proved that for $\,r_0\ge
3R/2\,$ and $\, T_0\ge 4+{\ts{4\over R\s^2}}\,$, then $\;\P( r_s\;\hbox{ increases towards }\;\ii)\ge
1/2\,$, together with the Markov property : we obtain 
$$ \P\Big( \tau_M <\L'_n\,\Big/\,E_n\Big)\,>\, 1/2- \e\, ,\;\hbox{ for } n \hbox{ larger than some }
n_0\,. $$ 

Therefore, noting that $\,E_{n+1}\subset \{ \tau_M>\L'_n\}$, we get $\;\P( E_{n_0+k+1}) < 
[1/2+ \e]^k\,$ and then $\;\P(\cap_n\,E_{n}) =0\,$, as wanted.  $\;\diamond$ 
\par\medskip \medskip

   {\bf 5)} \ {\it Lower estimates related to} $\,D_{3n+2}\,$ \medskip 

   The aim of the following lemma is to get a lower bound (in an optimal way, owing to Remark
\ref{rem.control2} below) on the duration of an excursion in the hole, ouside an event of small
probability. The idea is that when $\,T_\cdot\,$ vanishes, with an acceleration of order
$\,b^2_\cdot\,$ it takes a lapse of time of order $\,1/b_\cdot\,$ to make a non-infinitesimal move. 
\begin{lem}\label{lem.control1}\quad  Suppose $\;r_0=R\,$, $\,T_0>0\,$, and $\,b_0\,$ large enough.
$$ \hbox{Then we have }\quad \P\Big( D_1 > { {c\over b_0}}\wedge \tau_M\Big) 
>1- \,e^{- b_0/(64\s^2c)}\, , \quad \hbox{ for }\; c :=R^3/(512 M). $$ 
\end{lem}
\ub{Proof}\quad   Consider \quad $D_0' :=\inf\{ s>0\,|\,T_s=0\}\le D_0\,$, 
$\;D_1' :=\inf\{ s>0\,|\,r_s= {\ts{R\over 2}}\}\in [D'_0,D_1]$, \parn     
\centerline{ $\eta := \inf\{ s>0\,|\,\log({\ts{ b_s\over b_0}}) = \pm\log
2\}$, \, and $\;\la :=\inf\{ s>D'_0\,|\,T_s=\pm{\ts{b_0\over R}}\}\le D_0\,$.} \parn   
Set $\;E:=\{ D'_0<{\ts{c\over b_0}}\wedge \tau_M\}$,  $\;E':=\{ \s\, W^*[2c/b_0] <{\ts{1\over 2}}\}$,
\ and $\;E'':=\{ \s\, \tilde{W}^*[2cb_0/R^2] <{\ts{b_0\over 4R}}\}$, \parn 
where $\,W\,$ and $\,\tilde{W}\,$ are standard real Brownian motions used to write the stochastic
equations governing $\,b_\cdot\,$ and $\,T_\cdot\,$ respectively, and 
$\,W^*\,$, $\,\tilde{W}^*\,$ are their maximum processes : $\,W^*(s):=\max |W|([0,s])\,$ for any $s>0\,$
and similarly for $\,\tilde{W}^*$.\parn
Recall from the proofs of Lemmas \ref{lem.bdiv} and \ref{lem.estim2b} that \ 
$\;{\ds\log (b_s/b_0) = \s\,W\Big[\int_0^s \!\Big( 1+{r_t^2\over b_t^2}\Big)\,dt\,\Big] +\s^2\,
s}\,$, and (as already used for example in the proof of Lemma
\ref{lem.maj}) that for $\,s\ge D'_0\,$ :
$$ T_s= \s\,\tilde{W}\Big[\int_{D'_0}^s(T_t^2+1-{\ts{R\over r_t}})\, dt\Big] +{\ts{3\s^2\over
2}}\int_{D'_0}^sT_t\,dt  +\int_{D'_0}^s(r_t-{\ts{3R\over 2}})\,{\ts{b_t^2\over
r_t^4}}\,dt-\int_{D'_0}^s{\ts{R\over 2r_t^2}}\,dt\, . $$  
Therefore on $\,E\cap E'$, for $\,0\le s\le {\ts{c\over b_0}}\wedge \tau_M\wedge \eta\;$ we have : 
$$ \Big|\log (b_s/b_0)\Big| \le \s\, W^*\Big[ (1+{\ts{4\,M^2\over b_0^2}})\,{\ts{c\over b_0}}\Big] +
{\ts{\s^2c\over b_0}} < \5 + {\ts{\s^2c\over b_0}} <\log 2\,, \; \hbox{ whence }\;\eta\ge {\ts{c\over
b_0}}\wedge \tau_M\, , $$  
and on the other hand, on $\,E\cap E'\cap E''$, for $\,D'_0\le s\le{\ts{c\over b_0}}\wedge \tau_M\wedge
\la\wedge D'_1\;$ we have : 
$$ |T_s| < \s\,\tilde{W}^*\Big[ ({\ts{b_0^2\over R^2}}+1)\,{\ts{c\over b_0}}\Big] +
{\ts{3\,\s^2\over 2}}\,{\ts{b_0\over R}}\,{\ts{c\over b_0}} + 2M\,{\ts{64\,b_0^2\over
R^{4}}}\,{\ts{c\over b_0}} + {\ts {2\over R}}\,{\ts{c\over b_0}} < {\ts{b_0\over 4R}}+ {\ts{3\s^2c\over
2R}} +{\ts{128 Mcb_0\over R^4}}+ {\ts{2c\over Rb_0}} < {\ts{b_0\over R}} \, , $$  
since $\; c=R^3/(512 M)\,$, for $\;b_0\,$ large enough.  Since $\,c\le R^2/2\,$, this implies that  
$$ r_s\,=\,R+\int_0^sT_t\,dt \,>\, R\,-\,{\ts{b_0\over R}}\,{\ts{c\over b_0}}
\;\ge \,{\ts{R\over 2}}\;,\quad \hbox{ still for }\; \,D'_0\le s\le{\ts{c\over b_0}}\wedge \tau_M\wedge
\la\wedge D'_1\,. $$ 
Hence  we find that $\;E\cap E'\cap E''\subset \{D_1'>{{c\over b_0}}\wedge \tau_M\}$, and then also 
$\;E'\cap E''\subset \{D_1'>{{c\over b_0}}\wedge \tau_M\}$. \par

$$ \P\Big(\!\tilde{\tau}>{c\over b_0}\!\Big) = \P\Big(\!\tilde{\tau}\wedge{c\over
b_0}<\tilde{\tau}\!\Big)  = \P\Big(\! \sup_{s<\tilde{\tau}\wedge{c\over b_0}} |\log({\ts{ b_s\over b_0}})|
<\log 3 , \sup_{s<\tilde{\tau}\wedge{c\over b_0}} |T_s| < {\ts{b_0\over R}} , 
\sup_{s<\tilde{\tau}\wedge{c\over b_0}} |r_s-R| <{\ts{R\over 2}}\!\Big) . $$ 

   Now, as already used in the proof of Lemma \ref{lem.maj}, we have for $\,b_0\,$ large enough : 
$$ \P[ E']\,>\, 1-8\s\,\sqrt{c/(\pi b_0)}\,e^{- b_0/(16\s^2c)}\, >\, 1-\5\,e^{-b_0/(16\s^2c)} $$ 
$$ \hbox{ and }\quad 
\P[ E'']\,>\, 1-16\s\,\sqrt{c/(\pi b_0)}\,e^{- b_0/(64\s^2c)}\,>\,1-\5\,e^{-b_0/(64\s^2c)}\, , $$ 
whence finally : 
$$\P\Big( D_1 > { {c\over b_0}}\wedge \tau_M\Big) \ge \P\Big( D_1' > { {c\over b_0}}\wedge \tau_M\Big)
\ge \P(E'\cap E'') >\,1- \,e^{- b_0/(64\s^2c)}\, . \;\;\diamond $$  
\par \medskip 

   Let us apply now the preceding lemma, to get jointly a lower bound on all proper
hitting times at which an excursion in the hole begins. 
\begin{lem}\label{lem.estimDk}\quad  Fix $\,\e>0$, and recall that $\;\tau_M:=\inf\{ s\,|\,r_s= M\}$. Then
there exists a lower bound $\,b(\s,R,M,\e)\,$ such that for $\,b_0\ge b(\s,R,M,\e)\,$ we have : 
$$ \P\Big( D_{3n-1} \ge \tau_M\wedge \Big( {\ts{1\over (1+\e)\s^2}}\log(1+({\ts{(1+\e)
c\,\s^2\over  b_0^{1+(\e/2)}}})\, n)\Big)\;\hbox{ for all }\; n\in\N\Big) >1-\e\,. $$ 
\end{lem}
\ub{Proof}\quad  Set $\;\d_k:= e^{(1+\e)\s^2\,D_{3k+2}}\,$, for $\,k\in\N\,$, \ and consider the
following events, indexed by the integer $\,n\ge 0\,$ : 
$$ B_n := \Big\{ b_0^{1-(\e/2)}\, e^{(1-\e)\s^2s}\le b_s \le b_0^{1+(\e/2)}\, e^{(1+\e)\s^2s} \quad\hbox{
for all }\; s\le D_{3n-1} \Big\} , $$ 
$$  A_{n}:=\Big\{ D_{3k-1} \ge \tau_M\wedge \Big( {\ts{1\over (1+\e)\s^2}} \log(1+({\ts{(1+\e)
c\,\s^2\over b_0^{1+(\e/2)}}})\,k)\Big) \;\;\hbox{ for all integers }\;  k\le n \Big\} , $$ 
and 
$$  A_n' := \Big\{ D_{3n-1} \ge {\ts{1\over (1+\e)\s^2}} \log(1+({\ts{(1+\e) c\,\s^2\over 
b_0^{1+(\e/2)}}})\,n) \Big\}\; ,\quad B_n' := B_n\cap \{\tau_M\ge D_{3n-1}\} . $$  
Finally let us consider also 
$$ B :=  \Big\{ b_0^{1-(\e/2)}\, e^{(1-\e)\s^2s}\le b_s \le b_0^{1+(\e/2)}\, e^{(1+\e)\s^2s}
\quad\hbox{ for all }\; s\le\tau_M\Big\} . $$ 
We have 
$$ A'_{n+1} = \Big\{ D_{3n+2} \ge {\ts{1\over (1+\e)\s^2}}\log(1+({\ts{(1+\e) c\,\s^2\over
b_0^{1+(\e/2)}}}) (n+1))\Big\}  = \Big\{ \d_n \ge 1+ ({\ts{(1+\e) c\,\s^2\over b_0^{1+(\e/2)}}}) (n+1)
\Big\} $$ 
$$ \supset \Big\{  \d_n\ge \d_{n-1}+ {\ts{(1+\e) c\,\s^2\over b_0^{1+(\e/2)}}}\Big\} \bigcap 
\Big\{ \d_{n-1} \ge 1+({\ts{(1+\e) c\,\s^2\over b_0^{1+(\e/2)}}})\,n \Big\} $$ 
$$ \supset \Big\{  (D_{3n+2}-D_{3n-1})\,\d_{n-1} \ge {c\over b_0^{1+(\e/2)}}\Big\} \bigcap 
\Big\{ \d_{n-1} \ge 1+({\ts{(1+\e) c\,\s^2\over b_0^{1+(\e/2)}}})\,n \Big\} , $$ 
and 
$$ \Big\{ (D_{3n+2}-D_{3n-1})\,\d_{n-1} \ge {c\over b_0^{1+(\e/2)}} \Big\} \supset 
\Big\{ D_{3n+2}-D_{3n-1} \ge {c\over b_{D_{3n-1}}}\Big\} \bigcap B_n\, , $$
whence, using that $\; A_{n+1}= (A'_{n+1}\cap A_n) \cup (\{D_{3n+2} >\tau_M\}\cap A_n)\,$ :  
$$ A_{n+1}\cap B'_n\;\supset\,\Big\{ D_{3n+2}-D_{3n-1} >{c\over b_{D_{3n-1}}}\;\hbox{ or }\;
D_{3n+2} >\tau_M\Big\}\bigcap A_n \cap B'_n\,.  $$  
 
   Now by Lemma \ref{lem.control1} and the strong Markov property we have (denoting by $\,\FF_\cdot\,$
the natural filtration of the diffusion, and setting $\,\a:= (64 \s^2c)\1$) : 
$$ \P\Big( D_{3n+2}-D_{3n-1} \ge {c\over b_{D_{3n-1}}}\;\hbox{ or }\; D_{3n+2} >\tau_M \,\Big|\,
\FF_{D_{3n-1}}\Big) \,\ge \, \E\Big( 1-e^{-\a\,b_{D_{3n-1}}}\,\Big|\,\FF_{D_{3n-1}}\Big)  . $$ 
Hence 
$$ \P(A_{n+1}\cap B'_n)\,\ge\, \P\Big(\Big\{ D_{3n+2}-D_{3n-1} \ge {c\over b_{D_{3n-1}}}
\;\hbox{ or }\;D_{3n+2} >\tau_M\Big\}\bigcap A_n\cap B'_n\Big) $$ 
$$ \ge\, \E\Big[\Big( 1-e^{-\a\,b_{D_{3n-1}}}\Big)\times 1_{ A_{n}\cap B'_n}\Big] \ge 
\E\Big[\Big( 1-e^{-\a\,b_0^{1-(\e/2)}\,\d_{n-1}^{({1-\e\over 1+\e})}}\Big)\times 1_{A_{n}\cap
B'_n}\Big] $$ 
$$ \ge \Big( 1-\exp\Big[{-\a\,b_0^{1-(\e/2)}\,((1+\e) c\,\s^2\,n/b_0^{1+(\e/2)})^{1-\e\over
1+\e}}\Big]\Big) \times \P( A_{n}\cap B'_n) $$ 
$$ = \Big( 1-\exp\Big[{-{\ts{1\over 64}}(\s^2c)^{-2\e\over1+\e}
(1+\e)^{1-\e\over1+\e}\,b_0^{\e(2-\e)\over1+\e}\, n^{1-\e\over 1+\e}}\Big]\Big) \times \P ( A_{n}\cap
B'_n)\, . $$  Set $\;\e_n:= \exp\Big[{-{\ts{1\over 64}}(\s^2c)^{-2\e\over1+\e}
(1+\e)^{1-\e\over1+\e}\,b_0^{\e(2-\e)\over1+\e}\, n^{1-\e\over 1+\e}}\Big]\,$, for $\,n\ge 1\,$, and
$\;\e_0:= e^{-\a\,b_{0}}\,$, so that we have clearly, for $\,b_0\ge b(\s,R,\e)$, on one hand \ 
\if{ (setting $\;\kappa := {\ts{1\over 64}}(\s^2c)^{-2\e\over1+\e}(1+\e)^{1-\e\over1+\e}$) :  
$$\sum_{n\ge 0}\limits\e_n < \e_0+\int_0^\ii \!\exp\!\Big[{-\kappa\,b_0^{\e(2-\e)\over1+\e} 
t^{1-\e\over 1+\e}}\Big] dt = e^{-\a\,b_{0}}+ b_0^{-\e(2-\e)\over1-\e}
\!\int_0^\ii\! e^{-\kappa\,t^{1-\e\over 1+\e}}dt < b_0^{-\e(2-\e)} <\e\,, $$ 
}\fi   $\;\sum_{n\ge 0}\limits\e_n < b_0^{-\e(2-\e)} <\e , $ \ 
and on the other hand :
$$ \P(A_{n+1}\cap B'_n)\,\ge\, (1-\e_n)\times \P(A_{n}\cap B'_n)\,\ge\,\P(A_{n}\cap B'_n)
-\e_n\,, $$  or equivalently 
$$ \P\Big( A_{n}\cap B'_n\cap (A_{n+1})^{\bf c}\Big) \le \e_n \, . $$ 
Note that $\;B'_n\supset\{\tau_M\ge D_{3n-1}\}\bigcap B\,$.
\quad  Therefore 
$$ \P(A_{n+1}\cap B) \,=\, \P\Big(A_n\cap\{ D_{3n-1} >\tau_M\}\cap B\Big) + \P\Big( A_{n+1}\cap
\{\tau_M\ge D_{3n-1}\}\cap B\Big) $$ 
$$ = \P\Big(\! A_n\cap\{ D_{3n-1}\! >\!\tau_M\}\cap B\Big) + \P\Big(\! A_{n}\cap \{\tau_M\!\ge \!
D_{3n-1}\}\cap B\Big) - \P\Big(\! A_{n}\cap\{\tau_M\!\ge\! D_{3n-1}\}\cap B\cap (A_{n+1})^{\bf c}\Big) $$ 
$$  \ge\, \P\Big(A_{n}\cap \{ D_{3n-1} >\tau_M\}\cap B\Big) + \P\Big( A_{n}\cap \{\tau_M\ge D_{3n-1}\}\cap
B\Big) - \P\Big( A_{n}\cap B'_n\cap (A_{n+1})^{\bf c}\Big) $$ 
$$  \ge\, \P\Big(A_{n}\cap \{ D_{3n-1} >\tau_M\}\cap B\Big) + \P\Big( A_{n}\cap \{\tau_M\ge D_{3n-1}\}\cap
B\Big)  - \e_n \,=\, \P(A_{n}\cap B) -\e_n \, , $$ 
whence finally by Lemma \ref{lem.estim2b}, for $\,b_0\ge b(\s,R,\e)\,$ : 
$$\P\Big(\bigcap_{n\ge 0}\limits A_n\Big) \,\ge\; \P(B) -\sum_{n\in\N}\limits\e_n >
1-2b_0^{-\e^2/2}-b_0^{-\e(2-\e)}>1-\e\;.\;\;\diamond $$  
\par\medskip 

   {\bf 6)} \ {\it Estimates related to} $\,|a_\cdot|/b_\cdot\,$ \medskip 

   We need next to control the integral $\;{\ds\int{dt\over r_t\,b_t^2}}\;$, which occurs in the
It\^o expression of the crucial quantity $\,{\ds{a_s\over b_s}}\,$. \ 
The following lemma estimates its contribution due to an excursion in the hole.
\begin{lem}\label{lem.control2}\quad  During any excursion in the hole, id est during any  
proper time interval $\,[D_{3n},D_{3n+2}]\,$, we have the following control : \quad 
${\ds \int_{D_{3n}}^{D_{3n+2}}{dt\over r_t}\,\le\, {\pi\,R\over \min_{[D_{3n},D_{3n+2}]}\limits
b}\,}\,$. 
\end{lem}
\ub{Proof}\quad  Firstly, it is sufficient to consider the proper time interval $\,[D_{3n},D_{3n+1}]\,$,
since the estimates are exactly the sames on the other half $\,[D_{3n+1},D_{3n+2}]\,$. \ Since \parn 
$ T_t\,r_t^{3/2} = - \sqrt{(R-r_t)b_t^2+Rr_t^2+(a_t^2-1)r_t^3}\,\le - \sqrt{R-r_t}\; b_t\; \hbox{ for }
D_{3n}\le t\le D_{3n+1}\, ,$  \  we get 
$$ \int_{D_{3n}}^{D_{3n+1}}{dt\over r_t}\,\le\, -\int_{D_{3n}}^{D_{3n+1}} \sqrt{r_t\over R-r_t}\times 
{d\,{r_t}\over b_t} \le {1\over\min_{} b[D_{3n},D_{3n+2}]}\times \int_0^R \sqrt{r\over R-r}\,\,dr $$
$$ =\; {\pi\,R\over 2\,\min_{} b[D_{3n},D_{3n+2}]}\; . \;\;\diamond $$ 

   \begin{rem}\label{rem.control2}\quad  {\rm The same argument yields also the
following estimate on the duration of an excursion in the hole (showing that the estimate from below
in Lemma \ref{lem.control1} is essentially optimal) : 
$$ D_{3n+2}-D_{3n}=\int_{D_{3n}}^{D_{3n+2}}{dt}\,\le {2\over\min_{} b[D_{3n},D_{3n+2}]}\times \int_0^R
\sqrt{r\over R-r}\,\,r\, dr = {3\pi\,R^2\over 4\,\min_{} b[D_{3n},D_{3n+2}]}\; . $$  }
\end{rem}

   We can now deduce the control on the integral $\;{\ds\int{dt\over r_t\,b_t^2}}\;$, which we
shall need below. 
\begin{lem}\label{lem.int1/rb2}\quad For any $\,\e\in ]0,\5[\,$, \ there exists a lower bound
$\,b(\s,R,M,\e)\,$ such that for $b_0\ge b(\s,R,M,\e)\,$ we have : 
${\ds \;\int_0^{\tau_M}{dt\over b_t^2}+\int_0^{\tau_M} {R\,dt\over r_t\, b_t^2} \le b_0^{2\e-2}}$, with
probability larger than $\,1-\e\,$. 
\end{lem}
\ub{Proof}\quad  Recall from Lemma \ref{lem.estim2b} that with large probability we have
$\;b_s\ge b_0^{1-(\e/2)}\, e^{(1-\e)\s^2s}\;$ for $\,s\le\tau_M\,$,  whence 
$\;{\ds\int_{0}^{\tau_M} b_t\2dt\le {b_0^{\e-2}\over 2(1-\e)\s^2}}\,$. \ 
Moreover, using Lemma \ref{lem.control2} we have  : 
$$ \int_0^{\tau_M} {R\, dt\over r_t\, b_t^2} - \int_0^{\tau_M} {dt\over b_t^2} \,
< \, R\sum_{D_{3n}<\tau_M} \int_{D_{3n}}^{D_{3n+2}} {dt\over r_t\, b_t^2} \, 
<\, \pi R^2 \sum_{D_{3n}<\tau_M} ( \min b[D_{3n}, D_{3n+2}] )^{-3} $$
$$ <  \pi R^2\, b_0^{{3\e\over 2}-3} \sum_{D_{3n}<\tau_M} e^{-3 (1-\e) \s^2 D_{3n}} \le \pi R^2\,
b_0^{{3\e\over 2}-3} \sum_{D_{3n}<\tau_M} \exp\Big[ -3 ({\ts{1-\e\over 1+\e}}) \log(1+({\ts{(1+\e)
c\,\s^2\over b_0^{1+(\e/2)}}})\,n)\Big] $$  
\centerline{(by Lemma \ref{lem.estimDk}, on an event of probability larger than $\,1-\e\,$)}
$$ <\,\pi R^2\,b_0^{{3\e\over 2}-3} \sum_{n\in\N} (1+n)^{-3 ({\ts{1-\e\over 1+\e}}) } <\,b_0^{2\e-3} 
\quad (\hbox{for } b_0 \hbox{ large enough). } \;\;\diamond $$ 
\par \smallskip   

   We establish then the crucial control on $\,a_s/b_s\,$. 
\begin{lem}\label{lem.rM}\quad Almost surely, if $\,r_\cdot\,$ is bounded, then $\;a_s/b_s\;$ converges as
$\,s\to\ii\,$, to some random limit $\,\ell\in\R\,$. \  Moreover for any $\,\e > 0\,$ and any 
$\,M\ge 3R/2\,$ there exists a lower bound $\,b(\s,R,M,\e)\,$ such that for $\,b_0\ge b(\s,R,M,\e)\,$ we
have : 
$$ \P\Big( \sup_{0\le s< \tau_M}\limits\,\Big|{a_{s}\over b_{s}} - {a_{0}\over b_{0}}\Big|\,< 
({\ts{|a_{0}|\over b_{0}}}+1)\,\e^2 \Big) > 1-\e\; . $$ 
\end{lem}
\ub{Proof}\quad  $a)$ \ Let us prove the second assertion first. \par 
   Set \ $\;\b:=\inf\Big\{ s>0\,\Big|\,b_s<b_0^{1-{\e\over 2}} e^{(1-\e)\s^2s}\Big\}$, and
\ $\;\eta:=\inf\Big\{ s>0\,\Big|\,\Big|{a_{s}\over b_{s}} - {a_{0}\over b_{0}}\Big| > \e^2\Big\}$. \parn  
Fix $\,0<\e<1/3\,$ and $\,M\ge 3R/2\,$. Recall that we set $\;\tau_M:=\inf\{ s\,|\,r_s= M\}$, \ and that,
for $\,b_0\ge b(M,\e)\,$, \ the event \ $E_1:=\{ \b\ge \tau_M\}\;$ has probability at least $\,1-\e\;$ by
Lemma \ref{lem.estim2b}, and that the event \ ${\ds E_2:=\Big\{\int_0^{\tau_M} {dt\over
b_t^2}+\int_0^{\tau_M}{R\,dt\over r_t\, b_t^2} \le b_0^{2\e-2}\Big\}\, }$ has probability at least
$\,1-\e\,$ by Lemma \ref{lem.int1/rb2}. 
\if{ Moreover \  the event $\;E_3:=\Big\{ W^*[({\ts{|a_{0}|\over
b_{0}}}+1)^2 b_0^{3\e-2}]< ({\ts{|a_{0}|\over b_{0}}}+1) b_0^{2\e-1}\Big\}\,$  has clearly probability at
least $\,1-\e\,$, for any  Brownian motion $\,W\,$ and for $\,b_0\ge b(\e)\,$. \quad }\fi 
   Now, the equation governing $\,{a_{s}/ b_{s}}\,$ writes for some Brownian motion $\,W\,$
and for any $\,s\ge 0\,$ : 
$$ {a_{s}\over b_{s}} = {a_{0}\over b_{0}} + 
W\Big[\s^2\!\int_{S}^{s} ({a_t\over b_t})^2({r_t\over b_t})^2dt + 
\s^2\!\int_{0}^{s} \Big({R\over r_t\,b_t^2}-{1\over b_t^2}\Big)\,dt \Big] +
{\ts{\s^2\over 2}} \int_{0}^{s} ({a_t\over b_t})({r_t\over b_t})^2dt \, , $$ 
whence  
$$ {|a_{s}|\over b_{s}} - {|a_{0}|\over b_{0}}\le \Big|{a_{s}\over b_{s}} - {a_{0}\over b_{0}}\Big| \le  
W^*\Big[\s^2\int_{0}^{s} \Big|{a_t\, r_t\over b_t^2}\Big|^2dt + 
\s^2\int_{0}^{s} {R\over r_t\,b_t^2}\,dt \Big] +
{\ts{\s^2\over 2}} \int_{0}^{s} \Big|{a_t\over b_t}\Big|({r_t\over b_t})^2dt \, , $$ 
\centerline{(recall that $\,W^*\,$ denotes the two-sided maximum process of $\,W$ :
$W^*(s):=\max |W|([0,s])\,$) }\parn 
so that we have on $\, E_1\cap E_2\,$, for $\,b_0\ge b(\s,R,M,\e)\,$ and $\,0\le s <
\tau_M\wedge\eta\,$ : 
$$ \Big|{a_{s}\over b_{s}} - {a_{0}\over b_{0}}\Big| \le 
W^*\Big[ \s^2M^2 ({\ts{|a_{0}|\over b_{0}}}+\e^2)^2\!\int_{0}^{\tau_M} {dt\over b_t^2} +
\s^2\!\int_{0}^{\tau_M}\! {R\over r_t\,b_t^2}\,dt \Big] + {\ts{\s^2M^2\over 2}}\,({\ts{|a_{0}|\over
b_{0}}}+\e^2)\int_{0}^{\tau_M} {dt\over b_t^2} $$
$$ < W^*[ ({\ts{|a_{0}|\over b_{0}}}+1)^2 b_0^{3\e-2}] + ({\ts{|a_{0}|\over b_{0}}}+1)\,b_0^{3\e-2}
< 2\,({\ts{|a_{0}|\over b_{0}}}+1)\, b_0^{2\e-1} < \e^2 ({\ts{|a_{0}|\over b_{0}}}+1)\, , $$  
on an event $\,E_3\,$ of probability larger than $\,1-\e\,$, for $\,b_0\,$ large enough. \par 

   Hence we find that for $\,b_0\ge b(\s,R,M,\e)\,$ we have $\;\eta\ge\tau_M\,$ on $\,E_1\cap E_2\cap E_3\,$, whence :  
$$ \P\Big( \sup_{0\le s< \tau_M}\limits\,\Big|{a_{s}\over b_{s}} - {a_{0}\over b_{0}}\Big|\,<
({\ts{|a_{0}|\over b_{0}}}+1)\, \e^2\,\Big) > 1-3\e\; . $$ 

    $b)$ \ To prove the  first assertion of the statement, let us onsider then, for any $\,k\in\N^*$, the hitting time $\;\s_k:=\inf\{ s>0\,|\,b_s=b(\s,R,k,2^{-k})\}$, which is almost surely finite by Lemma
\ref{lem.bdiv}. Moreover, since $\,r_\cdot\,$ is bounded, there exists almost surely some random integer 
$\,\kappa\,$ such that $\,\tau_{\kappa}=\ii\,$. Applying the strong Markov property to the above, we get 
$$ \sum_{k\in\N^*}\P\Big( \sup_{\s_k\wedge\tau_k\le s< \tau_k}\limits\,\Big|{a_{s}\over b_{s}} -
{a_{\s_k}\over b_{\s_k}}\Big|\,\ge 4^{-k}\,({\ts{|a_{\s_k}|\over b_{\s_k}}}+1)\Big) \le 3\sum_{k\in\N^*}
2^{-k}<\ii\; , $$  showing by the Borel-Cantelli lemma that almost surely, there exists some random
integer $\,\kappa'\,$ such that $\;\sup_{\s_k\wedge\tau_k\le s< \tau_k}\limits\,\Big|{a_{s}\over b_{s}} -
{a_{\s_k}\over b_{\s_k}}\Big|\,< 4^{-k}\,({\ts{|a_{\s_k}|\over b_{\s_k}}}+1)\,\,$ for $\,k\ge\kappa'\,$.
Hence we get almost surely 
$\;\sup_{s\ge\s_k}\limits\,\Big|{a_{s}\over b_{s}} -{a_{\s_k}\over b_{\s_k}}\Big|\,<
4^{-k}\,({\ts{|a_{\s_k}|\over b_{\s_k}}}+1)\,\,$ for $\,k\ge k_0:=\kappa\vee\kappa'$. This implies at once
that 
$\;\sup_{s\ge\s_{k_0}}\limits\,\Big|{a_{s}\over b_{s}}\Big|\,<\,
(2\,{\ts{|a_{\s_{k_0}}|\over b_{\s_{k_0}}}}+1)\,$, and then that 
$\;\sup_{s\ge\s_k}\limits\,\Big|{a_{s}\over b_{s}} -{a_{\s_k}\over b_{\s_k}}\Big|\,<
4^{-k}\,(2\,{\ts{|a_{\s_{k_0}}|\over b_{\s_{k_0}}}}+1)\,\,$ for $\,k\ge k_0\,$, thereby showing that
$\,a_s/b_s\,$ satisfies the Cauchy criterion, and then converges as $\,s\to\ii$. $\;\diamond$ 
\par\medskip \medskip 

   {\bf 7)} \ End of description of the second case in Theorem \ref{the.piege}. \par 

\begin{lem}\label{lem.rr}\quad  Set $\;{\ds\rr:=\limsup_{s\to\ii}\limits r_s}\,$. Then almost surely, if $\;r_\cdot\,$ is bounded, $\;a_\cdot/b_\cdot\,$ converges to $\;\pm\,{1\over\rr}\,\sqrt{1-{R\over\rr}}\,$.
\end{lem}
\ub{Proof}\quad  Lemma \ref{lem.rM} insures that $\;a_s/b_s\,$ converges to some $\,\ell\in\R\,$ as
$\,s\to\ii$, and we can find a sequence $\,s_n\,$ of proper
times  increasing to $\ii\,$ such that $\;\rr =\lim_{n\to\ii}\limits\, r_{s_n}\;$ and $\;T_{s_n}=0\,$ for
every $\,n$, which by the pseudo-norm relation implies : \quad
${\ds {a^2_{s_n}\over b^2_{s_n}} = (1-{\ts{R\over r_{s_n}}})(r\2_{s_n}+b\2_{s_n})}\,$, so that using 
Lemma \ref{lem.bdiv} we find as wanted : $\,\ell^2=(1-{\ts{R\over \rr}})\rr\2$. $\;\diamond$ 
\par\medskip \medskip

   We prove next that when the relativistic diffusion becomes eventually captured by a neighbourhood of the hole, it tends to stay in some asymptotic random plane of $\,\R^3\,$. 
\begin{lem}\label{lem.planas}\quad  Almost surely, if $\,r_\cdot\,$ is bounded, the direction
$\,\vec{b}_s/b_s\,$ of the angular momentum $\,\vec{b}_s\,$ converges in $\,\S^2\,$ as $\,s\to\ii\,$ :
the trajectories are asymptotically planar. 
\end{lem} 
\ub{Proof}\quad As we already saw in the third part of the proof of Theorem \ref{the.martp}, it is
easily deduced from Proposition \ref{pro.systeqst} that  
$$ d\Big({\vec{b}_s\over b_s}\Big) = d\Big(\t_s\wedge{\dot\t_s\over U_s}\Big) = -\Big({r_s\over
b_s}\,\s\,d\b_s\Big)\, {\dot\t_s\over U_s} - \Big({\s^2\,r_s^2\over 2\,b_s^2}\,ds\Big)\,
\t_s\wedge{\dot\t_s\over U_s}\; . $$  
Then for any $\,k\in\N\,$ and $\,s\ge k\,$ we have (for some standard Brownian motion $\,B_k$) :
$$ \Big|{\vec{b}_s\over b_s} - {\vec{b}_k\over b_k}\Big| = \sup_{u\in\S^2} \int_k^s \Big[ 
{r_t\over b_t}\Big({\dot\t_s\over U_s}\cdot u\Big)\s\,d\b_t + {\s^2\,r_t^2\over
2\,b_t^2}\Big({\vec{b}_t\over b_t}\cdot u\Big) dt \Big] \le 
\s B_k^*\Big[\int_k^\ii {r_t^2\over b_t^2}\,dt\Big] + \int_k^\ii {\s^2r_t^2\over 2\,b_t^2}\,dt \, . $$ 

   Now Lemma \ref{lem.bdiv} insures that $\;{\ds\int_k^\ii {r_t^2\over b_t^2}\,dt\le e^{-\s^2k/2}}\;$
for $\,k\,$ larger than some finite random $\,\kappa\,$, and the Borel-Cantelli lemma insures that 
$\;B_k^*(e^{-\s^2k/2})\le e^{-\s^2k/6}\;$ for $\,k\,$ larger than some finite random $\,\kappa'\,$. 
Hence for $\,k\ge \kappa\vee\kappa'\,$ we get  
$\;\sup_{s\ge k}\limits \Big|{\vec{b}_s\over b_s} - {\vec{b}_k\over b_k}\Big| \le \s e^{-\s^2k/6}
+ \s^2e^{-\s^2k/2}$, showing that $\,\vec{b}_s/b_s\,$ satisfies almost surely the Cauchy criterion and
thus converges in $\S^2$. $\;\diamond$ 
\par\medskip \smallskip 

   {\bf 8)} \ End of the proof of Theorem \ref{the.piege}. \par \medskip 

    The following lemma proves the statement $(ii)$ of Theorem \ref{the.piege}. 
\begin{lem}\label{lem.rrr}\quad  Fix $\,\e\in ]0,1[$, and suppose $\,r_0\in ]R,3R/2[\,$ and $\,T_0=0$.
Fix $\,\a,\b\,$ such that $R\le\a < r_0< \b < 3R/2\,$.  Then  for $\,b_0\,$ large enough, with probability larger than $\,1-\e\,$  \ ${\ds\rr =\limsup_{s\to\ii}\limits r_s}\,$ belongs to $\,]\a,\b[\,$. 
\end{lem}
\ub{Proof}\quad  Set $\;f(r) := r\1\sqrt{1-{\ts{R\over r}}}\;$. The function $\,f\,$ is continuous and
srtictly increasing from $\,{\ds [R,{3R\over 2}]}\,$ onto $\,{\ds [0,{\ts{2\over 3\sqrt{3}\,R}}]}$. 
For $\,b_0\,$ large enough, by the pseudo-norm relation we have \ 
${\ds\Big|{{a_{0}\over b_{0}}}\Big| = f(r_0)\times\sqrt{1+r_0^2/b_0^2}\,\in\, ]f(\a),f(\b)[}\,$, \ and then by Lemma  \ref{lem.rM} the event  \parn 
${\ds E:= \Big\{ f(\a) <  \inf_{0\le s< \tau_{3R/2}}\limits\, \Big|{a_{s}\over b_{s}}\Big|<\sup_{0\le s< \tau_{3R/2}}\limits\, \Big|{a_{s}\over b_{s}}\Big| < f(\b)\Big\}}\,$ has probability larger than $\,1-\e\,$. \ 
Now the pseudo-norm relation implies that $\;{\ds f(r_s)\,1_{\{ r_s\ge R\}} \le\Big|{a_{s}\over b_{s}}\Big|}\,$ for any $\,s\ge 0\,$.  \ So that (by continuity of $\,f\1\,$) $\,r_s\,$ does never hit $\,\b\,$ on $\,E\,$, and then $\,E\subset\{\tau_{3R/2} =\ii\}$. \ Finally by Lemma \ref{lem.rr} and by continuity of $\,f\1$, we have $\,f(\a)<f(\rr)<f(\b)\,$ and then $\,\a <\rr <\b\,$ on $\,E\,$.   $\;\diamond$ 
\par\medskip \medskip 

   The statement $(i)$ in Theorem \ref{the.piege} is seen as in the proof of Theorem (\ref{the.asymp},2).
Indeed, we observed there that if $\;r_0\ge 3R/2\,$, $\,T_0\ge 4+{4\over R\s^2}\,$, and $\;\la
:=\inf\Big\{ s>0\,\Big|\, T_s=2+{2\over R\s^2}\,\Big\}$, then $\;1/|T_{s\wedge\la}|\,$ is a
supermartingale. \ Now this implies, if $\,T_0\,$ is large enough : 
$$ (2+{\ts{2\over R\s^2}})\1\,\P (\la <\ii )\le \liminf_{s\to\ii}\limits\, \E (|T_{s\wedge\la}|\1
1_{\{\la <\ii\}}) \le \liminf_{s\to\ii}\limits\, \E(|T_{s\wedge\la}|\1) \le T_{0}\1\,, $$  
and then that $\;\;\P (\lim_{s\to\ii}\limits r_s = +\ii ) \ge \P (\la =\ii ) \ge 1-(2+{\ts{2\over
R\s^2}})/T_0 > 1-\e\,$.  \par \smallskip 

   Then the irreducibility of the relativistic diffusion follows at once from the proof of Theorem (\ref{the.asymp},2). With Statements $(i)$ and $(ii)$, this proves that both cases in Theorem \ref{the.piege} indeed occur with strictly positive probability, from any initial condition. Of course we already proved the dichotomy in Lemma \ref{lem.dicho}.  \par   

   Finally the descrition of the first case in Theorem \ref{the.piege} follows directly from Lemma
\ref{lem.cv} and from Theorem (\ref{the.asymp},3). $\;\diamond$

\subsection{Proof of Corollary \ref{cor.piege}}\indf \label{sec.PC3}
   Let us begin by a lemma which specifies at which times the relativistic diffusion approaches 
the top of its excursions outside the hole. 
\begin{lem}\label{lem.piege}\quad  Let us consider the following stopping times, for any $\,n\in\N\,$ 
and any $\,\e>0\,$ :  \parn  
\centerline{$ D_{n}':=\min\{ s>D_{3n+1}\,|\, T_s=0\}\,$ \  and \ $\; D_{n}^\e:=\min\{ s>D_{3n+1}\,|\, T_s\le \e\,b_s\} < D_{n}'$.}  \parn 
 Then almost surely, when $\,r_\cdot\,$ is bounded, we have :  $\;\lim_{n\to\ii}\limits\,r_{D_{n}'} = \rr\,$, and  $\;\lim_{n\to\ii}\limits\,r_{D_{n}^\e} = \rr_\e <\rr\,$,  \parn 
where $\,\rr_\e\,$ is the unique solution between $0$ and ${3R\over 2}$ of the equation 
$\, (1-R/\rr_\e)\,\rr_\e\2 = \ell^2-\e^2$. 
\end{lem}
\ub{Proof}\quad This is only an additional precision to the second case in Theorem \ref{the.piege} : 
applying the pseudo-norm relation at time $\; D_{n}'\,$, we get \parn 
\centerline{$ (1-R/r_{D_{n}'})\, r\2_{D_{n}'}\sim (1-R/r_{D_{n}'})(r\2_{D_{n}'}
+b\2_{D_{n}'}) = a^2_{D_{n}'}/b^2_{D_{n}'} \lra \ell^2\,$ as $\,n\to\ii\,$,} \parn 
so that, by Lemma \ref{lem.limsup}, owing to the function $\,f\,$ used in the proof of Lemma
\ref{lem.rrr}, we must have the unique possibility : $\,\lim_{n\to\ii}\limits r_{D_{n}'} =\rr\,$. \ Similarly, 
at time $\; D_{n}^\e\,$ we get \parn
\centerline{$ (1-R/r_{D_{n}^\e})\, r\2_{D_{n}'}\sim (1-R/r_{D_{n}^\e})(r\2_{D_{n}^\e}
+b\2_{D_{n}^\e}) = {a^2_{D_{n}^\e}\over b^2_{D_{n}^\e}} - {T^2_{D_{n}^\e}\over b^2_{D_{n}^\e}}
\lra \ell^2-\e^2\,$ as $\,n\to\ii\,$,} \parn 
so that we must have the unique possibility : $\,\lim_{n\to\ii}\limits r_{D_{n}^\e} =\rr_\e\,$. 
$\;\diamond$   \par\medskip 

   We shall need to control the angular contribution around the top of the excursions outside the hole.  This is the aim of the following lemma. \ Notations are as above. 
\begin{lem}\label{lem.intb}\quad  We have  \quad  ${\ds  \limsup_{n\to\ii} \limits\, \int_{D_{n}^\e}^{{D_{n}'}} b_s\,ds \,\le  \,{\ts{12\over 3R-2\rr}}\, ({\ts{3R\over \sqrt{2}}})^4\,\e\, }$ on the intersection of the event $\,\{\rr <{3R\over 2}\}\,$ and of an event of probability larger than $\,1-\e\,$. 
\end{lem}
\ub{Proof}\quad   Set $\,E_0 := \{\rr <{3R\over 2}\}$. Fix $\,M=M(\e)\ge {3R\over 2}\,$ such that the event $\,E_0\sm \{ \sup r_\cdot < M\}\,$ has probability smaller than $\,\e/8\,$. Set 
$\,E_1=E_1(\e):= \{ \sup r_\cdot < M\}$. By Lemma \ref{lem.bdiv} we can find some $\,\g=\g(\e)\,$ such that, setting $\,E_2=E_2(\e):= \{ b_s  > \g\,e^{(1-\e)\s^2 s}\;\hbox{ for all } s\ge 0\}$,  the event $\,E_0\sm E_2\,$ has probability smaller than $\,\e/8\,$. Consider then the stopping time $\,\nu=\nu(\e)\,$ at which $\,b_\cdot$ hits the lower bound $\,b(\s,R,M,\e/4)\,$ of Lemma \ref{lem.estimDk}, which is finite on $\,E_2\,$. Applying now Lemma \ref{lem.estimDk} and the strong Markov property (applied at time $\,\nu$), we find some integer $\,k=k(\e)\,$ such that, setting  $\;E_3=E_3(\e) := \bigcap_{n>k} \Big\{ D_{3n+1} \ge {\ts{1\over (1+\e)\s^2}} \log( n/k)\Big\} $ and $\,E_4=E_1\cap E_2\cap E_3\,$, the event $\,E_0\sm E_4\,$ has probability smaller than $\,\e/2\,$. \par 
   Set $\,\b_n:= b_{D_{3n+1}}\,$.  \ We have obviously \
$$ E_4\subset B = B(\e) := \bigcap_{n>k} \Big\{ \sqrt{\b_n} \, > \g\,\exp[ {\ts{1-\e\over 2(1+\e)}} \log(n/k) ] \Big\} . $$  

Consider also for all $\,n > k\,$ : \ 
$\;{\ds B_n = B_n(\e) := \bigcap_{m=k+1}^{n} \Big\{ \sqrt{\b_n} > \g\,\exp[ {\ts{1-\e\over 2(1+\e)}} \log(n/k) ]\Big\} }$, and let us proceed now somewhat as in Lemma \ref{lem.control1}, to control the variations of $\,b_\cdot\,$ and $\,T_\cdot\,$ between $\,D_{3n+1}\,$ and $\,D'_{n}\,$.  \par   
   Set $\;\eta_n:= \inf\Big\{ s>D_{3n+1}\,\Big|\,|\log (b_s/\b_n)|=\log 2\Big\}$, and 
$\;\tau := \inf\{ s>0\,|\, r_s=M\}$ (which is infinite on $\,E_4$). \ 
For $\,D_{3n+1}\le s\le (D_{3n+1}+\b_n^{-1/2})\wedge \eta_{n}\wedge \tau\;$ and $\,n\,$ large enough, we have : 
$$  \Big|\log (b_s/\b_n)\Big| \le \s\, W_n^*\Big[ (1+{\ts{4\,M^2\over \b_n^2}})\,{\ts{1\over \b_n^{1/2}}}\Big] + {\ts{\s^2\over  \b_n^{1/2}}} < \5 + {\ts{\s^2\over  \b_n^{1/2}}} <\log 2\, \; \hbox{ on } \;
B_n\cap F_n\, , $$  
where $\;F_n:= \Big\{ \s W_n^*[ 2 \b_n^{-1/2}\,] < \5 \Big\} $. \   Whence \quad 
$\,B_n\cap F_n\subset \{ \eta_n > (D_{3n+1}+\b_n^{-1/2})\wedge\tau\}$. \parn  
Recalling Remark \ref{rem.control2}, we observe then that  for large enough  $\,n\,$ we have \parn 
${\ds  D_{3n+2}\wedge (D_{3n+1}+ \b_n^{-1/2})\wedge\tau \le D_{3n+1}+ {{3\pi R^2\over 4\,\b_n}} <  D_{3n+1}+ \b_n^{-1/2}\;}$ on $\,B_n\cap F_n\,$, whence for any constant $\,\a$ : \ 
${\ds B_n\cap F_n\subset \{\eta_n> (D_{3n+2}+\a/\b_n)\wedge\tau \} }$ for large enough  $\,n\,$.  \medskip 

   Then for any deterministic $\,q,j \in\N^*$, set $\, A_{q,j} :=\{( j+1) 2^{-q} < {3R\over 2}-\rr\} \subset E_0\,$, and  fix a deterministic integer $\,m\,$ such that each $\; A_{q,j}^m := A_{q,j}\cap \Big\{ \sup_{s > D_{3m+1}} \limits r_s < \rr + j\,2^{-q}\Big\}$ satisfies  \parn  
$\P(A_{q,j}\sm A_{q,j}^m)< {2^{-q} \e\over 3R}\,$.  \ Set  $\;\tau_{q,j}^m := \inf\{ s> D_{3m+1}\,|\, r_s= {3R\over 2}- {j \over 2^{q}}\}$ (which is infinite on $\,A_{q,j}^m$). \par\medskip  

\if{
   For any deterministic $\,h>0\,$, set $\;A^h :=\{ 3R-2\rr> 4h\} \subset E_0\,$, and  fix a deterministic integer $\,m\,$ such that $\; A_m^h := \Big\{ \sup_{s > D_{3m+1}} \limits r_s < \rr +h < {3R\over 2}-h\Big\}$ satisfies $\P(A^h\sm A_m^h)<\e/2\,$.  \ Set also $\,\tau_m := \inf\{ s> D_{3m+1}\,|\, r_s= {3R\over 2}-h\}$ (which is infinite on $\,A_m^h$). \par\smallskip 
}\fi  

   Consider \  $ \L^\e_n:=\inf\{ s> D_{n}^\e\,|\,T_s\ge \b_n\}$, and fix some $\,\a >0\,$ (to be specified below).  \ The equation governing $\,T_\cdot\,$ implies that for $\,D_{n}^\e \le s\le \eta_n \wedge D'_{n}\wedge  \L^\e_n \wedge (D_{n}^\e+ \a/\b_n^{})\wedge \tau_{q,j}^m\,$  and for $\,n\,$ large enough we have : 
$$ T_s = \e\,b_{D_{n}^\e} + \s\,\bar{W}_n\Big[\int_{D_{n}^\e}^s(T_t^2+1-{\ts{R\over r_t}})\, dt\Big] + {\ts{3\s^2\over 2}}\int_{D_{n}^\e}^s T_t\,dt  +\int_{D_{n}^\e}^s(r_t-{\ts{3R\over 2}})\,{\ts{b_t^2\over r_t^4}}\,dt  - \int_{D_{n}^\e}^s {\ts{R\over 2 r_t}}\,dt   \,  $$  
$$ \le\, 2\e\b_n + \s\,\bar{W}_n^*\Big[(\b_n^2+1)\, \a/\b_n\Big] + {\ts{3\s^2\over 2}}\, \a
 \,<\, 2\e\b_n + \b_n ^{3/4}   + {\ts{3\s^2\over 2}}\, \a \,<\, \b_n $$ 
on  $\, F_n' :=\Big\{ \s \bar{W}_n^*[ 2\a\,\b_n] < \b_n ^{3/4} \Big\} $, so that 
$\; F_n'\subset \{ \L_n^\e > \eta_n \wedge D_{n}'\wedge (D_{n}^\e+\a/\b_n^{})\wedge \tau_{q,j}^m\} $. \par \smallskip 

   The equation of $\,T_\cdot\,$ implies also (in the same way) that, for large enough  $\,n\,$ and  for \  $ D_{n}^\e \le s\le \eta_n \wedge D'_{n}\wedge  \L^\e_n \wedge (D_{n}^\e+ \a/\b_n^{})\wedge \tau_{q,j}^m\,$ : 
$$ T_s \le\, 2\e\b_n + \s\,\bar{W}_n^*\Big[2 \a\,\b_n\Big] + {\ts{3\s^2\over 2}}\, \a -  
{{j\over 2^q}}\, ({\ts{2\over 3R}})^4\, {\b_n^2\over 4}\, (s- D_{n}^\e )  \,, $$  
while on $\,F'_n\,$, for $\, s= D_{n}^\e+ \a/\b_n\,$, taking $\;\a := {\ts{3\e 2^q\over j}}\times ({\ts{3R\over \sqrt{2}}})^4$, the right hand side above is 
$\,< 0\,$ for $\,n\,$ large enough. \quad  Hence we find that \par \smallskip 
\centerline{ $ B_n\cap F_n\cap F_n'\subset \{ D'_n\wedge \tau_{q,j}^m < D_{n}^\e+ \a/\b_n \} $, \ 
for large enough  $\,n\,$.} \par \smallskip 

   Therefore so far we get on $\, B_n\cap F_n\cap F_n'\,$, almost surely for large enough  $\,n\,$ : 
$$  \int_{D_{n}^\e\wedge \tau_{q,j}^m }^{{D_{n}'}\wedge \tau_{q,j}^m } b_s\,ds \,\le\, 2\a = {{6\, \e\, 2^q\over j}}\, ({\ts{3R\over \sqrt{2}}})^4\,  . $$

   Moreover, since the standard Brownian motions $\,W_n\,$ and $\,\bar{W}_n\,$ are independent from the $\,\s$-field $\,\FF_{D_{3n+1}}$, we have  : 
$$ \P ( F_n^{\bf c}\,|\,\FF_{D_{3n+1}} ) \, <\, e^{-\sqrt{\b_n}/(16\s^2)} \quad \hbox{ and } 
\quad \P ( (F'_n)^{\bf c} \,|\,\FF_{D_{3n+1}} ) \, <\, e^{-\sqrt{\b_n}/(4\s^2\a)}\,. $$  
Therefore 
$$ \sum_{n>k}\, \P ( F_n^{\bf c}\cap B) \, \le\,\sum_{n>k}\, \P ( F_n^{\bf c}\cap B_n \, <\,\sum_{n>k}\, \exp\Big(- {\ts{\g\over 16\s^2}}\, \exp[ {\ts{1-\e\over 2(1+\e)}} \log(n/k) ]\Big) < \ii\,  , $$ 
and the Borel-Cantelli lemma implies that $\,B\,$ is almost surely included in $\,\liminf_n\limits\, F_n\,$. This is obviously the same for $\,F'_n\,$, proving that  
almost surely for any $\,q,j\,$ we have on $\,B\cap A_{q,j}^m\,$, for large enough $\,n\,$ : 
$$  \int_{D_{n}^\e}^{{D_{n}'}} b_s\,ds \,\le\,  {{6\, \e\, 2^q\over j}}\, ({\ts{3R\over \sqrt{2}}})^4\,  . $$ 
Finally, since the event ${\ds E_0\sm \Big\{\limsup_{n\to\ii} \limits\, \int_{D_{n}^\e}^{{D_{n}'}} b_s\,ds \,\le  \,{\ts{12\over 3R-2\rr}}\, ({\ts{3R\over \sqrt{2}}})^4\,\e\Big\} }$ is the increasing union (for $\,q\in\N^*$) of  the events 
$\; {\ds \bigcup_{j=1}^{3R\,2^{q-1}} \limits A_{q,j}\sm \Big\{\limsup_{n\to\ii} \limits\, \int_{D_{n}^\e}^{{D_{n}'}} b_s\,ds \,\le  \,{\ts{6\, \e\, 2^q\over j}}\, ({\ts{3R\over \sqrt{2}}})^4\Big\} }$, \ 
we get 
$$ \P\Big[ E_0\sm \Big\{\limsup_{n\to\ii} \limits\, \int_{D_{n}^\e}^{{D_{n}'}} b_s\,ds \,\le  \,{\ts{12\over 3R-2\rr}}\, ({\ts{3R\over \sqrt{2}}})^4\,\e\Big\} \Big] \le \P[E_0\sm B] + \sup_{q\in\N^*} \sum_{j=1}^{3R\,2^{q-1}}  \P[A_{q,j}\sm A_{q,j}^m]  <  \e \, , $$ 
which concludes the proof. $\;\diamond$ 
\bigskip  

   We can now establish the second assertion  of Corollary \ref{cor.piege}.  (The first assertion is just Lemma \ref{lem.piege}.) \par 
      Set $\; g(r):=R-r+\ell^2r^3\,$, so that $\,g(\rr) = 0\,$ and 
$\,g'(r)=(\sqrt{3}\,|\ell| r+1)(\sqrt{3}\,|\ell| r-1)\,$. \parn 
We have $\,\rr\le{3R\over 2} \LRa \rr\le {1\over|\ell|\sqrt{3}}\,$, and $\,\rr ={3R\over 2} \LRa \rr =
{1\over|\ell|\sqrt{3}}\,$. \ So that the root $r=\rr$ of $\,g\,$ is double if and only if 
$\,\rr ={3R\over 2}\,$. So that the integral $\;{\ds\int_0^{\rr}{dr\over\sqrt{r\, g(r)}}}\;$
converges if and only if $\,\rr < {3R\over 2}\,$, thereby justifying our assumption  $\,\rr <{3R\over 2}\,$. \parn 

   Then the the pseudo-norm relation and Lemmas \ref{lem.bdiv} and \ref{lem.rr} imply the convergence of \ 
$ {T_s^2\over b_s^2} + (1-{R\over r_s})\,r_s\2\;$ towards $\,\ell^2=(1-{R\over \rr})\,\rr\2\,$. \ By the definition of $\,b\,$, this can be written : 
$$ U_s\,ds = {\rm sgn}(T_s)\Big(r_s\times\Big[R-r_s+(\ell^2+\e_s)r_s^3\Big]\Big)^{-1/2} dr_s\, , $$
with $\,\lim_{s\to\ii}\limits\,\e_s=0\,$. \ Lemma \ref{lem.planas} means that we can restrict to the
limit plane orthogonal to $\,\vec{\b} :=\lim_{s\to\ii}\limits\,\vec{b_s}/b_s\in\S^2$. Indeed, we have 
$\; \dot\t_s/U_s = (\vec{b_s}/b_s) \wedge \t_s\,$, whence 
$\; \dot\t_s = U_s\,(\vec{\b}\wedge\t_s+\vec{\e}_s)\,$, with $\,\vec{\e}_s\bot\t_s\,$ and
$\,\lim_{s\to\ii}\limits\,\vec{\e}_s=0$ ; so that  writing $\;\t_s=
\la_s\,\vec{\b} + \vec{V}_s\,$ with $\,\vec{V}_s\bot\vec{\b}\,$, \  we have as $\,s\to\ii\,$ : \quad  
$|\vec{V}_s| = |\vec{\b}\wedge\vec{V}_s| = |\vec{\b}\wedge\t_s| = |\dot\t_s/U_s-\vec{\e}_s|\lra 1\,$,  whence $\,\la_s\lra 0\,$. \par  \smallskip   
   Moreover \ $\vec{\b}\wedge\vec{V}_s = \vec{\b}\wedge\t_s\,$ and then, setting 
$\;{\ds\vec{v}_s:=\,{\vec{V}_s / |\vec{V}_s|}}\,$, \ we have on one hand :   
$$ |\t_s-v_s| \le |\la_s|+| |\vec{V}_s|-1| \lra 0\;\hbox{ almost surely as }\; s\to\ii\, , $$ 
and on the other hand : 
$$ d\vec{v}_s= d\,{\vec{V}_s\over|\vec{V}_s|} = {d\t_s-d\la_s\,\vec{\b}\over |\vec{V}_s|} -
{\vec{V}_s\over |\vec{V}_s|^{2}}\,d|\vec{V}_s| = 
\vec{\b}\wedge{\vec{V}_s\over|\vec{V}_s|}\,U_s\,ds - {\vec{\b}\over|\vec{V}_s|}\,d\la_s 
+ {U_s\over|\vec{V}_s|}\,\vec{\e}_s\,ds - {\vec{V}_s\over |\vec{V}_s|^{2}}\,d|\vec{V}_s| $$   
$$ = \,\Big(\vec{\b}\wedge{\vec{v}_s}\Big)\,U_s\,ds + (\vec{\e}_s\cdot[\vec{\b}\wedge {\vec{v}_s}])\,
\Big( \vec{\b}\wedge{\vec{v}_s}\Big) \,{U_s\over|\vec{V}_s|}\, ds \qquad \Big(\hbox{since }\; 
d\vec{v}_s\bot\{\vec{\b},\vec{v}_s\}\Big) $$ 
$$ =\, \Big( \vec{\b}\wedge{\vec{v}_s}\Big)\,\Big( 1+ det\Big[{\ts{\vec{\e}_s\over|\vec{V}_s|}}\,,
\,\vec{\b}\,, \,{\vec{v}_s} \Big]\Big)\, U_s\,ds\, . $$

   Let us now denote by $\,\phi_s\,$ the angular coordinate of $\,\vec{v}_s\,$ in the constant plane $\,\vec{\b}^\bot\,$ ; then the preceding equation writes equivalently (using the expression of $\,U_s\,$ seen at the beginning of this proof, and choosing the orientation in the plane 
$\,\vec{\b}^\bot\,$ in order to have a positive sign for the remainder of this proof) :  
$$ d\phi_s\,=\, \Big( 1+ det\Big[{\ts{\vec{\e}_s\over|\vec{V}_s|}}\,,\,\vec{\b}\,, \,{\vec{v}_s}
\Big]\Big)\, U_s\,ds\, = \Big( 1+ det\Big[{\ts{\vec{\e}_s\over|\vec{V}_s|}}\,,\,\vec{\b}\,, \,{\vec{v}_s}
\Big]\Big)\, {{\rm sgn}(T_s)\, dr_s\over\sqrt{r_s\,[R-r_s+(\ell^2+\e_s)r_s^3]}}\, , $$ 
whence (setting $\,\d_s:=det\Big[{\ts{\vec{\e}_s\over|\vec{V}_s|}}\,,\,\vec{\b}\,, \,{\vec{v}_s}\Big]$)
for any $\,n\in\N\,$ and $\,\e > 0 $ : 
$$ \int^{D_{n}^\e}_{D_{3n+1}} d v_s \, =  \int^{D_{n}^\e}_{D_{3n+1}} d\phi_s \, =   \int_{D_{3n+1}}^{D_{n}^\e}  {(1+\d_s)\, dr_s\over\sqrt{r_s\,[R-r_s+(\ell^2+\e_s)\,r_s^3]}}\,, $$  
and idem with $\,D_{n}'\,$ instead of $\,D_{n}^\e\,$. \ Since we have \  $ {\ds  \lim_{n\to\ii}\, \Big| \int^{D_{n}^\e}_{D_{3n+1}} d \t_s - \int^{D_{n}^\e}_{D_{3n+1}} d v_s \Big| \, =0 \, } $ by the beginning of this proof, we have henceforth to deal with \ 
$ {\ds \int_{D_{3n+1}}^{D_{n}^\e}  {(1+\d_s)\, dr_s\over\sqrt{r_s\,[R-r_s+(\ell^2+\e_s)\,r_s^3]}}\,}  $ \parn  (and analogously with $\,D_{n}'\,$ instead of $\,D_{n}^\e\,$). 
Recall that $\;\lim_{n\to\ii}\limits\, r_{D_{n}^\e} =\rr_\e <\rr\;$ by Lemma \ref{lem.piege}. \par 

   Now, pushing somewhat further the observation we already used for proving Theorem \ref{the.martp},
we can use the strictly increasing radius $\,r_s\,$ as alternative ``time'' variable on
the proper time interval $\,[D_{3n+1},D_{n}']\,$. \ So we set
$\,s=s_n(r)\,,\,\tilde{\e}_r^n:=\e_{s_n(r)}\,,\, \tilde{\d}_r^n:=\d_{s_n(r)}\,$ on
$\,[D_{3n+1},D_{n}']\,$, \ and we get almost surely :  
$$ \int^{D_{n}^\e}_{D_{3n+1}} d\phi_s \,  =  \int_{0}^{r_{D_{n}^\e}}  {(1+\tilde{\d}_r^n)\,
dr\over\sqrt{r\,[R-r+(\ell^2+\tilde{\e}_r^n)\,r^3]}}\,. $$  
Note that we have $\,|\tilde{\e}_r^n| = |\e_{s(r)}| \ra 0\,$ and $\,|\tilde{\d}_r^n| = \O(|\vec{\e}_{s(r)}|) \ra 0 $, uniformly on $[D_{3n+1},D_{n}']$ as $\,n\to\ii\,$, and that 
we have for large enough $\,n\,$ :  $\;r_{D_{n}^\e} <(\rr_\e +\rr)/2 < \rr\,$.  \parn 
Hence, using the function $\,g\,$ introduced at the beginning of this proof, we see that \parn 
$ 0\le r\le r_{D_{n}^\e} \,\Ra\, g(r)\ge g((\rr_\e +\rr)/2)  >0\,$, and then that dominated convergence holds, for any $\,\e>0\,$, showing that we have almost surely 
$$ \lim_{n\to\ii} \,\int^{D_{n}^\e}_{D_{3n+1}} d\phi_s \,  =  \int_{0}^{\rr_{\e}}  { dr\over
\sqrt{r\,[R-r+\ell^2\,r^3]}}\,. $$  

  Since it is clear that $\;{\ds \int_{0}^{\rr_{\e}}  { dr\over \sqrt{r\,[R-r+\ell^2\,r^3]}}\;} $  converges to $\;{\ds \int_{0}^{\rr}  { dr\over \sqrt{r\,[R-r+\ell^2\,r^3]}}\;} $ as $\,\e\sea 0\,$, we are left with the remainder $\;{\ds \int_{D_{n}^\e}^{D_{n}'} d\phi \;} $, which we must control uniformly.  Note that such control is not obvious at all, since we have a root of the denominator precisely at time 
$\,D'_n\,$ (at which we are at the top of the excursion), according to the fact that around this same time the angular move is much more rapid than the radial one. 

   Now we have, for large enough $\,n\,$ (using that $\, \e<{\ts{2\over R}} \Ra (1-R/\rr_\e)\,\rr_\e\2 > - ({\ts{R\over 2}})\2 \Ra \rr_\e > {\ts{R\over 2}}  \Ra r_{D_{n}^\e} \ge {\ts{R\over 3}}\, $) : 
$$ \Big|  \int_{D_{n}^\e}^{D_{n}'} d\phi_s \,\Big|    \,=\,   \int_{r_{D_{n}^\e}}^{r_{D_{n}'}} {(1+\tilde{\d}_r^n)\,dr\over
\sqrt{r\,[R-r+(\ell^2+\tilde{\e}_r^n) r^3]}}  \,
 \le\, 2\int_{r_{D_{n}^\e}}^{r_{D_{n}'}}  {dr\over \sqrt{r\,[R-r+(\ell^2+\tilde{\e}_r^n) r^3]}} $$ 
$$ = \, 2\int_{r_{D_{n}^\e}}^{r_{D_{n}'}}  {b_{s_n(r)}\,d(s_n(r))\over r^{2}}\, 
\le\,  {\ts {18\over R^{2}}} \int_{r_{D_{n}^\e}}^{r_{D_{n}'}}  b_{s_n(r)}\,d(s_n(r)) \,   =\, 
 {\ts {18\over R^{2}}} \int_{D^\e_n}^{{D_{n}'}} b_s\,ds \, .  $$ 
 
    Therefore, applying Lemma \ref{lem.intb}, we conclude that we have with probability larger than $1-\e$ :  \parn 
\centerline{ $ {\ds \limsup_{n\to\ii} \Big| \int_{D_{n}^\e}^{D_{n}'} d\phi_s \,\Big| \,\le  \,{\ts{12\over 3R-2\rr}}\,({\ts{3R\over \sqrt{2}}})^4\,\e\;} $, \  and then } 
$$  \limsup_{n\to\ii}  \Big|  \int_{D_{3n+1}}^{D_{n}'} d\phi_s - \int_{0}^{\rr}  {dr\over \sqrt{r\,[R-r+\ell^2\,r^3]}}\,\Big| \,\le \,{\ts{12\over 3R-2\rr}}\, ({\ts{3R\over \sqrt{2}}})^4\,\e +  \int_{\rr_\e}^{\rr}  {dr\over \sqrt{r\,[R-r+\ell^2\,r^3]}} \, , $$ 
which goes to 0 as $\,\e\sea 0\,$, while the left hand side does not depend upon $\,\e\,$. This proves that indeed we have almost surely : 
$$  \lim_{n\to\ii} \, \int_{D_{3n+1}}^{D_{n}'} d\t_s \, = \, \lim_{n\to\ii} \, \int_{D_{3n+1}}^{D_{n}'} d\phi_s \, = \, \int_{0}^{\rr}  {dr\over \sqrt{r\,[R-r+\ell^2 \,r^3]}}\,. \;\;\diamond $$ 

 \section{Appendix : \ Study of timelike and null geodesics}\label{sec.geod}
 \subsection{Timelike geodesics} \label{Tlgeod}  \indf     
    The case $\,\s =0\,$ in Section \ref{sec.S} corresponds to geodesics ; precisely it is the case of timelike geodesics having speed 1. The equations we obtained in Section \ref{sec.S} are here  simplified as follows : 
$$ a \; ,\, b \; \hbox{ constant }\; ; \qquad \ddot r_s = (r_s- {\ts{3\over 2}}R)\,{b^2\over r_s^4} -
{R\over 2\, r_s^2}\; . $$ 
Integrating this last equation (after multiplication by $\,\dot r_s$) leads simply, up to some constant, 
to the unit pseudo-norm relation $\;^t\! e_0\, g\, e_0=1\,$ : 
$$ \dot r_s^2 = a^2-(1-{\ts{R\over r_s}})(1+{\ts{b^2\over r^2_s}})\; , $$ 
or equivalently to : $\qquad {\ds | s | = \int^{r(s)}_{r(0)}{r\, dr\over\sqrt{a^2r^2-(1-{R\over
r})(r^2+b^2)}}}\;$.
\par\smallskip\par  

   The equations relative to $\,\f_s\,$ and to $\,\psi_s\,$ yield easily a real constant $\,k\,$ such that : 
$$ \dot \f_s^2 = {b^2\over r_s^4} - {k^2\over r_s^4\,\sin^2\f_s}\quad ,\qquad \dot \psi_s = {k\over
r_s^2\,\sin^2\f_s}\; . $$
They are equivalent to the equation $\quad{{d\over ds}\,(r^2_s\,\dot \t_s)^{} = -\,b^2\,r\2_s\;\t_s}\,$. \quad (This is also a consequence of Proposition \ref{pro.systeqst} in Section \ref{sec.bhh} below, taking there $\,\s =0\,$.) \ They are also equivalent to the constancy of the angular momentum $\,\bv\,$. This means in particular that $\;{\ds \ddot\t = -U^2\,\t -2(\dot r /r)\,\dot\t}\;$ along geodesics, and therefore that every geodesic is included in some plane containing the origin $\{ r=0\}$. \par\medskip  

   Set $\;u:= 1/r\;$ and $\; P(u) := (1-Ru)(1+b^2u^2)\,$. \ The variation of $P$ as $\,r\,$ increases
from $R$ to $\ii$ is as follows : \parn 
- \ if $\,b\le R\sqrt{3}\,$ : \ $ a^2-P(u)\,$ decreases (as $r$ increases from $R$ to $\ii$) 
from $a^2$ to $a^2-1\,$ ; \parn 
- \ if $\,b > R\sqrt{3}\,$ : (as $r$ increases from $R$ to $\ii$) \ $ a^2-P(u)\,$ decreases 
firstly from $a^2$ to $a^2-P(u_1)\,$, then increases from $a^2-P(u_1)\,$ to $a^2-P(u_2)\,$, then decreases from $a^2-P(u_2)\,$\parn  to $\,a^2-1\,$, \ where \quad 
$ {2\over 3R}\, >\, u_1:= {1+\sqrt{1-3R^2b\2}\over3\, R}\, >\, u_2:= {1-\sqrt{1-3R^2b\2}\over3\, R}
\, >\, 0\,$.\par  
   Note that for $\,b > R\sqrt{3}\,$
$\; P(u_1) = {8\over 9}+ {2\over 27R^2}(b^2-3R^2)(1+\sqrt{1-3R^2b\2}\,)\, \in\, ]{8\over 9},+\ii\, [\,$,
\parn   and  
$\; P(u_2) = {8\over 9}+ {2\over 27R^2}(b^2-3R^2)(1-\sqrt{1-3R^2b\2}\,)\, \in\, ]{8\over 9},\min\{
1,P(u_1)\} [\,$. \par 

   Besides, $\;{r\over\sqrt{a^2r^2-(1-{R\over r})(r^2+b^2)}}\;$ is integrable near a simple root of 
$\,a^2-P(1/r)\,$, but not near a double root. \par 
   Hence, using the relation $\;\dot r^2 = a^2-P(u)\,$ (which does not allow any value for $\,r_0$), we get the following classification of timelike geodesics. (See [L-L], \$ 100, problems 1,2, for a partial
resolution, and [M-T-W], section 25.5, for a more explicit one). \par 
   Note that we restrict first to the case of $\,{\cal S}_0\,$. The extension to the full Schwarzschild space $\,{\cal S}\,$ is then easy : see Remark \ref{rem.fin} below.  \par \smallskip 

\noindent {\bf 1.} \ $P\,$ monotone \par 
   \ub{Case 1.1} : $\,b\le R\sqrt{3}\,$ and $\,|a|\ge 1 $ : $\,a^2-P\,$ has no root, so that
$\;(s\mapsto r_s)\,$ runs, in one direction or in the other one, an increasing trajectory from $R$ to
$+\ii\,$, slowing down (in the increasing case)  till the limit speed $\,\sqrt{a^2-1}\,$. In the
particular case $\,|a|=1\,$, we have $\,r_s\sim (9Rs^2/4)^{1/3}\;$ as $s\to\ii\,$. \par 
   \ub{Case 1.2} : $\,b\le R\sqrt{3}\,$ and $\,a^2 <1$ (but $\,b\not= R\sqrt{3}$ or $a^2\not= 8/9$) : 
$\,a^2-P\,$ has a simple root, so that $\;(s\mapsto r_s)\;$ runs a bounded trajectory, which has 2 ends at $R$, and increases firstly then decreases, with a unique maximum  at $R_0\in ]R,\ii [\,$. 
$\,r_0>R_0\,$ is impossible. \par 
   \ub{Case 1.3} : $\,b = R\sqrt{3}\,$ and $\,a^2=P(u_1)=8/9$ : $\,a^2-P\,$ has a triple root at
$\,{1\over u_1} =3R\,$, so that we have here either a geodesic which runs monotonically (during an
infinite time) the interval $\,[R,3R[\,$, in one direction or in the other one, or a geodesic included
in a circle (centred at 0) of radius $3R$. Indeed, it is easily verified (looking at $\ddot r$) that such
circular geodesics  correspond precisely to multiple roots of $\,a^2-P\,$. 
$\,r_0>3R\,$ is impossible. \par 

\noindent {\bf 2.} \ $P\,$ non-monotone \par 
   \ub{Case 2.1} : $\,b > R\sqrt{3}\,$ and $\,a^2\ge\max\{ P(u_1),1\}\,$, $\,a^2>P(u_1)$ : $\,a^2-P\,$
has no root, so that $\,\dot r\,$ does not vanish and we are brought back to the case 1.1, except that there is an acceleration phase on the  interval $r\in\, ]{1\over u_1},{1\over u_2}[\,$. \par 
   \ub{Case 2.2.1} : $\,b > R\sqrt{3}\,$ and $\,a^2< P(u_2) $ : $\,a^2-P\,$ has a unique (simple) root,
$\,\dot r\,$  vanishes at $r=R_0\in ]R,{1\over u_1}[\,$, \ and we are brought back to the case 1.2. \par 
   \ub{Case 2.2.2} : $\,b > R\sqrt{3}\,$ and $\,P(u_1)<a^2<1 $ : $\,a^2-P\,$ has a unique (simple) root, $\,\dot r\,$ vanishes at $r=R_2\in ]{1\over u_2},\ii[\,$, \ and we are brought back to the case 1.2, alternatively with acceleration and slackening phases, the unique maximum being here at $R_2\,$. \par 
   \ub{Case 2.3} : $\,b > R\sqrt{3}\,$ and $\,a^2=P(u_1)\ge 1$ : $\,a^2-P\,$ has a double root at
$\{ r={1\over u_1}\}$ (and no other root), so that the trajectory $\,(s\mapsto r_s)\,$ needs an infinite
time to reach this  level. Such geodesics run monotonically, in one direction or in the other one, either the interval $\,[R,{1\over u_1}[\,$, or the interval $\,]{1\over u_1},\ii ]\,$. There are again
here also geodesics included in a circle (centred at 0) of radius ${1\over u_1}$. \par 
   \ub{Case 2.4} : $\,b > R\sqrt{3}\,$ and $\,1\le a^2<P(u_1)$ : $\,a^2-P\,$ has a two (simple)
roots, $\,\dot r\,$ vanishes at $r=R_0\in ]R,{1\over u_1}[\,$ and at $r=R_1\in ]{1\over u_1},{1\over u_2}[$ ; \ we are brought back to the case 1.2 if $\,r_0\in [R,R_0]$ ; $\,r_0\in ]R_0,R_1[$ is impossible ; and if $\,r_0\in [R_1,\ii]$, then $\;(s\mapsto r_s)\;$ runs an unbounded 
trajectory the 2 ends of which are at $\ii$, with a unique minimum at $R_1\,$. In this last subcase we
have $\,r_s\sim |s|\,\sqrt{a^2-1}\;$  (and even $\,r_s= |s|\sqrt{a^2-1}+\log (1+|s|)+r_0+ o(1)\,$) as
$s\to\ii\,$ if $\,a^2>1$, and $\,r_s\sim (9Rs^2/4)^{1/3}\,$ if $\,a^2=1$. Such (projection of) geodesic
runs  approximately a parabola or a branch of hyperbola.  \par 

   \ub{Case 2.5.1} : $\,b > R\sqrt{3}\,$ and $\,a^2=P(u_1)< 1$ : $\,a^2-P\,$ has a double root at $\{
r={1\over u_1}\}$ and a simple root at $\,r=R_2>{1\over u_2}\,$. Hence the trajectory $\,(s\mapsto r_s)\,$
needs an infinite time to reach the level $\{ r={1\over u_1}\}$, so that there are on one hand
(projection of) geodesics which run monotonically the interval $\,[R,{1\over u_1}[\,$, in one direction
or in the other one, and on the other hand (projection of) geodesics which run (during an infinite time)
from $\,{1\over u_1}\,$ to $\,{1\over u_1}\,$ via a unique maximum at $R_2$. Again, there are also here geodesics included in a circle (centred at 0) of radius ${1\over u_1}$. $\,r_0>R_2\,$ is impossible. \par 
   \ub{Case 2.5.2} : $\,b > R\sqrt{3}\,$ and $\,a^2=P(u_2)$ : $\,a^2-P\,$ has a double root at $\{
r={1\over u_2}\}$ and a simple root at $\,r=R_0<{1\over u_1}\,$. Hence we are brought back to the
case 1.2 if $\,r_0\in [R,R_0]$, and we have also centred circular geodesics of radius ${1\over u_2}$.
$r_0\in ]R_0,{1\over u_2}[\cup ]{1\over u_2},\ii[$ is impossible. \par 

   \ub{Case 2.6} : $\,b > R\sqrt{3}\,$ and $\,P(u_2)< a^2< \min\{ P(u_1),1\} $ : 
$\,a^2-P\,$ has 3 distinct simple roots, so that $\,\dot r\,$ vanishes at $r=R_0\in ]R,{1\over u_1}[\,$, at $r=R_1\in ]{1\over u_1},{1\over u_2}]\,$ and at $r=R_2\in [{1\over u_2},\ii[$ ; we are brought back to the case 1.2 if \hbox{$\,r_0\in [R,R_0]$ ;} if $\,r_0\in [R_1,R_2]\,$, then $(s\mapsto r_s)$ oscillates periodically, increasing from $R_1$ to $R_2$ then decreasing from $R_2$ to $R_1$, running approximately an ellipse. $\,r_0\in]R_0,R_1]\cup ]R_2,\ii[\,$ is impossible. \par 
\medskip 
   Note that the set of limiting radii $\,{1\over u_1}\,$ such that there exists a geodesic which winds
asymptotically (either from inside or from outside) around a circle (centred at 0) of radius $\,{1\over
u_1}\,$ equals the interval $\,]3R/2,3R]\,$. Note also that the set of radii ${1\over u_1}$ or
${1\over u_2}$ of circles (centred at 0) which contain geodesics equals the interval $\,]3R/2,\ii[\,$.
\par 
This gives a geometrical intrinsic meaning to the particular radii $\,3R$ and $3R/2\,$. 

\begin{rem}\label{rem.fin}\quad {\rm The extension to the full space $\,{\cal S}\,$ of the preceding classification and description of timelike geodesics in the restricted space $\,{\cal S}_0\,$ is more or less straightforward. Indeed the results of Theorem \ref{the.bh} are clearly valid for geodesics as well, with the major simplification that $\,a_s$ and $\bv_s\,$ are constant when $\s=0$. \par  
Otherwise, there exist  timelike geodesics included in the cylinder $\{ r=R\}$ : looking at the above  and taking $\,a=0<|k|<b\,$, we find easily such solutions, which satisfy
$\;{\ds \f_s= {\rm Arccos}\Big(\sqrt{1-k^2/b^2}\,\sin[\pm b(s-s_0)/R^2]\Big)\,}$. \ This completes for the space $\cal{S}$ the picture of timelike geodesics we have just drawn above for the strict Schwarzschild space ${\cal S}_0\,$. \par
   Note finally the following intrinsic characterisation of the radius $R$ in the space $\cal{S}$ : \ 
$R\,$ is the minimal radius which can be reached by a timelike geodesic which does not hit the
singularity.  }\end{rem} 

\subsection{Null geodesics} \label{Ngeod}  \indf      
      For these null geodesics, or light rays, the unit pseudo-norm relation is replaced by 
$$  \a^2-\dot r^2 = (1-R/r)\,r\2 \, .   $$ 
Note that the proper time does not make sense any more, so that the new ``time''-parameter or abscissa $\,\la\,$ makes sense only up to an affine transform $\,(\la\mapsto q\la+q')\,$, and the constant parameters $\,a\,$ and $\,b\,$ do not make sense both anymore, but only their quotient $\,\a:=a/b\,$. This unique ``impact parameter''  $\,\a\,$ of the null geodesic is of course a constant of the geodesic. As for the timelike geodesics, the equations relative to $\,\t=(\f,\psi)\,$ show that they are planar, so that by means of a trivial change of axis we may consider that $\,\f\equiv \pi/2\,$. \par 

   Thus every null geodesic is determined by the equations (in its own plane, the derivatives being relative to the abcissa $\,\la\,$) (see also [M-T-W], page 674) : 
$$ \a \;\, {\rm constant}\; ; \quad \dot r^2  + (1-R/r)\,r\2 = \a^2 \; ;\quad \dot t  = \a/(1-R/r)\; ;\quad \dot \psi = r\2\, . $$ 
Eliminating the abcissa $\,\la\,$, we get : 
$$ d\psi = \pm\, r\,[ \a^2-(1-R/r)r\2]^{-1/2}\, dr\; ;\quad dt = -\a\,r\, (1-R/r)\1\,[ \a^2-(1-R/r)r\2]^{-1/2}\, dr\, . $$

   As $\,r\,$ increases from 0 to $\,3R/2\,$,  $\,[ \a^2-(1-R/r)r\2]\,$ decreases from  infinity to $\,\a^2-{4\over 27R^2}\,$, and then increases from $\,\a^2-{4\over 27R^2}\,$ to $\,\a^2\,$ as $\,r\,$ increases from $\,3R/2\,$ to infinity. Therefore, owing to their projection on the coordinate $\,r\,$,  we find three cases for the null geodesics :\par \smallskip  
   - \ub{Case 0} : \  $\;|\a| = {2\over 3\sqrt{3}\,R}\,$ : \quad The null geodesic can be either included 
in $\{ r=3R/2\}$, or it can be asymptotic to $\{ r=3R/2\}$, either growing strictly from $\,r=0\,$ to $\, r=3R/2$,  or growing strictly from $\, r=3R/2\,$ to infinity (or in the reverse direction). \par 
   - \ub{Case 1} : \  $\;|\a| > {2\over 3\sqrt{3}\,R}\,$ : \quad The null geodesic runs from 0 to infinity (with an asymptotic velocity : $\,r_\la\sim\a\la$), or in the reverse direction.  \par 
   - \ub{Case 2} : \  $\;|\a| < {2\over 3\sqrt{3}\,R}\,$ : \quad The null geodesic can either be reminiscent from a parabola, with a minimal radius $\,\rr'>3R/2\,$, or indefinitely oscillate between 
   $\,r=0\,$ and a maximal radius $\,\rr \in\, ]R,3R/2[\,$. \par\smallskip 
   
   This last sort of null geodesics, which are recurrent at the singularity $\,r=0\,$, is the most interesting for us here, since the confined diffusion trajectories of Section \ref{sec.NumbC} seem to have their shape asymptotic to the shape of one such geodesic. Indeed, considering the impact parameter $\,\a\,$ such that $\,|\a| < {2\over 3\sqrt{3}\,R}\,$, the maximal radius $\,\rr <3R/2\,$ solves $\,\a^2\rr^2+R/\rr =1\,$, so that $\,\a\,$ stands for $\,\ell = \lim_sa_s/b_s\,$ in Theorem (\ref{the.piege}, 2) of Section \ref{sec.NumbC}, and the angular deviation $\,\Psi\,$ during each increase from $\,r=0\,$ to $\,r=\rr\,$ (or decrease from $\,r=\rr\,$ to $\,r=0\,$) has exactly the expression found for confined diffusion paths in Corollary \ref{cor.piege} (and Remark \ref{rem.desc}) of Section \ref{sec.NumbC}. \par 
   Moreover, \quad ${\ds t= -\a\int _\rr^r{ r \, dr \over (1-R/r)\sqrt{\a^2r^2-1+R/r}}}\;$ \ shows that each confined null geodesic has its graph invariant under the symmetry $\,(r,t)\mapsto (r,-t)\,$. Indeed this is clear in $\,{\cal S}_0\,$, and remains true everywhere by analytic continuation. 
   In the Kruskal-Szekeres coordinates, this means invariance under the symmetry $\,(u,v)\mapsto (u,-v)\,$. The gluing we defined at the singularity $\,r=0\,$ implies at once that null geodesics must also be symmetric with respect to $\{ u=v=0\}$. \par 
   As a conclusion, we see that the confined null geodesics run indefinitely a closed analytic curve, which in the coordinate plane $\,(r,\psi)\,$ appears as a figure eight (symmetrical and centred at the origin), and in the Kruskal-Szekeres coordinate plane $\,(u,v)\,$ looks like a pair of round brackets (symmetrical with respect to the coordinate axes and joining the two branches $\{ r=0\}$). 

\section{References} \label{Ref} 
\bigskip \medskip 

\vbox{ \noindent 
{\bf [Bi]} \ Bismut J.-M.  \ \  {\it M\'ecanique al\'eatoire. } \par \smallskip \hskip 32mm    
Lecture Notes in Mathematics n$^o$ 866, Springer, Berlin 1981. }
\bigskip 

\vbox{ \noindent 
{\bf [D]} \ Debbasch F.  \ \  {\it A diffusion process in curved space-time.}  
\par \hskip 25mm \quad J. Math. Phys. 45, n$^o$ 7, 2744-2760, 2004. }
\bigskip 

\vbox{ \noindent 
{\bf [DF-C]} \ De Felice F. , Clarke C.J.S. \ \  {\it Relativity on curved manifolds. }
\par
\smallskip \hskip 11mm    Cambridge surveys on mathematical physics, Cambridge university press, 1990. }
\bigskip 

\vbox{ \noindent 
{\bf [D1]} \ Dudley R.M. \ \ {\it Lorentz-invariant Markov processes in relativistic phase space.}
\par \smallskip \hskip 33mm  Arkiv f\"or Matematik 6, n$^o$ 14, 241-268, 1965. }
\bigskip 

\vbox{ \noindent 
{\bf [D2]} \ Dudley R.M. \ \ {\it A note on Lorentz-invariant Markov processes.}
\par \smallskip \hskip 33mm  Arkiv f\"or Matematik 6, n$^o$ 30, 575-581, 1967. }
\bigskip 

\vbox{ \noindent 
{\bf [D3]} \ Dudley R.M. \ \ {\it Asymptotics of some relativistic Markov processes.}
\par \smallskip \hskip 33mm  Proc. Nat. Acad. Sci. USA n$^o$ 70, 3551-3555, 1973. }
\bigskip 

\vbox{ \noindent 
{\bf [D4]} \ Dudley R.M. \ \ {\it Recession of some relativistic Markov processes.}
\par \smallskip \hskip 33mm  Rocky Mountain J. Math. n$^o$ 4, 401-406, 1974. }
\bigskip 

\vbox{ \noindent 
{\bf [E]} \ Elworthy D.  \ \  {\it Geometric aspects of diffusions on manifolds. }
\par \smallskip 
\hskip 29mm  \'Ecole d'\'EtŽ de ProbabilitŽs de Saint-Flour XV-XVII, 1985-87,  277-425, \par \hskip 29mm  Lecture Notes in Math. 1362, Springer, Berlin 1988.}
\bigskip 

\vbox{ \noindent 
{\bf [Em]} \ \'Emery M.  \ \  {\it On two transfer principles in stochastic differential geometry. }
\par \smallskip 
\hskip 29mm   S\'eminaire de Probabilit\'es XXIV, 407-441, Springer, Berlin 1990. }
\bigskip 

\vbox{ \noindent 
{\bf [E-F-LJ]} \  Enriquez N. , Franchi J. , Le Jan Y.  \  \ {\it Canonical lift and exit law of the \par 
\hskip 43mm  fundamental diffusion associated with a Kleinian group.}  
\par \hskip 43mm   S\'eminaire de Probabilit\'es XXXV, 206-219, Springer 2001.} 
\bigskip 

\vbox{ \noindent 
{\bf [F-N]} \ Foster J. , Nightingale J.D. \ \  {\it A short course in General Relativity. }
\par
\smallskip \hskip 61mm    Longman, London 1979. }
\bigskip 

\vbox{ \noindent 
{\bf [H]} \ Hsu E.P.  \ \  {\it Stochastic analysis on manifolds. } \par \smallskip \hskip 23mm    
Graduate studies in Mathematics vol. 38, A.M.S., Providence 2002. }
\bigskip 

\vbox{  \noindent 
{\bf [I-W]} \ Ikeda N. ,  Watanabe S. \quad {\it Stochastic differential equations and diffusion
processes.}
\par \hskip 56mm North-Holland Kodansha, 1981. }
\medskip 

\vbox{ \noindent 
{\bf [K]} \ Kochubei A.N. \quad  {\it On diffusions with invariant generating operators.} \par
\smallskip \hskip 21mm  Theory Prob. Appl. 34, 597-603, 1989. (Translated from Russian Journal.) }
\bigskip 

\vbox{  \noindent 
{\bf [L-L]} \ Landau L. ,  Lifchitz E. \quad {\it Physique thŽorique, tome II : ThŽorie des champs.}
\par \hskip 54mm \'Editions MIR de Moscou, 1970. }
\medskip 

\vbox{ \noindent 
{\bf [M]} \ Malliavin P. \ \  {\it GŽomŽtrie diffŽrentielle stochastique. } \par
\smallskip \hskip 30mm  Les Presses de l'UniversitŽ de MontrŽal, MontrŽal 1978. }
\bigskip 

\vbox{ \noindent 
{\bf [Ms]} \ Markus L. \ \  {\it Global Lorentz geometry and relativistic Brownian motion. } \par
\smallskip \hskip 23mm  in \ Pitman research notes, Math. series 150, 273-286, 1986. }
\bigskip 

\vbox{ \noindent 
{\bf [M-T-W]} \ Misner C.W , Thorne K.S. , Wheeler J.A. \qquad    {\it Gravitation.} \par
\smallskip \hskip 15mm   W.H. Freeman and Company, New York, 1973. }
\bigskip 

\vbox{ \noindent 
{\bf [S]} \ Stephani H. \ \ {\it General Relativity. } \qquad Cambridge university press, 1990. }
\bigskip 

\centerline{--------------------------------------------------------------------------------------------}
\medskip

  \ub{A.M.S. Classification} : \ Primary 58J65, secondary 53C50, 60J65, 83C57, 83A05, 83C10. \par 

  \ub{Key Words} : \ Relativistic diffusion, Brownian motion, Stochastic flow, General relativity, 
Lorentz manifold, Schwarzschild space, Black hole, Asymptotic behaviour. 

\medskip 

\centerline{--------------------------------------------------------------------------------------------}

\bigskip 

\vbox{  
Jacques FRANCHI : \quad Universit\'e Louis Pasteur, I.R.M.A., 7 rue Ren\'e Descartes, \\ 
67084 Strasbourg cedex. France. \quad franchi@math.u-strasbg.fr 

\smallskip \noindent 

Yves LE JAN : \quad Universit\'e Paris Sud, Math\'ematiques, B\^atiment 425, 91405 Orsay. France. 
\quad  yves.lejan@math.u-psud.fr
}
\medskip 

\centerline{--------------------------------------------------------------------------------------------}

\end{document}